\documentclass[11pt]{article}
\usepackage[utf8]{inputenc}
\usepackage{float} 
\usepackage{amsmath}
\usepackage{amssymb}

\usepackage{array}
\usepackage{caption}
\usepackage{subcaption}
\usepackage{bbm}
\usepackage{makeidx}
\usepackage[linktocpage=true]{hyperref}
\usepackage[capitalise]{cleveref}
\usepackage{amsthm}
\usepackage{color}
\usepackage{mathtools}
\usepackage{graphicx} 
\usepackage[LGR,T1]{fontenc}
\usepackage{siunitx}
\usepackage{mathrsfs}
\usepackage{scalerel}
\usepackage{forest}
\usepackage{algorithm,algpseudocode}
\usepackage{lipsum}
\usepackage{xcolor}
\usepackage[title]{appendix}
\hypersetup{
    colorlinks,
    linkcolor={blue!80!black},
    citecolor={green!50!black},
}
\usepackage{accents}
\usepackage{apptools}
\usepackage{tikz}
\usetikzlibrary{arrows, positioning, automata}
\usepackage{comment}
\usepackage[shortlabels]{enumitem}

\usepackage[mathscr]{euscript}
\DeclareSymbolFont{rsfs}{U}{rsfs}{m}{n}
\DeclareSymbolFontAlphabet{\mathscrsfs}{rsfs}

\AtAppendix{\counterwithin{lemma}{section}}

\usepackage{mathtools}

\theoremstyle{plain}

\newtheorem{theorem}{Theorem}
\newtheorem{proposition}{Proposition}[section]

\newtheorem{lemma}{Lemma}[section]
\newtheorem{corollary}{Corollary}[section]
\newtheorem{assumption}{Assumption}[section]

\newtheoremstyle{myremark} 
    {\topsep}                    
    {\topsep}                    
    {\rm}                        
    {}                           
    {\bf}                        
    {.}                          
    {.5em}                       
    {}  

\theoremstyle{myremark}
\newtheorem{remark}{Remark}[section]

\allowdisplaybreaks
\usepackage[top=1in,bottom=1in,left=1in,right=1in]{geometry}
\usepackage[utf8]{inputenc}

\def\bfone{{\boldsymbol 1}}
\def\bbone{\mathbbm{1}}
\def\cuF{\mathscrsfs{F}}
\def\cuP{\mathscrsfs{P}}
\def\cuE{\mathscrsfs{E}}
\def\cuD{\mathscrsfs{D}}

\def\mmse{{\sf mmse}}
\def\TAP{\mathrm{TAP}}
\def\GD{\mathsf{GD}}
\def\data{\bD}

\def\<{\langle}
\def\>{\rangle}
\def\dd{\mathrm{d}}

\def\top{\intercal}

\DeclareMathOperator*{\plim}{p-lim}
\DeclareMathOperator*{\plimsup}{p-lim\,sup}
\def\AMP{\mathsf{AMP}}
\def\op{\mbox{\rm \tiny op}}

\def\htheta{\hat{\Theta}}

\def\hm{\hat{m}}

\def\cC{\mathcal{C}}

\def\cA{\mathcal{A}}

\def\GOE{\mbox{\rm GOE}}

\def\mmse{\mathsf{mmse}}

\def\bhY{\hat{\bY}}
\def\obY{\overline{\bY}}
\def\obG{\overline{\bG}}

\def\oG{\overline{G}}

\def\Lamp[#1]{\boldsymbol{\Lambda}_{\mathrm{AMP}}^{(#1)}}
\def\lalg[#1]{\Lambda_{\mathrm{alg}, #1}}

\def\ker{{\mathrm{ker}}}
\def\proj{{\boldsymbol{P}}}

\def\sB{\mbox{\rm \tiny B}}
\def\sTAP{\mbox{\rm \tiny TAP}}
\def\sNMF{\mbox{\rm \tiny NMF}}

\def\ep{\eps}

\def\odelta{\overline{\delta}}
\def\Info{{\sf I}}
\def\Sym{{\mathbb S}}
\def\Spher{{\mathsf S}}
\def\Orth{{\mathbb O}}

\def\supp{{\rm supp}}
\def\ed{\stackrel{{\textrm d}}{=}}

\def\rP{{\rm P}}
\def\de{{\rm d}}
\def\prob{{\mathbb P}}
\def\sF{{\sf F}}
\def\integers{{\mathbb Z}}

\def\EE{\mathbb{E}}

\def\NN{\mathbb{N}}

\def\PP{\mathbb{P}}

\def\RR{\mathbb{R}}

\def\ZZ{\mathbb{Z}}

\def\cB{{\mathcal B}}

\def\bA{\mathbf{A}}
\def\bB{\mathbf{B}}

\def\bD{\mathbf{D}}

\def\bF{\mathbf{F}}
\def\bG{\mathbf{G}}
\def\bH{\mathbf{H}}
\def\bI{\mathbf{I}}

\def\bM{\mathbf{M}}

\def\bP{\mathbf{P}}
\def\bQ{\mathbf{Q}}

\def\bT{\mathbf{T}}

\def\bW{\mathbf{W}}
\def\bX{\mathbf{X}}
\def\bY{\mathbf{Y}}
\def\bZ{\mathbf{Z}}

\def\ba{\boldsymbol{a}}
\def\bb{\boldsymbol{b}}

\def\bg{\boldsymbol{g}}

\def\bm{\boldsymbol{m}}

\def\bp{\boldsymbol{p}}

\def\bs{\boldsymbol{s}}
\def\bt{\boldsymbol{t}}
\def\bu{\boldsymbol{u}}
\def\bv{\boldsymbol{v}}
\def\bw{\boldsymbol{w}}
\def\bx{\boldsymbol{x}}
\def\by{\boldsymbol{y}}
\def\bz{\boldsymbol{z}}

\def\bA{\boldsymbol{A}}
\def\bB{\boldsymbol{B}}

\def\bD{\boldsymbol{D}}

\def\bF{\boldsymbol{F}}
\def\bG{\boldsymbol{G}}
\def\bH{\boldsymbol{H}}
\def\bI{\boldsymbol{I}}

\def\bM{\boldsymbol{M}}

\def\bP{\boldsymbol{P}}
\def\bQ{\boldsymbol{Q}}

\def\bT{\boldsymbol{T}}

\def\bW{\boldsymbol{W}}
\def\bX{\boldsymbol{X}}
\def\bY{\boldsymbol{Y}}
\def\bZ{\boldsymbol{Z}}

\def\normal{{\mathsf{N}}}
\def\Ball{{\sf B}}

\def\btheta{\boldsymbol{\theta}}

\def\blambda{\boldsymbol{\lambda}}

\def\bgamma{\boldsymbol{\gamma}}
\def\bOmega{\boldsymbol{\Omega}}

\def\beps{\boldsymbol{\varepsilon}}
\def\bepsilon{\boldsymbol{\varepsilon}}
\def\bzero{\boldsymbol{0}}

\def\bzeta{\boldsymbol{\zeta}}

\def\bnu{\boldsymbol{\nu}}
\def\bby{\bar{\by}}

\def\bphi{{\boldsymbol \varphi}}

\def\hbm{\hat{\boldsymbol m}}
\def\hbp{\hat{\boldsymbol p}}

\def\obp{\overline{\boldsymbol p}}
\def\hbM{\hat{\boldsymbol M}}
\def\hby{\hat{\boldsymbol y}}
\def\tbX{\widetilde{\boldsymbol X}}
\def\id{{\boldsymbol I}}
\def\bfzero{{\boldsymbol 0}}
\def\cL{{\mathcal L}}

\def\oR{{\overline{R}}}
\def\sP{{\sf P}}
\def\sM{{\sf M}} 

\def\Lipc{{\rm Lip}_0}
\def\Lips{{\rm Lip}_*}
\def\reals{{\mathbb R}}
\def\sT{{\sf T}}
\def\sconv{\mbox{\rm\tiny conv}}
\def\sreg{\mbox{\rm\tiny reg}}
\def\salg{\mbox{\rm\tiny alg}}
\def\swr{\mbox{\rm\tiny wr}}

\def\sAMP{\mbox{\rm\tiny AMP}}
\def\sNGD{\mbox{\rm\tiny NGD}}
\def\sKL{\mbox{\rm\tiny KL}}
\def\sIT{\mbox{\rm\tiny IT}}
\def\Law{{\rm Law}}
\def\Ball{{\sf B}}
\def\tcA{\tilde{\mathcal A}}

\renewcommand{\P}{\mathbb{P}}
\newcommand{\E}{\mathbb{E}}

\newcommand{\eps}{\varepsilon}
\newcommand{\Var}{\operatorname{Var}}

\newcommand{\argmin}{\operatorname{argmin}}

\newcommand{\toP}{\overset{\P}{\to}}
\newcommand{\Cov}{\operatorname{Cov}}
\newcommand{\sign}{\operatorname{sign}}

\newcommand{\diag}{\operatorname{diag}}

\newcommand{\Unif}{\operatorname{Unif}}

\newcommand{\Tr}{\operatorname{Tr}}

\newcommand{\RN}[1]{%
  \textup{\uppercase\expandafter{\romannumeral#1}}%
}

\newcommand\iidsim{\stackrel{\mathclap{iid}}{\sim}}

\newcommand{\RNum}[1]{\uppercase\expandafter{\romannumeral #1\relax}}



\makeatletter
\newcommand*{\rom}[1]{\expandafter\@slowromancap\romannumeral #1@}
\makeatother

\title{Posterior Sampling in High Dimension via\\ 
Diffusion Processes}
\author{
	Andrea Montanari\footnotemark[2]
	\thanks{Department of Electrical Engineering, Stanford University} 
	\and 
	Yuchen Wu\thanks{Department of Statistics, Stanford University}
}
\date{}
\begin{document}
\maketitle

\begin{abstract}
Sampling from the posterior is a key technical problem in Bayesian statistics. 
Rigorous guarantees are difficult to obtain for 
Markov Chain Monte Carlo algorithms of common use. 
In this paper, we study an alternative class of algorithms based on diffusion processes
and variational methods.
The diffusion is constructed in such a way that, at its final time, it approximates the 
target posterior distribution. The drift of this diffusion is given by
the posterior expectation of the unknown parameter vector
$\btheta$ given the data and the additional noisy observations. 

In order to construct an efficient sampling algorithm, we use
a simple Euler discretization of the diffusion process, and leverage message
passing algorithms and variational inference techniques to approximate the posterior 
expectation oracle. 

We apply this method to posterior sampling in two canonical problems
in high-dimensional statistics: sparse regression and low-rank matrix estimation within the spiked model. 
In both cases we develop the first algorithms with accuracy guarantees in the
regime of constant signal-to-noise ratios.
\end{abstract}

\tableofcontents

\section{Introduction}\label{sec:intro}

Consider the standard setup of Bayesian inference,
whereby the joint distribution of observed data $\data$ and unobserved parameter $\btheta$ 
is defined by
\begin{align}
	\btheta \sim \pi(\,\cdot\,), \qquad \data \sim \rP(\,\cdot\, \mid \btheta),\label{eq:FirstBayes}
\end{align}  
where $\pi(\,\cdot\, )$ is the \emph{prior distribution}.
Bayesian techniques draw inferences from  the
\emph{posterior distribution}\footnote{In the application of Bayes formula below, we are identifying 
 $\rP(\,\cdot\, \mid \btheta)$ with its density with respect to a reference measure.}
\begin{align}
\mu_{\data}(\de\btheta) &:= \prob(\de\btheta |\data)\;
\propto\;  \rP(\data \mid \btheta) \pi(\de\btheta )\, .\label{eq:FirstPosterior}
\end{align}
A substantial amount of research has been devoted to developing
approximation methods \cite{csillery2010approximate,blei2017variational} and sampling algorithms
for the Bayes posterior. Among these, Markov Chain Monte Carlo (MCMC)  
  \cite{gilks1995markov,gamerman2006markov,stuart2010inverse,kroese2013handbook} methods
  play a special role because of their versatility. However,
    MCMC often suffers from  slow mixing 
     \cite{mossel2006limitations,woodard2013convergence}, especially when the posterior
     is multi-modal.  In these circumstances, it might be impossible 
     to run the chain long enough to produce a sample with approximately correct distribution,
     and this can lead to erroneous inference. 
Rigorous upper bounds for the mixing times have been established under various settings
 \cite{levin2017markov,yang2016computational,dalalyan2017theoretical,dwivedi2018log,durmus2019high},
 but proving such bounds is very challenging and existing guarantees only cover a small
 fraction of practical applications.

In this paper we study a different approach to posterior sampling,
which relies on two ideas: 
\begin{itemize}
\item[$(i)$] Construct a non-homogeneous stochastic process
(more precisely, a diffusion process) such that the marginal distribution of which at time
$t$ converges, as $t\to\infty$, to $\mu_{\data}(\de\btheta)$ (possibly up to a linear transformation).
The drift of this diffusion process is given by the posterior 
expectation of (a linear function of) $\btheta$ given data $\data$ and additional observations that are corrupted by Gaussian noise.
\item[$(ii)$] Use a combination of variational methods and message passing algorithms
to approximate the posterior mean of $\btheta$ given $\data$ and additional
noisy  observations,
 which can be computed efficiently.
\end{itemize}
Given an accurate approximation of the posterior mean as per point $(ii)$,
the diffusion process mentioned above can be approximated by applying a standard Euler discretization.
To ensure the stability of this discretization, we prove that the approximate posterior mean is sufficiently regular as a function of the data. 
Proving such regularity results will be an important 
technical contribution of the paper.

Such sampling approach was first proposed in \cite{alaoui2022sampling}  to address 
a problem from statistical physics:
sampling from the Sherrington-Kirkpatrick Gibbs measure at a  high-temperature
(see Section \ref{sec:Related} for more pointers to this line of work).
However, this Gibbs measure has no statistical meaning and several of the techniques developed 
in \cite{alaoui2022sampling} are tailored to a `very noisy' regime that is not of 
statistical interest.
The present paper describes the first rigorous application of this
 approach  to  sampling problems in high-dimensional statistics.

The approach of \cite{alaoui2022sampling} was developed as 
an algorithmic version of the \emph{stochastic localization} (SL) process, introduced by Eldan
as a proof technique
\cite{eldan2013thin,eldan2020taming,eldan2022analysis,chen2022localization}.
These works provide a wealth of mathematical results about this process.

In a parallel development, the SL diffusion process 
 has attracted considerable amount of work in machine learning.
Under the name of `denoising diffusion,' it
has become the state-of-the-art technique in generative AI
\cite{song2019generative,ho2020denoising,song2020improved,songscore}.
This connection is discussed in detail in
 \cite{montanari2023sampling}.
 We refer to Remark \ref{rem:GenerativeDiff} for a very succinct description.

In the rest of this introduction we will outline the construction of the general
sampling algorithm and describe two specific sampling problems to which we apply our approach. 
We will then provide an informal summary of our results. 
%
%
\subsection{The diffusion process and the associated sampling algorithm}
\label{sec:IntroDiffusion}

Stochastic localization is defined in \cite{eldan2013thin,eldan2020taming,chen2022localization}  
as a measure-valued stochastic process. Here we follow an alternative viewpoint
that, after suitable generalizations, is actually equivalent to the original one
 \cite{montanari2023sampling}.

Let $\bH\in\reals^{N \times n}$ be a matrix designed by the statistician:
in particular, $\bH$ can depend on the data $\data$. 
In what follows, whenever we condition on $\data$,
it is understood that we condition on $\bH$ as well. In many cases of interest, 
the choice $\bH=\id_n$ is adequate ad hence the discussion below can be simplified.
 
For  $\btheta \sim \mu_{\data}$, we define the linear observation process 
$\{\by(t)\}_{t\ge 0}$ via
\begin{align}
\by(t) = t\, \bH\, \btheta +\bG(t)\, ,\label{eq:GeneralLinearObs}
\end{align}
where $\{\bG(t)\}_{t\ge 0}$ is a standard Brownian motion in $\RR^N$ and is independent 
of $\btheta$. 

The posterior distribution of $\btheta$ given 
$\data$ and $\by(t)$ takes the form:
\begin{align}
\mu_{\data, t}(\de \btheta) & =  \prob(\de\btheta |\data,\by(t))\nonumber\\ 
& = \frac{1}{Z(\data, \by(t), t)}
 \exp\left\{ \langle \by(t), \bH \btheta \rangle - \frac{t}{2}\|\bH\btheta\|_2^2 \right\}
\mu_{\data}(\dd \btheta) \, ,\label{eq:FormulaMu}
\end{align}
where $Z(\data, \by(t), t)$ is a normalizing constant.
We define two posterior mean functions as follows:
\begin{align}
\bm(\by, t)  := \E[\bH\btheta|\data,t\bH\btheta+\bG(t)=\by]\, , \;\;\;\;\;\;\;
\bm_{\btheta}(\by, t) := \E[\btheta|\data,t\bH\btheta+\bG(t)=\by]\, .
\end{align}
Of course,  $\bm(\by, t)  = \bH\bm_{\btheta}(\by, t)$, but distinguishing the two
functions is useful when considering estimators that approximate each one of them.
The process $(\by(t))_{t\ge 0}$ admits an alternative characterization as the unique solution
of the following stochastic differential equation (see, e.g., \cite{alaoui2022sampling} 
or  \cite[Theorem 7.1]{liptser1977statistics} 
for a derivation):
\begin{align}
	\by(t) = \bm(\by(t), t) \dd t + \dd \bB(t), \qquad \by(0) = \mathbf{0}_N\, ,\label{eq:general-SDE}
\end{align}
with $\{\bB(t)\}_{t \geq 0}$ being a standard Brownian motion in $\RR^N$. 

As $t\to\infty$, it holds that $\by(t)/t = \bH\btheta+(\bG(t)/t) \overset{a.s.}{\to} \bH\btheta$, where 
we recall that $\btheta\sim\mu_{\data}$.  
Since $\bH$ is known, the desired sample $\btheta$ 
can often be computed from $\bH\btheta$.
The simplest example is when  
 $\bH$ has full column rank.  
 We can then generate an approximate sample
 via $\btheta^{\salg}= \bH^{\dagger}\by(T)/T$ for a large enough $T$.  
 More generally, we can set $\btheta^{\salg}= \bm_{\btheta}(\by(T), T)$.
 Further discussion of this point is 
 deferred to Section \ref{sec:MainLinear}.

 In order to sample from $\mu_{\data}$, it is therefore sufficient to approximate
  the process \eqref{eq:general-SDE}.  In a nutshell, we proceed as follows:
\begin{enumerate}
	\item We discretize the stochastic differential equation (SDE) \eqref{eq:general-SDE} in time,
	over the time interval $[0,T]$, using an Euler scheme. 
	\item We construct an efficient algorithm that provides 
	 an approximation $\hbm(\by, t)$  of the posterior expectation $\bm(\by, t)$. 
	\item Let $(\hby(t))_{t\in [0,T]}$ denote the approximate diffusion, we
	generate a sample $\btheta^{\salg}$ by outputting $\hbm_{\btheta}(\hby(T), T)$,
	where $\hbm_{\btheta}(\,\cdot\,, T)$ is a suitable approximation of  $\bm_{\btheta}(\,\cdot\,, T)$.
\end{enumerate}

The most complex part of the above procedure is the computation of the posterior mean
$\hbm(\by, t)$.
In the sequel, we offer general sampling guarantees under two conditions:
 $(1)$~A probabilistic bound on the distance between 
 $\bm(\by, t)$, and $\hbm(\by, t)$, and $(2)$~Approximate Lipschitz continuity of the mapping
 $\by \mapsto \hbm(\by, t)$.

Our approach to construct $\hbm(\by, t)$ is based on the Approximate Message Passing
(AMP) algorithm and on variational inference methods.
Namely, we construct a suitable free energy functional $(\bm,\by,t)\mapsto\cuF(\bm;\by,t)$
for which we expect 
\begin{align}
\bm(\by,t) \approx \arg\min_{\bm}\cuF(\bm;\by,t)\, .
\end{align}
While the free energy functional $\cuF(\bm;\by,t)$ is --in general-- non-convex, 
we show that a suitable AMP algorithm can be used to find an approximate local minimum with the desired properties.

\begin{remark}\label{rem:GenerativeDiff}
As mentioned in the previous section, the `denoising diffusions' 
method in generative AI
\cite{song2019generative,ho2020denoising,song2020improved,songscore}
is based on the same diffusion process \eqref{eq:general-SDE}.
In generative AI
 the target distribution $\mu$
is not (as in our case) a Bayesian posterior, but instead the distribution of a certain data 
source, e.g., natural images of a certain type.
As a consequence, the form of $\mu$ is not known explicitly, 
but we have instead access to samples $\btheta_1$, \dots $\btheta_{N_{\mathrm{train}}}\sim_{iid} \mu$. 
This samples are used to construct an estimate the posterior expectation $\hbm(\by,t)$.
New samples are generated by implementing numerically a discretization of
 \eqref{eq:general-SDE}.

A second (less important) difference is that denoising diffusions are typically
derived using a different argument than stochastic localization,
based on `time reversal.' As a consequence, the resulting diffusion does
not coincide with  \eqref{eq:general-SDE}, but is equivalent to it by a change of variables,
as detailed in \cite{montanari2023sampling}.
\end{remark}

\subsection{High-dimensional linear regression and spiked matrix models}
\label{sec:IntroModels}

Most of our technical effort will be devoted to implementing the above general program in
two canonical problems: Bayesian inference  in
 spiked matrix models and in high-dimensional regression models. We briefly outline these settings below. 
 
 \paragraph{Spiked matrix model.} Given a signal-to-noise parameter $\beta > 0$, we observe an
 $n \times n$ symmetric matrix $\bX$ generated as follows:
\begin{align}\label{eq:model}
	\bX = \frac{\beta}{n}\btheta \btheta^{\sT} + \bW. 
\end{align}
Hence, for this example we have $\data = \bX$. 
Here, $\bW  \sim \GOE(n)$ and is independent of $\btheta$, i.e., $\bW$ is an $n \times n$ 
symmetric matrix with independently distributed entries above the diagonal: 
$\{W_{ii}: i \in [n]\} \iidsim \normal(0, 2 / n)$, $\{W_{ij}: 1 \leq i < j \leq n\}
 \iidsim \normal(0, 1 / n)$. 
As for $\btheta$, we assume a product prior:
 \begin{align}
  (\theta_{i})_{i\le n} \iidsim \pi_{\Theta}\, . \label{eq:PriorTheta}
  \end{align}
Throughout, we will denote by  $\pi_{\Theta}^{\otimes n}(\dd \btheta)= \pi_{\Theta}(\dd \theta_1)\cdots
 \pi_{\Theta}(\dd \theta_n)$ the product distribution over $\RR^n$ with marginal $\pi_{\Theta}$.
  
Model
   \eqref{eq:model} is the symmetric version  of the \emph{spiked model} 
first introduced in \cite{johnstone2001distribution}.
In the asymmetric (rectangular) version, data takes the form 
$\tbX = \bu\btheta^{\sT}+\bZ\in\reals^{n\times p}$, where  
$(u_i)_{i\le n} \iidsim \pi_{U}$, $(\theta_{i})_{i\le p} \iidsim \pi_{\Theta}$,
and $\bZ$ is a noise matrix. While we carry out our analysis in the symmetric
setting for simplicity, the generalization to asymmetric matrices is
straightforward. 

Under model \eqref{eq:model} with prior \eqref{eq:PriorTheta}, given observation $\bX$, the posterior takes the form: 
\begin{align}
	\mu_{\bX}(\dd\btheta) \propto \exp \left( \frac{\beta}{2} \< \btheta, \bX \btheta \> 
	-\frac{\beta^2}{4n}\|\btheta\|_2^4\right) \pi^{\otimes n}_{\Theta}(\dd \btheta)\, .\label{eq:PosteriorFirst}
\end{align}
The problem of sampling from $\mu_{\bX}$ is already interesting if
the prior is the uniform distribution over $\{+1,-1\}^n$, i.e.,
$\pi_{\Theta} = \Unif(\{+1,-1\})$. 
 In this case, model \eqref{eq:model} is also known as $\ZZ_2$-synchronization and is closely
 related to the stochastic block model
  \cite{singer2011angular,deshpande2017asymptotic,abbe2018group}.
 Sampling from this measure is known to  be $\#$P-complete in the worst case
\cite{sly2012computational,galanis2016inapproximability}.

The symmetric spiked model \eqref{eq:model} has attracted considerable amount of work
as an idealized setting for low-rank matrix  estimation. 
   Consider, to be definite, the case $\E[\Theta]:=\int \theta\, \pi_{\Theta}(\de\theta)=0$.
   Then,  it is known that non-trivial estimation of $\btheta$
   is possible  if $\beta>\beta_{\sIT}$, with $\beta_{\sIT}$ a constant that was first rigorously characterized
   in  \cite{lelarge2019fundamental}. 
   On the other hand, for  $\beta<\beta_{\sIT}$, 
   the posterior measure is very close to the prior, hence sampling
   from the posterior in this case is not interesting. 
 Polynomial-time algorithms that
   achieve non-trivial estimation
  are known to exist for $\beta>\beta_{\swr} := 1/\E[\Theta^2]$, while
  they are conjectured not to exist for $\beta<\beta_{\swr}$ ($\beta_{\swr}$
  is referred to as the weak recovery threshold), see e.g. 
  \cite{deshpande2017asymptotic,lelarge2019fundamental,montanari2021estimation}.
  Finally, the results of \cite{lelarge2019fundamental,montanari2021estimation}
  characterize a set $J(\pi_{\Theta}) \subseteq (\beta_{\swr},\infty)$ such that, 
   for $\beta \in J(\pi_{\Theta})$ efficient algorithms (Bayes AMP) are known to achieve
   asymptotically Bayes-optimal error. On the other hand, this is conjectured to be
  hard for $\beta \in  (\beta_{\sIT},\infty)\setminus\overline{J}(\pi_{\Theta})$. 
  in Section \ref{sec:MainSpiked} ---under some regularity assumptions on $\pi_{\Theta}$,
 $J(\pi_{\Theta})$ includes all $\beta>\beta_{\salg}(\pi_{\Theta})$ for a constant 
 $\beta_{\salg}(\pi_{\Theta})<\infty$. We use $\beta_{\salg}(\pi_{\Theta})>0$ to denote the 
 minimum value of such constant.
 Since posterior sampling is not easier than posterior estimation, 
 $\beta_{\salg}(\pi_{\Theta})$ is also the infimum value of $\underline{beta}$
 such that we expect to be able to sample efficiently for all $\beta\in (\underline{\beta},\infty)$.

 \paragraph{High-dimensional regression.}
In this context, we observe a covariate matrix $\bX \in \RR^{n \times p}$ and a response vector 
$\by_0\in \RR^n$,
which we assume follows a linear model:
\begin{align}
\label{eq:LR}
	\by_0 = \bX \btheta + \bepsilon \in \RR^n. 
\end{align}
Here, $\btheta \in \RR^p$ is an unknown vector of coefficients, 
and $\bepsilon \sim \normal(\bzero,\sigma^2\id_n)$ is a noise vector. 
We further assume $\btheta$ follows a product prior: $(\theta_{i})_{i\le p} \iidsim \pi_{\Theta}$, 
where $\pi_{\Theta}$ is a fixed distribution over $\RR$, and $\bX$ is independent of $\btheta, \beps$.
In this problem, the observed data is $\data = (\bX, \by_0)$. 
Given observations $(\bX,\by_0)$, the posterior distribution of $\btheta$
reads
\begin{align}
	\mu_{\bX,\by_0}(\dd\btheta) \propto \exp \left( -\frac{1}{2\sigma^2} \big\|\by_0-\bX\btheta\big\|_2^2\right)
	\pi^{\otimes p}_{\Theta}(\dd \btheta)\, .\label{eq:PosteriorFirst_LR}
\end{align}
If $n/p\to\infty$, it is easy to see that the above posterior concentrates around the 
true coefficients vector (under suitable conditions on the design matrix $\bX$) in the sense
that $\E[\int \bfone(\|\tilde\btheta-\btheta\|^2_2 / p>\eps)\mu_{\bX,\by_0}(\dd \tilde\btheta)] \to 0$ for  arbitrarily small positive constant $\eps$.
Hence posterior sampling is significantly easier in this regime \cite{bontemps2011bernstein}. 

In contrast, in this work we consider the proportional asymptotics in which $n,p\to\infty$
with $n/p\to\delta\in (0,\infty)$. In this regime, the asymptotics
of the Bayes optimal estimation error was characterized in \cite{barbier2019optimal} 
 for design matrices $\bX$ with i.i.d. entries. 
For this setting, the results of \cite{barbier2019optimal} characterize
 a set $K(\pi_{\Theta},\sigma^2)\subseteq \reals$ such that, for $\delta\in K(\pi_{\Theta},\sigma^2)$
 Bayes AMP achieves the Bayes optimal error, while this is conjectured to be hard for
 $\delta\in \reals\setminus \overline{K}(\pi_{\Theta},\sigma^2)$. 
 As discussed in Section \ref{sec:MainLinear}, under suitable conditions
 on $\pi_{\Theta}$, $K(\pi_{\Theta},\sigma^2)$ contains a semi-infinite interval 
  $(\delta_{\salg},\infty)$, for a certain (minimal)
 critical value $\delta_{\salg} = \delta_{\salg}(\pi_{\Theta},\sigma^2)$.
For general $\pi_{\Theta}$ and the noise level $\sigma^2$,
$K(\pi_{\Theta},\sigma^2)$ is a strict subset of $\reals$.
In particular $\delta_{\salg}$ is the infimum value of $\underline{\delta}$
such that we expect efficient sampling to be possible for all 
$\delta\in(\underline{\delta},\infty)$.


\subsection{Contributions}
 
 We present two types of contributions: $(i)$~A general framework to 
 construct algorithms for Bayes posterior sampling and prove correctness of these algorithms;
 $(ii)$~Instantiations of this framework to sample from spiked matrix models
 and high-dimensional linear regression. Proving correctness for the algorithms
 of point $(ii)$ will require to establish that the proposed posterior mean estimate is accurate
 and Lipschitz continuous in its argument.
 
 More specifically, a summary of our results is given below.

\paragraph{A general framework to prove approximate sampling.} 
We develop a general algorithmic framework that uses an oracle to approximate
 the posterior mean of $\bH\btheta$ (given $\data$  and an additional observation $\by(t)$),
 and outputs approximate samples from the posterior $\mu_{\data}$.
Our algorithm follows the outline given in Section \ref{sec:IntroDiffusion}.
We provide a general upper bound on the distance between the distribution
of samples generated by our algorithm and the target posterior distribution.
The upper bound is a function of the posterior mean estimation accuracy and the Lipschitz modulus
 of the oracle, see Section \ref{sec:GenLinear} for details.
 
\paragraph{Posterior sampling for spiked matrix model.} 
In Section \ref{sec:MainSpiked} we describe a specific implementation of 
our algorithm when the goal is to sample from the posterior distribution \eqref{eq:PosteriorFirst}
under the spiked model \eqref{eq:model}.
In this case, we set $\bH=\id_n$, and  
 the overall algorithm has complexity $O(n^2)$.

We prove that there exists 
a constant $\beta_0$, such that, for $\beta\ge \beta_0$, the samples generated
by the algorithm have distribution that is close to the actual posterior 
\eqref{eq:PosteriorFirst}. Here, ``close'' means that the Wasserstein distance between these
two distributions satisfies $W_2(\mu_{\bX},\mu^{\salg}_{\bX}) = o_P(\sqrt{n})$.
We prove this result for two classes of priors:
$(i)$~Discrete priors $\pi_{\Theta}$ supported on a finite number of points; 
$(ii)$~Continuous
priors with density $\pi_{\Theta}(\de\theta) = \exp(-U(\theta))\de\theta$, with $\|U''\|_{\infty}<\infty$
(but not necessarily log-concave\footnote{Let us emphasize that the overall log-prior
$\sum_{i\le d}U(\theta_i)$ is not a bounded (in $d$) perturbation of a convex function 
and hence fast mixing does not follow from the standard Holley-Stroock criterion.}).

As mentioned above, we require $\beta$ to be larger than a constant for
non-trivial estimation to be possible.
In such constant $\beta$ regimes, the posterior  $\mu_{\bX}$ 
is  non-trivially aligned with the true coefficients $\btheta$, 
hence the sampling task is non-trivial.

\paragraph{Posterior sampling for linear regression.} 
In Section \ref{sec:MainLinear}, we introduce an algorithm that samples from
the posterior distribution \eqref{eq:PosteriorFirst_LR} under the linear
model \eqref{eq:LR}. 
In this case, we choose $\bH=\bX$ and $\data = (\bX, \by_0)$, and  
 the associated algorithm has complexity $O(np)$.
 
For our algorithm, we prove that there exists $\delta_0=\delta_0(\pi_{\Theta},\sigma^2)$,
 such that if $\delta>\delta_0$, then the proposed algorithm produces
 samples with distribution that is close to the target posterior $\mu_{\bX,\by_0}$, namely, 
 $W_2(\mu_{\bX,\by_0},\mu^{\salg}_{\bX,\by_0}) = o_P(\sqrt{p})$.
 Again, in this regime non-trivial estimation is possible, but
 at the same time sampling is non-trivial.

\begin{remark}[Conjectured regime of validity]
 We expect stochastic localization sampling to be successful for a broader range of 
values of $\beta$ and $\delta$. In particular, we expect sampling to be possible
for all $\beta>\beta_{\salg}(\pi_{\Theta})$ (for the spiked matrix model)
or all $\delta>\delta_{\salg}(\pi_{\Theta})$ (for the linear regression model).
As we will see, proving Lipschitz continuity of the AMP estimator or the variational inference method 
is the current bottleneck towards reaching this goal.
\end{remark}

\subsection{Notations}
\label{sec:notation}

For $n \in \NN_{>0}$, we define the set $[n] := \{1,2, \cdots, n\}$. 
We denote by $\Sym(n)$ the collection of all $n \times n$ symmetric real matrices, 
by $\Sym_{\ge}(n)$ the subset of positive semidefinite matrices, and by $\Orth(n)$ the collection of $n \times n$ orthogonal matrices. 
We denote by $\bA^{\dagger}$ the pseudoinverse of matrix $\bA$.  We denote the Euclidean
ball in $\reals^n$ with center $\ba$ and radius $r$ by $\Ball^n(\ba,r)$.
Given a differentiable map $\bF:\reals^m\to\reals^n$, we denote by 
$\bD\bF(\bx)\in\reals^{n\times m}$ its differential (Jacobian matrix)
at $\bx\in \reals^m$. 

 We denote by $\plim$ limit in probability. The set of probability measures over the
 Borel measurable space $(\reals^n,\cB_{\reals^n})$ is denoted by $\cuP(\reals^n)$,
 and the set of probability measures with finite second moment is denoted by $\cuP_2(\reals^n)$.
 We write $\Law(X)$ to represent the probability distribution of a random variable (or vector) 
 $X$. 

The $W_2$ Wasserstein distance of two  probability measures $\mu, \nu$ over $\reals^n$
is denoted by 
$$W_2(\mu,\nu) :=\inf_{P\in\cC(\mu,\nu)}\E_{(\bX,\bY) \sim P}[\|\bX-\bY\|^2_2]^{1/2},$$
where the infimum is taken over all joint distributions with marginals $\mu$ and $\nu$. We will also consider
the normalized Wasserstein distance $W_{2,n}(\mu,\nu):=W_{2}(\mu,\nu)/\sqrt{n}$.

For a $k$-th order tensor $\bT$, we define its operator norm as 
\begin{align*}
	\|\bT\|_{\op} = \sup_{\bv_1,\dots \bv_k:\|\bv_i\|=1} \langle \bT, 
	 \bv_1\otimes \cdots\otimes \bv_k\rangle\, .
\end{align*}
%
%
%
\section{Further related work}
\label{sec:Related}

As mentioned in the introduction, 
MCMC methods are
the dominant approach to sampling from Bayes posteriors. 
In the case of non-log-concave and possibly discrete posteriors such as 
Eq.~\eqref{eq:PosteriorFirst} of interest here, Gibbs sampling (a.k.a. Glauber dynamics)
would probably be the method of choice. 
This Markov chain updates one coordinate at each step according to the conditional
 probability distribution given all the other coordinates. 
 Classical methods to bound the mixing times of such Markov chains
 are based on the so-called Dobrushin condition   \cite{dobrushin1968description}
 and require (in the present case) $\beta\lesssim 1/\sqrt{n}$.
 Over the last two years, remarkable breakthroughs were achieved establishing Markov chain
 mixing results under the much weaker `spectral mixing conditions'
  \cite{bauerschmidt2019very,eldan2022spectral,anari2021entropic}.  
  However, existing results only apply to the case of $\ZZ_2$-synchronization
  (i.e., model \eqref{eq:model} with $\btheta\in\{+1,-1\}^n$). 
  In this setting, these papers established mixing time upper bounds under the condition
  $\beta(\lambda_{\max}(\bX)-\lambda_{\min}(\bX))\le 1$.
Applying classical results on extremal eigenvalues of spiked random matrices \cite{baik2005phase},
 this condition amounts to  $\beta<1/4$. 
 We note that this is not a regime of statistical
 interest, since for $\beta<1$ it is impossible to produce an estimator
 with non-vanishing correlation with the true signal  \cite{deshpande2017asymptotic}.

Variational methods 
\cite{jordan1999introduction,wainwright2008graphical,mezard2009information,blei2017variational} offer 
another popular framework for conducting approximate Bayesian analysis.
In the challenging regime with a constant signal-to-noise ratio as considered here,
naive mean field methods are known to yield incorrect estimates of the posterior mean.
Namely,  \cite{ghorbani2019instability} shows that, even in the case of 
$\integers_2$-synchronization the naive mean field estimate $\hbm_{\sNMF}$ suffers an error
$\|\hbm_{\sNMF}-\bm\|_2\ge c\sqrt{n}$ even in cases in which efficient algorithms to estimate $\bm$
exist.

However, a recent line of work proves that, for low-rank matrix estimation and high-dimensional regression,
asymptotically correct estimates of the posterior 
can be achieved if the mean field free energy is modified by adding a Thouless-Anderson-Palmer (TAP) 
correction \cite{fan2021tap,celentano2021local,qiu2022tap}. Namely,  \cite{fan2021tap}
proved (for $\integers_2$ synchronization) that all local minima of the TAP free energy 
below a certain level are concentrated in a small neighborhood of the 
posterior mean\footnote{To obtain a non-trivial posterior mean, the overall symmetry $+1/-1$ needs 
to be factored out.};
 \cite{celentano2021local}
proved that in fact the TAP free energy is locally strongly convex around this point and hence
there is only one such local minima (up to symmetries). Finally, \cite{qiu2022tap} studied TAP free energy in the context of high-dimensional linear regression, and showed that the TAP approximation consistently estimates the log-normalizing constant of the posterior distribution.  
 Our work leverages some of these results to construct a sampling procedure. 
 
The recent work \cite{koehler2022sampling} develops a sampling algorithm for Ising models that
merges ideas from MCMC and variational inference. However, even for the case of
$\ZZ_2$-synchronization, this approach falls short of providing an algorithm that is
 successful in the regime of $\beta$ that is of statistical interest. 
  
As a side remark, a substantial literature characterizes the optimal estimation error,
efficient estimation algorithms, and the potential computational barriers for the spiked model \eqref{eq:model}. 
These works also offer generalizations to the multi-rank setting and asymmetric matrices.
A subset of these works include 
\cite{deshpande2017asymptotic,dia2016mutual,miolane2017fundamental,el2018estimation,
lelarge2019fundamental,barbier2016mutual,barbier2019adaptive, 
chen2022hamilton,montanari2021estimation,montanari2022fundamental}.

As already mentioned, the present work  is closely related to 
\cite{alaoui2022sampling}, which first uses an algorithmic implementation of the stochastic 
 localization process, in conjunction with an AMP approximation of the posterior expectation to construct a sampling algorithm. 
 However, \cite{alaoui2022sampling} focuses uniquely on the Sherrington-Kirkpatrick model,
 i.e., assumes $\btheta\in\{+1,-1\}^n$, 
 	given by $\mu_{\bW}(\btheta) \propto \exp( \beta \< \btheta, \bW \btheta \>/2)$,
 	whereby $\bW\sim\GOE(n)$. 
 	Namely, there is no spike in this model.
 	For the Sherrington-Kirkpatrick model,  \cite{alaoui2022sampling}  established sampling guarantees for $\beta < 1 / 2$. 
 	This threshold was later improved to $\beta<1$ by  \cite{celentano2022sudakov}. Hardness results for
 	$\beta>1$ were also established in \cite{alaoui2022sampling} for a class of `stable' algorithms.
 	
Despite its significant role in  physics, 
the Sherrington-Kirkpatrick model is not of statistical interest:
 the data is ``pure noise'', and the probability measure is not a Bayes posterior. 
 To the best of our knowledge, this paper presents the first application of diffusion-based
  methods to high-dimensional statistical inference.

In concurrent works, guarantees for sampling via diffusions were recently established 
in \cite{lee2022convergence,chen2022sampling,chen2022improved}.
Two important differences with respect to the present paper are: $(1)$~The
results of \cite{lee2022convergence,chen2022sampling,chen2022improved} \emph{assume} the
 existence of an accurate posterior mean estimator, while
we construct it; $(2)$~In these papers, the posterior mean estimate is assumed to be very accurate, 
namely $\E\{\|\hbm(\by,t)-\bm(\by,t)\|^2\} = o_n(1)$. Under this assumption
it is possible to establish strong sampling guarantees.
 
In many circumstances only the weaker error bound
$\E\{\|\hbm(\by,t)-\bm(\by,t)\|^2\} = o_n(n)$ is accessible. 
Our proof technique allows us to derive useful sampling guarantees even in this more 
challenging scenario.

\paragraph{Follow-up work.} After a first version of this manuscript was posted online, several 
groups developed the approach proposed here in several other interesting directions. 
Among others, \cite{huang2024sampling} studied the mixed spherical spin glass from statistical physics,
and constructed a posterior mean estimator satisfying $\E\{\|\hbm(\by,t)-\bm(\by,t)\|^2\} = o_n(1)$.
Using the same construction as in \cite{alaoui2022sampling} and in the present paper,
this yields a sampling algorithm with guarantees in total variation distance. 
\cite{mei2023deep} studied a general class of sampling problems for graphical models, and
obtained algorithms with the guarantee $D_{\sKL}(\mu\|\mu^{\salg}) = o(n)$. 
A similar 
guarantee was established for posterior sampling in high-dimensional linear regression
\cite{cui2024sampling}, using again the same algorithm studied here. 

It is worth mentioning that the bound $D_{\sKL}(\mu\|\mu^{\salg}) = o(n)$
does not imply any non-trivial bound on Euclidean distances between distributions. 
Namely, it
is easy to construct probability measures
$\mu_1, \mu_2$ on $\reals^n$ such that, for some constant $R$,
and $\|\ba_1-\ba_2\|_2\ge 3R\sqrt{n}$:
\begin{enumerate}
\item $D_{\sKL}(\mu_1\|\mu_2) = o(n)$.
\item $\mu_1(\Ball(\ba_1,R\sqrt{n}))\ge 1-n^{-1}$, $\mu_2(\Ball(\ba_2,R\sqrt{n}))\ge 1-n^{-1}$.
\end{enumerate}
In words, the two probability measures are essentially supported on far apart balls, and yet
they are $o(n)$ apart in KL 
distance\footnote{For instance, let ${\sf U}_i$ be the uniform measure over 
$\Ball(\ba_i,R\sqrt{n}))$ and $mu_i = (1-p_i) \, {\sf U}_1+ p_i\, {\sf U}_2$
 for suitable choices of $p_1,p_2$.}. This cannot happen under the $W_2$ guarantee 
established in the present work.
%
%
\section{General sampling scheme}
\label{sec:GenLinear}

The sampling procedure outlined in Section \ref{sec:IntroDiffusion}
is made precise in  \cref{alg:general} below. 
This algorithm makes use of two oracles $\hbm(\by,t)$ and $\hbm_{\btheta}(\by, t)$ 
that separately approximate $\bm(\by, t)$ and $\bm_{\btheta}(\by, t)$.
%
\begin{algorithm}
\caption{General sampling scheme}\label{alg:general}
\textbf{Input: }Parameters $(L, \Delta)$;	
\begin{algorithmic}[1]
\State Set $\hat{\by}_0 = \mathbf{0}$;
	\For{$\ell = 0,1,\cdots, L - 1$}
		\State Draw $\bw_{\ell + 1} \sim \normal(0,\id_N)$, independent of everything so far;
		\State Update $\hby_{\ell + 1} = \hby_\ell + \Delta 
		\hbm (\hby_\ell,\Delta\ell)  + \sqrt{\Delta} \bw_{\ell + 1}$;	
	\EndFor
\State \Return $\btheta^{\salg} = \hbm_{\btheta}(\hat{\by}_L,L\Delta)$; 
\end{algorithmic}
\end{algorithm}

Note that the oracle $\hbm_{\btheta}(\by, t)$ is only used at the final time 
$t = T:=L\Delta$. 
Hence, we only need to approximate the posterior 
expectation of $\btheta$ given $\by(T)/T = \bH\btheta+(\bg/\sqrt{T})$,
where $\bg \sim\normal(0,\id_N)$.
For large $T$, this corresponds to very low noise.
Constructing such an oracle is relatively easy in a number of circumstances, 
as illustrated by the following two examples:
\begin{itemize}
\item If $\bH$ has full column rank (in particular, $N\ge n$), then we can construct an estimator by taking
$\hbm_{\btheta}(\by;T) = \bH^{\dagger}\by(T)/T$.
\item If $\bH$ does not have full column rank (e.g., $N<n$), but  $\pi_{\Theta}$
induces sparsity in $\btheta$. In this case, standard techniques from
 compressed sensing and high-dimensional sparse regression literature can be brought to bear \cite{tibshirani1996regression,donoho2006compressed,candes2008introduction,buhlmann2011statistics}. 
\end{itemize}
In contrast, the oracle $\hbm(\by, t)$ is 
required for all $t$ that appears in the discretization scheme. However, one can hope to exploit the
freedom to choose the matrix $\bH$ to simplify this task.

We provide a theoretical guarantee for \cref{alg:general} under the following assumptions
(probabilities below are with respect to the distribution of the stochastic localization process 
$\{\by(t)\}$, at fixed $\bD$):
\begin{enumerate}
	\item[\hypertarget{A1}{$\mathsf{(A1)}$}] (Posterior mean consistency) With probability at least 
	$1 - \eta$, it holds that
	\begin{align*}
	 \max_{\ell \in \{0\} \cup [L - 1]} 	\frac{1}{\sqrt{N}}
	 \big\| \bm(\by(\ell\Delta), \ell\Delta) - 
		\hbm(\by(\ell\Delta),\ell\Delta) \big\|_2 \leq \eps_1. 
	\end{align*}
	Further, with the same probability,
	$\| \bm_{\btheta}(\by(T), T) - \hbm_{\btheta}(\by(T),T) \big\|_2 \leq \eps_1\sqrt{n}$, where we recall that $T = L \Delta$.  
	\item[\hypertarget{A2}{$\mathsf{(A2)}$}] (Path regularity) With probability at least $1 - \eta$, 
	it holds that
	\begin{align*}
\max_{\ell \in \{0\} \cup [L - 1]}\sup_{t \in [\ell\Delta, (\ell + 1) \Delta]} \frac{1}{\sqrt{N}} 
\|\bm(\by(t), t) - \bm(\by(\ell \Delta), \ell\Delta)\|_2 \leq C_1\sqrt{\Delta} + \eps_2.
	\end{align*}
	\item[\hypertarget{A3}{$\mathsf{(A3)}$}] (Lipschitz continuity) There
	exists a sequence $\{r_\ell\}_{0\le \ell\le L} \subseteq \RR_+$, such that 
	letting $B(\ell) :=\{\by\in\reals^N:\;\|\by-\by(\ell\Delta)\|_2\le r_\ell\sqrt{N}\}$, then  
	the following holds with probability at least 
	$1 - \eta$:
	\begin{align*}
	\max_{\ell\in \{0\} \cup [L - 1]} 
		\sup_{\by_1\neq \by_2\in B(\ell)} \left[
		\frac{1}{\sqrt{N}} {\|\hbm(\by_1,\ell\Delta) - \hbm(\by_2,\ell\Delta)\|_2} - C_2 
		\frac{1}{\sqrt{N}} {\|\by_1-\by_2\|_2} \right]
		\leq \eps_3 \, .
	\end{align*}
	Further, we assume that $r_\ell> (C_1\sqrt{\Delta} + \eps_1 + \eps_2 + \eps_3) e^{C_2\ell\Delta}/C_2$
	for all $\ell \in \{0\} \cup [L]$.
	We also assume that with the same probability,
	$\| \hat  \bm_{\btheta}(\by_1, T) - \hbm_{\btheta}(\by_2,T) \big\|_2 / \sqrt{n} \leq C_2\|\by_1-\by_2\|_2 / \sqrt{N} + \eps_3$
	for all $\by_1,\by_2 \in B(L)$. 
\end{enumerate}
The dependence on constants $C_1, C_2$  will be tracked in the statement below.
%
\begin{theorem}\label{thm:gen}
Assume that $\|\hbm_{\btheta}(\by,T)\|_2 \le \oR\sqrt{n}$ for all $\by \in \RR^N$ (this can always be achieved by projection
onto the ball $\Ball^n(\bzero,\oR\sqrt{n})$)
and that conditions \hyperlink{A1}{$\mathsf{(A1)}$}, \hyperlink{A2}{$\mathsf{(A2)}$} and
\hyperlink{A3}{$\mathsf{(A3)}$} hold.
Letting $\mu_{\data}^{\salg}:=\Law(\btheta^{\salg})$ be the distribution of
 the samples generated by \cref{alg:general},  then we have:
 	\begin{align}\label{eq:ZeroGeneral}
 	\begin{split}
 		W_{2, n}(\mu_{\data}, \mu_{\data}^{\salg})
		\leq &  2(C_1 \sqrt{\Delta} + \eps_1 + \eps_2 + \eps_3)\cdot e^{C_2T} + 
		\frac{1}{n}\mu_{\data}\big(\|\btheta\|_2^2\cdot\bfone_{\|\btheta\|_2\ge \oR\sqrt{n}}\big)^{1/2}+ 10 \oR \eta \\
		&  + W_{2, n}\big(\mu_{\data}, \Law( \bm_{\btheta}(\by(T), T))\big).
	\end{split}
 	\end{align}
 	%
	%
	If in addition $\bH$ has full column rank, and
	 $\int(\|\btheta\|^2_2/n)^{c_0} \mu_{\data}(\de\btheta) \le R^{2c_0}$ for some $c_0>1$,
	 then 
	\begin{align}
	\begin{split}
	 W_{2, n}(\mu_{\data}, \mu_{\data}^{\salg})
		\leq & 2(C_1 \sqrt{\Delta} + \eps_1 + \eps_2 + \eps_3)\cdot e^{C_2T}+
		C(c_0) R\eta^{(c_0-1)/c_0}
		+{\frac{1}{\sqrt{T}} \Tr\big((\bH^{\sT}\bH)^{-1}\big)^{1/2}},\, .
		\label{eq:SecondGeneral}
	\end{split}
	\end{align}
	where $C(c_0)$ is a positive constant that depends only on $c_0$. 
\end{theorem}
The proof of this theorem is deferred to Appendix \ref{app:ProofGen}.

It is useful to compare Theorem \ref{thm:gen}
with concurrent in the literature on denoising diffusions.
The key difference is that we make minimal assumptions on the accuracy 
of the posterior mean. Assumption \hyperlink{A2}{$\mathsf{(A2)}$} only requires
an error bound of the form
$\E\{\|\hbm(\by,t)-\bm(\by,t)\|^2\} = o_N(N)$, while \cite{lee2022convergence,chen2022sampling,chen2022improved}.
require the much stronger condition $\E\{\|\hbm(\by,t)-\bm(\by,t)\|^2\} = o_N(1)$,
which implies stronger guarantees (in total variation). 
Since the $\E\{\|\hbm(\by,t)-\bm(\by,t)\|^2\} = o_N(1)$ 
is extremely difficult to prove in the present context (see 
\cite{huang2024sampling}), we cannot  use the results of 
\cite{lee2022convergence,chen2022sampling,chen2022improved}.

In fact,  our assumption is even weaker than $\E\{\|\hbm(\by,t)-\bm(\by,t)\|^2\} = o_N(N)$:
we only require $\|\hbm(\by,t)-\bm(\by,t)\|^2 = o_N(N)$ to hold with high probability
(as for the other assumptions.
This leads to the additional error term $\mu_{\data}\big(\|\btheta\|^2\cdot\bfone_{\|\btheta\|\ge \oR\sqrt{n}}\big)$
which however is easy to control in the applications we are interested in.
Similarly, we make no effort to improve the dependence on $T$ (appearing as $e^{C_2T}$)
because in the applications we will consider, it is safe to take $T=O(1)$.

Let us mention that guarantees in $W_2$ distance analogue to Theorem \label{thm:gen}
were proven in \cite{alaoui2022sampling}, although only stated in a special case.
Bounds in Wasserstein distance were proven concurrently to  \cite{alaoui2022sampling}
in \cite{de2022convergence}. We note however that the results of \cite{de2022convergence}
require to bound (in our notation) $\sup_{\by} \|\hbm(\by,t)-\bm(\by,t)\|/(1+\|\by\|)$
which is very difficult in the problems studied here.

%
%
\section{Main results: Spiked matrix model}\label{sec:MainSpiked}

In this section, we define the sampling algorithm and state our results for
spiked matrix models. 
We will present two theorems establishing that our approach  produces samples 
whose distribution is close to  the desired posterior
both for discrete and continuous priors.

\subsection{General setting}

Recall that the spiked matrix model is given in  Eq.~\eqref{eq:model},
with $\btheta \sim \pi_{\Theta}^{\otimes n}$. 
Without loss of generality, 
we assume $\pi_{\Theta}$ has unit second moment
$\int\theta^2\, \pi_{\Theta}(\de\theta) =1$.
We further assume that $\beta,\pi_{\Theta}$ are fixed and known.
If $\beta> \beta_{\swr}=1$
 is unknown, then $\beta$ can be consistently estimated by applying a non-linear function to
 the top eigenvalue of 
     $\bX$ \cite{baik2005phase}. 
On the other hand, if $\beta \le \beta_{\swr}$, then 
   either the posterior is close ---in a suitable sense--- to the prior
    (if $\beta\leq \beta_{\sIT}$),
   or estimation (hence sampling) is conjectured to be hard\footnote{See discussion in 
   Section \ref{sec:IntroModels}. } (for 
   $\beta_{\sIT}<\beta\leq \beta_{\swr}$).
  In addition, it is possible to consistently estimate whether $\beta > 1$ by
  using the top eigenvalue of $\bX$.  
  
  Throughout our study of this model, we will set $\bH = \id$, and hence we will not
   distinguish between $\bm$ and $\bm_{\btheta}$.

We describe two implementations of the general sampling 
algorithm as stated in \cref{alg:general}. 
This amounts to two different ways of constructing  $\hbm(\by, t)$:
\begin{enumerate}
\item For the case of a discrete prior $\pi_{\Theta}$, we construct $\hbm(\by, t)$
using a simple Bayes AMP
algorithm. The pseudocode for the overall procedure is given 
in \cref{alg:Spiked-Sampling-AMP}
\item For the case of a continuous prior $\pi_{\Theta}$, 
we construct $\hbm(\by, t)$ by implementing Bayes AMP
followed by local minimization of a suitable free energy functional to construct $\hbm(\by, t)$.
Pseudocode for this procedure is given in \cref{alg:Spiked-Sampling-continuous}. 
\end{enumerate}
We believe that the distinction between these two constructions is a proof artifact.
We introduce the free energy minimization step in the second construction in order to
show that $\by\mapsto \hbm(\by,t)$ is approximately Lipschitz continuous, but it should be possible to apply 
the first procedure to both cases.

We next present the Bayes AMP algorithm with spectral initialization, 
originally introduced in \cite{montanari2021estimation}.
 This algorithm is used for both of our constructions.
We begin by defining a scalar denoiser $\sF(\,\cdot\, ;\gamma):\reals\to\reals$
via
\begin{align}
\sF(z;\gamma):= \E[\Theta|\gamma\Theta+\sqrt{\gamma} G = z]\, , \label{eq:F}
\end{align}
where $(\Theta,G)\sim\pi_{\Theta}\otimes \normal(0,1)$. 
The function $\sF(z;\gamma)$ can be efficiently evaluated numerically, by computing a one-dimensional 
integral (or a sum if $\pi_{\Theta}$
is discrete).

For $\bz\in\reals^n$, we denote by $\sF(\bz;\gamma)$ the entrywise action
of $\sF(\,\cdot\,;\gamma)$, i.e., 
$$\sF(\bz;\gamma)=\big(\sF(z_1;\gamma),\dots,\sF(z_n;\gamma)\big)^{\sT} \in \RR^n.$$
Given inputs $(\bX, \by(t))$, we compute  
$\bv_1(\bX) \in \RR^n$ that is a top eigenvector of $\bX$.
The Bayes AMP estimate at time $t$ is then computed
recursively as follows:
\begin{align}\label{eq:general-AMP}
\begin{split}
	& \hat{\bm}_t^k = \sF(\bz_t^k;\gamma^k_t)\, , \;\;\;\bz_t^0 = \sqrt{n\beta^2(\beta^2 - 1)}\, \cdot\bv_1(\bX) =:\bnu\, ,\\
	& \bz_t^{k + 1} = \beta \bX \hat{\bm}_t^k + \by(t) - b_t^k \hat{\bm}_t^{k - 1}\, ,
\end{split}
\end{align}
Here, the sequence $\{\gamma^k_t\}_{k \geq 0}$ is defined by the state evolution recursion
%
\begin{align}\label{eq:AMPSE}
	\gamma_t^{k + 1} = \beta^2(1-\mmse(\gamma_t^k)\big)+ t ,  \qquad \;\;\;\; \gamma_t^0 = \beta^2 - 1,
\end{align}
and the memory (Onsager) coefficient is given by
$b_t^k = \beta^2 \mmse(\gamma_t^k)$, where
\begin{align}\label{eq:MMSE-DEF}
	\mmse(\gamma) := \E[(\Theta - \E[\Theta \mid \gamma\Theta + \gamma^{1/2} G])^2]\, .
\end{align}
Earlier work characterizes the optimality of Bayes AMP when $t=0$,
see \cite{montanari2021estimation} and references therein.
In order to recall the main result of this line of work, we introduce the following free energy functional:
\begin{align}\label{eq:Phi}
	\Phi(\gamma, \beta, t):=  \frac{1}{4\beta^2} (\gamma - t)^2- \frac{1}{2} (\gamma - t)+ \mathsf{I}(\gamma),
\end{align}
where $\mathsf{I}(\gamma) := I(\Theta;Y)$ is the mutual information
between $\Theta$ and $Y = \sqrt{\gamma} \Theta + G$ when
 $(\Theta,G)\sim\pi_{\Theta}\otimes \normal(0,1)$.
Explicitly, $\mathsf{I}(\gamma)  = \E \log \frac{\dd p_{Y \mid \Theta}}{\dd p_Y}(Y, \Theta)$.

In the next theorem, we characterize the optimality condition of Bayes AMP. 
The theorem is an analogous to  \cite[Corollary 2.3]{montanari2021estimation}, and we refer to
Appendix \ref{sec:AMP-optimal} for its proof.
\begin{theorem}\label{cond:first-stationary-point}
Assume that $\pi_{\Theta}$ has finite fourth moment,
$\beta>1$ and  the first stationary point of $\gamma \mapsto \Phi(\gamma, \beta, t)$  
on $(t,\infty)$ is also the unique global minimum of the same function over 
$\gamma\in (t, \infty)$. 
\textcolor{black}{We also assume that the sign of the spectral initialization $\bnu$ is
 chosen following Lines 2-4 in \cref{alg:Spiked-Sampling-AMP}.}
Then Bayes AMP asymptotically computes the 
posterior expectation, in the sense that \textcolor{black}{for all $t > 0$,} 
\begin{align}
\lim_{k\to\infty}\lim_{n\to\infty}\frac{1}{n}  \E\big[\|\hat{\bm}^k_t(\by,t)- \bm(\by,t)\|_2^2\big] = 0\, ,
\end{align}
where the expectation is taken over $\by \overset{d}{=} t \btheta + \sqrt{t} \bg$ for $(\btheta, \bg) \sim \pi_{\Theta}^{\otimes n} \otimes \normal(\mathbf{0}_n, \id_n)$. 


Further, if $\pi_{\Theta}$ has finite eighth moment, then
the above condition is always satisfied for all $\beta$ large enough.
We denote by $\beta_{\salg}(\pi_{\Theta})$ the smallest threshold such that the condition holds for all 
$\beta\in (\beta_{\salg}(\pi_{\Theta}),\infty)$.
\end{theorem}
It has been conjectured that if Bayes AMP does not reach 
the Bayes optimal estimation error, then no polynomial time algorithm 
does. The paper \cite{montanari2022equivalence} provides some rigorous support for this conjecture. 
Of course, if Bayes optimal estimation is impossible, so is sampling. 
%
%
\subsection{Discrete $\pi_{\Theta}$}
\label{sec:discrete-spiked-sampling}

We first investigate the case of $\pi_{\Theta}$ supported on finitely many points. 
We present our sampling algorithm as \cref{alg:Spiked-Sampling-AMP}. 
\begin{algorithm}
\caption{Posterior sampling under the spiked matrix model and a discrete prior}\label{alg:Spiked-Sampling-AMP}
\textbf{Input: }Data $\bX$, parameters $(\beta, K_{\AMP}, L, \Delta)$;	
\begin{algorithmic}[1]
\State Set $\hat{\by}_0 = \bfzero_n$;
\State Compute $\bnu$, a leading eigenvector of $\bX$, with norm $\|\bnu\|_2^2=n\beta^2(\beta^2 - 1)$;
\If{$\pi_{\Theta}$ is not symmetric}
        \State $\bnu \leftarrow \cA(\bnu)\bnu$;
\EndIf
	\For{$\ell = 0,1,\cdots, L - 1$}
		\State Draw $\bw_{\ell + 1} \sim \normal(0,\id_n)$ independent of everything so far;
		\State Let $\hat \bm(\hat \by_{\ell}, \ell \Delta)$ be the output of the Bayes AMP algorithm \eqref{eq:general-AMP} with $K_{\AMP}$ iterations;
		\State Update $\hat{\by}_{\ell + 1} = \hat{\by}_\ell + \hat{\bm}( \hat{\by}_\ell, \ell\Delta) \Delta + \sqrt{\Delta} \bw_{\ell + 1}$;	
	\EndFor
\If{$\pi_{\Theta}$ is symmetric}
        \State Set $s\sim\Unif(\{+1,-1\})$;
    \Else
       \State Set $s=+1$;
    \EndIf
\State \Return $\btheta^{\salg}  = s\cdot \hat{\bm}(\hat{\by}_L, L\Delta)$;
\end{algorithmic}
\end{algorithm}

Observe that \cref{alg:Spiked-Sampling-AMP} proceeds in slightly different ways depending on whether
$\pi_{\Theta}$ is symmetric around $0$ or not. If it is symmetric, then the
posterior $\mu_{\bX}$ is also symmetric under reflection $\btheta\mapsto -\btheta$:
we take account of this by explicitly symmetrizing the output samples at the end of the algorithm.  
 On the other hand, if $\pi_{\Theta}$ is not symmetric, then we need to design a separate algorithm $\cA$ that  
  aligns the initial spectral estimate with $\btheta$. 
We present theoretical guarantee for $\cA$ in \cref{lemma:sign} below, delaying the description of $\cA$ to Appendix \ref{sec:proof-of-lemma:sign}. 
We postpone the proof of \cref{lemma:sign} to Appendix \ref{sec:proof-of-lemma:sign} as well.
\begin{lemma}\label{lemma:sign}
	We assume that $\pi_{\Theta}$ has unit second moment (for this lemma, $\pi_{\Theta}$ is not necessarily discrete).
	If further $\pi_{\Theta}$ is not symmetric about the origin, then for any $\beta > 1$, 
	there exists an algorithm $\cA: \RR^n \to \{+1 ,-1\}$ with complexity
	$O(n)$, such that 
	\begin{align*}
	\lim_{n\to\infty}\prob\big(\cA(\bnu) =\sign(\< \btheta, \bnu \>)\big) = 1\, . 
	\end{align*}
\end{lemma} 
\begin{remark}
	\cref{lemma:sign} and the associated algorithm apply to
	 general $\pi_{\Theta}$ with unit second moment. For this lemma, we do not need to
	  assume  $\pi_{\Theta}$ to be either discrete or continuous. 
\end{remark}
  
We next state our theoretical guarantee for \cref{alg:Spiked-Sampling-AMP}.
In the theorem statement, we recall that $\mu_{\bX}^{\salg}:=\Law(\btheta^{\salg})$ is the distribution of the algorithm output 
$\btheta^{\salg}$ by \cref{alg:Spiked-Sampling-AMP}.
\begin{theorem}\label{thm:main2}
Assume that $\pi_{\Theta}$ is supported on finitely many points. 
Then there exists a constant $\beta_0(\pi_{\Theta})$ depending only 
on $\pi_{\Theta}$, 	such that for all $\beta \geq \beta_0(\pi_{\Theta})$, the following statement holds:
For any $\xi > 0$, there exist $K_{\AMP} , L\in \NN_+$ and $ \Delta \in \RR_{>0}$ 
that depend uniquely on $(\beta, \xi, \pi_{\Theta})$, such that if 
\cref{alg:Spiked-Sampling-AMP} takes as inputs $(\bX,\beta, K_{\AMP}, L, \Delta)$, then,
with probability $1-o_n(1)$ with respect to the choice of $\bX$, we have
	\begin{align*}
	W_{2}(\mu_{\bX}, \mu_{\bX}^{\salg})\le \xi\sqrt{n}\, .
	\end{align*} 	
\end{theorem}
The proof of this theorem is given in Appendix \ref{sec:appendix-spiked}.

\begin{remark}
\cref{alg:Spiked-Sampling-AMP} may return a vector
$\btheta^{\salg}$ that is outside the support of the prior. This does not contradict
Theorem \ref{thm:main2}, which guarantees that $\btheta^{\salg}$ is close
in $\ell_2$ distance to a sample $\btheta\sim \mu_{\bX}$ but does not necessarily lie inside the support.

If necessary, this can be remedied by a simple rounding procedure.
More precisely,
we can replace the last line of \cref{alg:Spiked-Sampling-AMP} by the following procedure:
Compute $\bm^{\salg}:= s\cdot \hat{\bm}(\hat{\by}_L, L\Delta)$, 
and for each $i\in [n]$, let $\underline\theta^{\salg}_i=\sup(x\in\supp(\pi_{\Theta}):\, x\le m_{i}^{\salg})$
and $\overline\theta^{\salg}_i=\inf(x\in\supp(\pi_{\Theta}):\, x>m_{i}^{\salg})$.
We them  return a vector with 
conditionally independent coordinates given $\bm^{\salg}$, such that, for each $i \in [n]$, 
$\theta^{\salg}_i\in \{\underline\theta^{\salg}_i,\overline\theta^{\salg}_i\}$, and $\E[\theta^{\salg}_i|\bm^{\salg}]=m_i^{\salg}$.
 This modification enjoys the same $W_2$ guarantees as the original one.
\end{remark}
%
%

\subsection{Continuous $\pi_{\Theta}$}

In this section, we discuss the case of a continuous prior $\pi_{\Theta}$. 

\subsubsection{Free energy minimization}
\label{sec:Free_En_Spiked}

As mentioned before, for the case of a continuous prior $\pi_{\Theta}$, we  modify our construction
of the approximate posterior mean $\hbm(\by,t)$. 
Specifically, denoting
by $\hbm_{t}^{K_{\AMP}}(\by)$ the output of the Bayes AMP algorithm \eqref{eq:general-AMP} at time $t$
with $K_{\AMP}$ iterations and side information vector $\by$, we run a gradient descent algorithm 
to locally minimize a modified version of the TAP free energy, with initialization
$\hbm_{t}^{K_{\AMP}}(\by)$.


We begin by introducing the free energy functional that we use. 
Let $\gamma_{\beta, t} := \argmin_{\gamma > 0} \Phi(\gamma, \beta, t)$, where we 
recall that $\Phi$ is defined in \cref{eq:Phi}. 
We further define 
\begin{align}
q_{\beta, t} := 1-\mmse(\gamma_{\beta,t})\, ,\label{eq:Qdef}
\end{align}
 where we recall that $\mmse$ is 
defined in \cref{eq:MMSE-DEF}. 
Whenever $t>0$ and the global minimizer of $\gamma \mapsto \Phi(\gamma, \beta, t)$ is unique,
this quantity has the following interpretation (see \cite{montanari2021estimation}):
\begin{align}
	q_{\beta, t} = \lim_{n \to \infty} \frac{1}{n}\E\left[ \|\E[\btheta \mid \bX, \by(t)]\|_2^2 \right]\, .
	\label{eq:LimitQ}
\end{align}
We then define the free energy functional using $q_{\beta, t}$. 
For $\bX\in\RR^{n\times n}$, $\by\in\RR^n$, $\beta,t\in\RR_{> 0}$, we define 
$\cuF_{\sTAP}(\,\cdot\,; \bX, \by, \beta, t) :\RR^n\to\RR$ via
\begin{align*}
	\cuF_{\sTAP}(\bm; \bX, \by, \beta, t) =  - \frac{\beta}{2} \langle \bm, \bX \bm \rangle + \frac{1}{2}\beta^2 (1 - q_{\beta, t})\|\bm\|_2^2 + \sum_{i = 1}^n h(m_i, \beta^2 q_{\beta, t}) - \langle \by, \bm \rangle,  
\end{align*} 
where
\begin{align*}
	 h(m, w) := \sup_{\lambda} \left[ \lambda m - \phi(\lambda, w) \right], \qquad \phi(\lambda, w) = \log \Big\{\int e^{\lambda \theta - w\theta^2 / 2} \pi_{\Theta}(\dd \theta)\Big\}
	 \, . 
\end{align*}
We denote by $\Spher^{n-1}_{\beta, t} := \Spher^{n-1}(\sqrt{n q_{\beta, t}})$,
the sphere  of radius $\sqrt{n q_{\beta, t}}$ in $\RR^n$. 
Motivated by \cref{eq:LimitQ}, we propose to minimize $\cuF_{\sTAP}(\bm; \bX, \by, \beta, t)$
 under two constraints: $\bm$ is on the sphere $\Spher^{n-1}_{\beta, t}$ and inside a 
 neighborhood of $\hbm_{t}^{K_{\AMP}}(\by)$.

  In order to implement this idea, 
  for $\bm \in \Spher^{n-1}_{\beta, t}$ we consider the system of coordinates 
  defined on $\{\bx\in \Spher^{n-1}_{\beta, t}: \; \<\bx,\bm\> >0\}$ by projection on the tangent plane
  from the center of the sphere. 
   Explicitly, let $\bT_{\bm} \in \RR^{n \times (n - 1)}$
  be a matrix with orthonormal columns such that $\bm^{\sT} \bT_{\bm} = \mathbf{0}$ and 
  define $\bphi_{\bm}:\RR^{n-1}\to\Spher^{n-1}_{\beta,t}$ via
\begin{align}
\label{eq:varphi-project}
	\bphi_{\bm}(\bw) =\sqrt{n q_{\beta, t}} \cdot \frac{\bm + \bT_{\bm} \bw}{\|\bm + \bT_{\bm} \bw \|_2}. 
\end{align}
To run gradient descent, we define $\hbm^{\AMP}:= \sqrt{n q_{\beta,t}}\hbm_{t}^{K_{\AMP}}(\by)/\|\hbm_{t}^{K_{\AMP}}(\by)\|_2$, and implement the following updates: 
\begin{align}
\label{eq:GD-tangent-space}
\begin{split}
	& \bw^0 = \mathbf{0}_{n - 1}, \\
	& \bw^{i + 1} = \bw^i + \zeta \cdot \nabla_{\bw} \cuF_{\TAP}(\bphi_{\hbm^{\AMP}}(\bw); \bX, \by, \beta,t) \big|_{\bw = \bw^i},
	 \qquad i = 0, 1, \cdots, K_{\GD} - 1\, ,  
\end{split}
\end{align}
with $K_{\GD}$ being the number of steps and $\zeta$ being the step size. 

The final approximation of the posterior mean $\hbm(\by,t)$ is given by $\bphi_{\hbm^{\AMP}}(\bw^{K_{\GD}})$.
The procedure stated here is used to construct $\hat\bm(\by, t)$ in \cref{alg:Spiked-Sampling-continuous}. 

\subsubsection{Sampling guarantees}

Next, we describe an algorithm that
yields approximately correct samples for a a broad class of priors $\pi_{\Theta}$ 
that have a probability density function. 
Without loss of generality, we write 
\begin{align}
\label{eq:pi-U-density}
	\pi_{\Theta}(\dd \theta) \propto e^{-U(\theta)} \dd \theta\, .
\end{align}
Our key assumption will be  that $U: \RR \to \RR$ has bounded second order derivative 
$\|U''\|_{\infty}:=\sup_{x\in\reals}|U''(x)| <\infty$. 

We state our algorithm as \cref{alg:Spiked-Sampling-continuous}.  
%
\begin{algorithm}
\caption{Posterior sampling under the spiked matrix model (continuous non-symmetric prior)}\label{alg:Spiked-Sampling-continuous}
\textbf{Input: }Data $\bX$, parameters $(\beta, K_{\AMP}, K_{\GD}, \zeta, L, \Delta)$;	
\begin{algorithmic}[1]
\State Set $\hat{\by}_0 = \bfzero_n$;
\State Compute $\bnu$, a leading eigenvector of $\bX$;
\State Normalize $\bnu$ such that $\|\bnu\|_2^2=n\beta^2(\beta^2 - 1)$;
\If{$\pi_{\Theta}$ is not symmetric}
        \State $\bnu \leftarrow \cA(\bnu)\bnu$;
    \EndIf
	\For{$\ell = 0,1,\cdots, L - 1$}
		\State Draw $\bw_{\ell + 1} \sim \normal(0,\id_n)$ independent of everything so far;
		\State Let $\hat \bm_{\ell \Delta}^{K_{\AMP}}(\hat \by_{\ell})$ be the output of the Bayes AMP algorithm \eqref{eq:general-AMP} with $K_{\AMP}$ iterations;
		\State Let $\hat\bm(\hat \by_{\ell},\ell\Delta)$ be the output of algorithm \eqref{eq:GD-tangent-space} with AMP output $\hat \bm_{\ell \Delta}^{K_{\AMP}}(\hat \by_{\ell})$, number of gradient descent steps $K_{\GD}$ and step size $\zeta$;
		\State Update $\hat{\by}_{\ell + 1} = \hat{\by}_\ell + \hat{\bm}( \hat{\by}_\ell, \ell\Delta) \delta + \sqrt{\Delta} \bw_{\ell + 1}$;	
	\EndFor
\If{$\pi_{\Theta}$ is symmetric}
        \State Set $s\sim\Unif(\{+1,-1\})$;
    \Else
       \State Set $s=+1$;
    \EndIf
\State \Return $\btheta^{\salg}  = s \cdot \hat{\bm}(\hat{\by}_L, L\Delta)$;
\end{algorithmic}
\end{algorithm}  

We recall that  $\cA$ is the algorithm of
 \cref{lemma:sign}, which works for both discrete and continuous priors. We 
 present our theoretical guarantee for \cref{alg:Spiked-Sampling-continuous} in
  the next theorem, deferring its proof  to Appendix \ref{sec:spiked-continious}.
\begin{theorem}
\label{thm:spiked-continuous}
 Assume $\pi_{\Theta}$ has form \eqref{eq:pi-U-density} with
 $\|U''\|_{\infty}<\infty$. 
Further, assume $\pi_{\Theta}$ has unit second moment and finite eighth moment.
Then there exists a constant $\beta_0(\pi_{\Theta})$ depending only 
on $\pi_{\Theta}$, 	such that for all $\beta \geq \beta_0(\pi_{\Theta})$, the following statement holds.
For any $\xi > 0$, there exist $K_{\AMP}, K_{\GD}, L\in \NN_+$ and $ \Delta, \zeta \in \RR_{>0}$ 
that depend uniquely on $(\beta, \xi, \pi_{\Theta})$, such that if 
\cref{alg:Spiked-Sampling-continuous} takes as inputs $(\bX,\beta, K_{\AMP}, K_{\GD}, \zeta, L, \Delta)$, then,
with probability $1-o_n(1)$ with respect to the choice of $\bX$, we have
	\begin{align*}
			W_{2} (\mu_{\bX}, \mu_{\bX}^{\salg}) \leq \xi \sqrt{n}. 
		\end{align*}  
\end{theorem} 
%
%
\section{Main results: Linear model}
\label{sec:MainLinear}

In this section, 
we define the sampling algorithm and state our results for high-dimensional linear regression.
An important difference with respect to the previous section is that here we 
use\footnote{Recall that $\bH$ is the matrix appears in process \eqref{eq:GeneralLinearObs}.}
$\bH=\bX$. 

We will describe the construction of the posterior mean estimator and state a theorem 
demonstrating that our algorithm approximately samples from the target posteriors,
provided that $n/p$ exceeds a constant dependent on $\pi_{\Theta}$.
Because of the choice $\bH=\bX$,  
the stochastic localization process is equivalent to a process along which the noise level in the 
response variables decreases with $t$. As a consequence, it will be sufficient to construct the 
posterior mean estimator at $t=0$.

In Appendix \ref{sec:low-SNR} we will also state a theorem for the case of a sufficiently small $n/p$.
We briefly discuss the differences in Remark \ref{rmk:n_small} below.

\subsection{General setting}
\label{sec:general-setting}

As announced, the response variables $\by_0$ 
(cf. Eq.~\eqref{eq:LR}) and stochastic localization process
$\bby(t)$ take the form
\begin{align}
	\by_0 &= \bX \btheta + \bepsilon\, ,\;\;\;\;\; \bepsilon \sim\normal(0,\sigma^2\id_n)\, ,\\
	\bby(t) &= t\bX \btheta + \bG(t)\, .
\end{align}
where we recall that $\btheta \sim \pi_{\Theta}^{\otimes p}$
is independent of $\bepsilon\sim \normal(\bzero,\sigma^2\id_n)$ and the standard 
Brownian motion $\bG(t)$.
We can conveniently aggregate $\bby(t)$ and $\by_0 $ using the sufficient statistics 
\begin{align}
\by(\sigma^{-2}+t) := \frac{1}{\sigma^2}\by_0 +\bby(t)\, .\label{eq:LRSsufficientStatistics}
\end{align}
One can verify that this takes the form $\by(t)= t\bX\btheta + \bG'(t)$, 
where $t\ge \sigma^{-2}$ and $\bG'(t)$ is a standard Brownian motion in $\RR^n$. 
In particular $(\sigma^{-2}+t)^{-1}\by(\sigma^{-2}+t)$ is distributed as $\by_0$, except that
the noise level $\sigma^2$ is replaced by $(\sigma^{-2}+t)^{-1}$.

Therefore, to construct the posterior mean at a general time $t$,
it is sufficient to consider the  $(\bX, \by_0)$ distributed according to the linear model \eqref{eq:LR},
which we will do in this section and the next one. 

Our main assumption is that the design matrix $\bX$ is random and has i.i.d. Gaussian entries.
We emphasize that (unlike $\btheta$ and $\bepsilon$), $\bX$ and $\by_0$ are given to
 the practioners as observations. 
\begin{assumption}
\label{assumption:LR}
	We assume $\bX$ has i.i.d. rows $(\bx_i)_{i \in [n]}$, with
	\begin{align*}
		\bx_i \sim \normal(0, \id_p / p)\, . 
	\end{align*}
	 As $n, p \to \infty$, we assume
	 $n / p \to \delta \in (0, \infty)$.
\end{assumption}
\begin{remark}
Assumption \ref{assumption:LR} is introduced to simplify the analysis of our algorithm.
We believe that the proposed sampling method can be generalized to a broader class of design matrices:
$(1)$~We expect that matrices  with independent, centered entries that have equal variance
will, under suitable moment conditions, fall into the same universality class as the
present work, and hence the same algorithm will apply 
\cite{bayati2015universality,barbier2019optimal,chen2021universality}; 
$(2)$~It should be possible to extend the present algorithm to 
matrices with unitarily invariant distribution or i.i.d. non-isotropic rows using tools from 
\cite{berthier2020state,fan2022approximate,dudeja2022spectral}.
\end{remark}

We will also make the following assumptions for technical reasons.
\begin{assumption}
\label{assumption:three-points}
	We assume $\pi_{\Theta}$ is supported on at least three distinct points. 
\end{assumption}

%
   
We next state the AMP algorithm for model \eqref{eq:LR}, which in many regimes effectively 
estimates the posterior expectation $\bm_{\btheta}(\bX,\by_0, \sigma^2) := \E[\btheta \mid \bX, \by_0]$.  
The general AMP iteration reads
\begin{align}
\label{eq:linear-AMP}
	\left\{
	\begin{array}{l}
		\bb^{k + 1} = \bX^{\top} f_k(\ba^k, \by_0) - \xi_k g_k(\bb^k), \\
		\ba^k = \bX g_k(\bb^k) - \eta_k f_{k - 1}(\ba^{k - 1}, \by_0), 
	\end{array} \right.
\end{align}
where at initialization we set $\ba^{-1} = \mathbf{0}_n$, and 
the scalar sequences $(\xi_k)_{k \geq 0}$ and $(\eta_k)_{k \geq 0}$ for the case of Bayes AMP are 
defined in \cref{eq:xik-etak} below.
After $k$ iterations, algorithm \eqref{eq:linear-AMP} outputs the following estimate of 
$\E[\btheta|\bX,\by_0]$:
\begin{align}
\hbm^k_{\btheta} = \hbm^k_{\btheta} (\bX,\by_0, \sigma^2) = g_k(\bb^k) \in \RR^p\, .  
\end{align}
In the above iteration, we assume $f_k: \RR^2 \to \RR$, $g_k: \RR \to \RR$ are Lipschitz continuous and act on
matrices/vectors row-wise. 

The Bayes AMP algorithm is defined by the following optimal choice of the nonlinearities: 
\begin{align}
f_k(x, y) = \frac{\sigma^2(y - x)}{\sigma^2 + E_k},  \;\;\;\;\;\;
 g_k (x) = \E[\Theta \mid \gamma_k \Theta + \sqrt{\gamma_k} G = \sigma^{-2}x]\, .
\end{align}	
Here, expectations are with respect to
 $(\Theta, G) \sim \pi_{\Theta} \otimes \normal(0,1)$ and 
 the values of $\{E_k\}_{k\ge 0}$ are determined by the following
state evolution recursion, with initialization $E_{-1}=1$:
\begin{align*}
	E_{k}= \mmse \Big( { \delta (\sigma^2 + E_{k - 1})^{-1} }  \Big)\, ,
\end{align*}
and $\gamma_k = \delta /(\sigma^2 + E_{k-1})$.
We recall that the function $\mmse:\RR_{\ge 0}\to\RR_{\ge 0}$ is defined in Eq.~\eqref{eq:MMSE-DEF}.
Finally the memory (Onsager) coefficients  $\xi_k,\eta_k$ are given as follows:
\begin{align}
\label{eq:xik-etak}
	& \xi_k = -\frac{\delta\sigma^2}{\sigma^2 + E_k}, \qquad \eta_k = \sigma^{-2}\E\left[ \Var[\Theta \mid \gamma_k  \Theta + \sqrt{\gamma_k} G] \right]. 
\end{align}
The optimality of Bayes AMP can be characterized by using the following free energy 
functional
\begin{align}
\Phi(\gamma,\sigma^2,\delta)  = \frac{\sigma^2 \gamma}{2} - \frac{\delta}{2} \log \frac{\gamma}{2\pi \delta} + \Info(\gamma), 
\end{align}
We copy the AMP consistency result from \cite{barbier2019optimal} below for readers' convenience. 
\begin{theorem}[\cite{barbier2019optimal}]\label{thm:AMP-LinReg}
Assume that the first stationary point of $\gamma\mapsto \Phi(\gamma,\sigma^2,\delta)$ on
$\gamma\in (0,\infty)$  is also the unique global minimum of the same function.
Then, we have
\begin{align}
\label{eq:lemma5.1-eq}
	\lim_{k\to\infty}\plimsup_{n, p \to \infty}  \frac{1}{{p}} \EE\big[ \|\bm_{\btheta} (\bX,\by_0, \sigma^2) - 
	\hbm_{\btheta}^k(\bX,\by_0, \sigma^2) \|_2^2 \big] = 0\, .
\end{align}
\end{theorem}
The following lemma is a consequence of \cite{barbier2019optimal}. 
We refer to Appendix 
\ref{sec:proof-lemma:PosteriorAMP_LR} for a proof.
\begin{lemma}\label{lemma:PosteriorAMP_LR}
If $\pi_{\Theta}$ has unit second moment and bounded eight moment, 
then  there exists $\odelta_{\salg}= \odelta_{\salg}(\pi_{\Theta})<\infty$ 
that depends only on $\pi_{\Theta}$, 
such that the above condition 
holds for all $\sigma^2$ and all $\delta)\in (\odelta_{\salg},\infty)$ 
We will denote by $\odelta_{\salg}$ the infimum of such threshold. 
\end{lemma} 
%
%
\subsection{Free energy and mean estimation}
\label{sec:LR_Free}

We use a modified TAP free energy to approximate the posterior mean.
The landscape of local minima of this TAP free energy for linear regression was recently studied in \cite{celentano2023mean},
and we will build upon their findings.

In order to define the free energy functional, we define the moments' set 
\begin{align*}
	\Gamma := \left\{ (m, s) \in \RR^2: \mbox{there exists }(\lambda, \gamma) \in \RR^2,
	 \mbox{ such that }m = \E_{\lambda,\gamma} [\Theta]\mbox{ and }s = \E_{\lambda,\gamma}[ \Theta^2 ] \right\}\, .
\end{align*}
where $\E_{\lambda, \gamma}$ denotes expectation with respect to the 
tilted probability measure
$$\pi_{\lambda, \gamma}(\de \theta) \propto e^{-\gamma \theta^2 / 2 + \lambda \theta}
 \pi_{\Theta}(\de \theta).$$
If $\pi_{\Theta}$ has  bounded  support, then standard theory of exponential
families \cite{brown1986fundamentals} implies that $\Gamma$ is open and convex.
We identify  $\Gamma^p$ with a subset of $\RR^p \times \RR^p$ in the obvious way and define,
 for $(\bm, \bs) \in \Gamma^p$, 
\begin{align}
\label{eq:TAP-main-no-t}
\begin{split}
	\cuF_{\sTAP}(\bm, \bs; \by_0, \bX, \sigma^2) := & \frac{n}{2} \log 2\pi \sigma^2 + D_0(\bm, \bs) + \frac{1}{2\sigma^2} \|\by_0 - \bX\bm\|_2^2 \\
	& + \frac{n}{2} \log \left( 1 + \frac{S(\bs) - Q(\bm)}{\sigma^2} \right) \, .
\end{split} 
\end{align}
We ignore the dependency of $\cuF_{\sTAP}(\cdots)$ on $(\by_0, \bX, \sigma^2)$ when there is no confusion. 
In the above display, 
\begin{align*}
	 S(\bs) &:= \frac{1}{p} \sum_{j = 1}^p s_j, \qquad Q(\bm): = \frac{1}{p} \sum_{j = 1}^p m_j^2, \qquad D_0(\bm, \bs) := \sum_{j = 1}^p -h(m_j, s_j), \\
	 -h(m, s) &= \sup_{(\lambda, \gamma) \in \RR^2}  \Big\{ -\frac{1}{2}\gamma s + \lambda m - \log 
	\int e^{-\gamma \theta^2 / 2 + \lambda \theta} \pi_{\Theta}(\de\theta)\Big\}\\
	&=\inf_{\nu}\Big\{D_{\sKL}(\nu\|\pi_{\Theta}):\; \int \theta\, \nu(\de\theta) = m,\, 
	 \int \theta^2\, \nu(\de\theta) = s\Big\}\, .
\end{align*}
Following \cite{celentano2021local}, we can find a local minimum using a two-stages algorithm
similar to the one described in the previous section, or in \cite{alaoui2022sampling,alaoui2023sampling}:
\begin{enumerate}
\item Run AMP for $K_{\sAMP}$ iterations, cf. \cref{eq:linear-AMP}, to produce an estimate $\hbm^{K_{\sAMP}}$.
\item Run natural gradient descent (NGD) for $K_{\sNGD}$ iterations with initialization
$\hbm^{K_{\sAMP}}$ and stepsize $\eta$.
\end{enumerate}
We refer to \cite[Section 3.3]{celentano2023mean} for further details on the NGD part.
We denote the estimator thus defined by $\hbm_{\btheta}(\bX,\by,\sigma^2)$,
omitting the dependence on algorithm parameters $K_{\sAMP}, K_{\sNGD}, \eta$.

%
%
\subsection{Sampling algorithm and guarantees}
\label{sec:SamplingLinear}

The sufficient statistics $\by(t)$ of Eq.~\eqref{eq:LRSsufficientStatistics}
 satisfies the following SDE for $t\ge \sigma^{-2}$: 
\begin{align}
\de \by(t) = \bX\bm_{\btheta}(\bX, \by(t) / t, 1 / t)\de t+\de\bB(t)\, ,\;\;\; \by(\sigma^{-2}) = \frac{1}{\sigma^2}\by_0\, , 
\end{align}
where $\bm_{\btheta}(\bX, \by(t) / t, 1 / t) = \EE[\btheta \mid \bX, \by(t)] = \EE[\btheta \mid \bX, \by_0, \bar\by(t - \sigma^{-2})]$.
Once again, we apply a Euler discretization to this process, as detailed in  \cref{alg:LR-high-SNR}.

\begin{algorithm}[ht]
\caption{Posterior sampling for linear models in a high SNR regime}\label{alg:LR-high-SNR}
\textbf{Input: }Data $(\bX, \by_0)$, parameters $(K_{\AMP}, K_{\mathsf{NGD}}, \eta, \delta, \sigma^2, L, \Delta)$;	
\begin{algorithmic}[1]
\State Set $\hat{\by}_0 = \by_0/\sigma^2$;
\For{$\ell = 0,1,\cdots, L - 1$}
	\State Draw $\bw_{\ell + 1} \sim \normal(0,\id_n)$ independent of everything so far;
	\State Let $t_{\ell} := \sigma^{-2}+\ell\Delta$;
	\State Denote $\hbm_{\btheta}(\bX, \,\hat\by_{\ell}/t_{\ell},\, 1/t_{\ell})$ the output of
	the estimation algorithm of Section \ref{sec:LR_Free}, with omitted parameters 
	$(K_{\AMP}, K_{\mathsf{NGD}}, \eta)$.
	\State Update $\hat{\by}_{\ell + 1} = \hat{\by}_\ell + \bX \hbm_{\btheta}(\bX, \,\hat\by_{\ell}/t_{\ell},\, 1/t_{\ell})\Delta + \sqrt{\Delta} \bw_{\ell + 1}$;	
	\EndFor
\State \Return $\btheta^{\salg}  = \hbm_{\btheta}(\bX,\, \hat\by_{L}/t_{L}, \,1/t_{L})$;
\end{algorithmic}
\end{algorithm}

 As before, we denote by $\mu_{\bX, \by_0}^{\salg}$ the law of
 $\btheta^{\salg}$, the output of \cref{alg:LR-high-SNR}, and by $\mu_{\bX, \by_0}$
 the true posterior distribution. 
 The following theorem characterizes the closeness between $\mu_{\bX, \by_0}^{\salg}$ and $\mu_{\bX, \by_0}$.
 We prove the theorem in Appendix \ref{sec:proof-thm:LR-main}.  
\begin{theorem}
\label{thm:LR-main}
Under Assumptions \ref{assumption:LR} and \ref{assumption:three-points}, further assume that 
$\pi_{\Theta}$ has bounded
support. Then there exists $\delta_0=\delta_0(\pi_{\Theta})<\infty$  that depends only on 
$\pi_{\Theta}$,
such that for any $\delta \geq \delta_0$, the following holds:
For any $\xi > 0$, there exist parameters $(K_{\AMP}, K_{\mathsf{NGD}}, \eta, L, \Delta)$ that 
depend only on $(\delta, \sigma^2, \pi_{\Theta})$, such that
with probability $1 - o_n(1)$ over the randomness of $(\bX, \by_0)$, 
		\begin{align*}
			W_{2, p}(\mu_{\bX, \by_0}, \mu_{\bX, \by_0}^{\salg}) \leq \xi. 
		\end{align*}
\end{theorem} 
%

\begin{remark}\label{rmk:n_small} 
In Appendix \ref{sec:low-SNR} we prove that a theorem analogous to 
\cref{thm:LR-main} holds when $\delta<\delta_1(\pi_{\Theta},\sigma^2)$ for a different 
(smaller) constant $\delta_1$.  

Recall that ---for a general prior $\pi_{\Theta}$--- there are intervals of $\delta$
for which Bayes optimal estimation is believed to be computationally hard (cf.
Section \ref{sec:IntroModels}). Hence, we do not expect to be able to sample efficiently
from the posterior for values of $\delta$ that fall inside these intervals.
\end{remark}

%
 
\section{The use of non-linear observations}
\label{sec:GenNonlinear}

We can further generalize the linear observation model \eqref{eq:GeneralLinearObs},
by admitting non-linear observations and non-isotropic noise.
Such an observation process takes the form:
\begin{align}
\by(t) = \int_0^t\bQ(s)\bF(\btheta, s)\, \de s + \int_{0}^t \bQ(s)^{1/2}\de\bG(s)\, .
\end{align}
Here $\bQ:\reals_{\ge 0}\to\Sym_{\ge 0}(N)$ is a function taking values in the cone of positive semidefinite 
matrices, and $\bF:\reals^n \times \reals_{\ge 0} \to\reals^N$.

It is straightforward 
to generalize the algorithm and analysis of the previous section to this case.
Instead of doing this, we discuss a specific construction that is well suited to 
the spiked model analyzed in this paper.
To be specific, we let 
\begin{align}
\obY(t) = \frac{t}{n}\btheta \btheta^{\sT} + \obG(t)\, ,\label{eq:Quadratic}
\end{align} 
where  $\obG(t)$ is a \emph{symmetric Brownian motion}, i.e., a stochastic process taking values
in $\reals^{n\times n}$, with $\obG(t)=\obG(t)^{\sT}$ and such that 
$(\oG_{ij}(t))_{1 \le i\le j\le n}$ is a collection of independent Brownian motions (independent of $\btheta$),
which are time scaled so that $\E \{\oG_{ii}(t)^2\}= 2t/n$,    $\E \{\oG_{ij}(t)^2\} =t/n$
for $i<j$. 

This observation process is  of the same nature as the original 
observation that defines the model, cf. Eq.~\eqref{eq:model}. In particular, 
the process \eqref{eq:Quadratic} does not break the symmetry $\btheta\to -\btheta$.
These remarks can be further formalized by noting that $\bY(\beta^2+t) := \beta \bX+\obY(t)$ 
is a sufficient statistics for  $\btheta$ given $\bX,\obY(t)$. Of course, 
$\{\bY(t)\}_{t \ge \beta^2}$ takes the form
\begin{align}
\bY(t) = \frac{t}{n}\btheta \btheta^{\sT} + \bG(t)\, ,
\end{align} 
where $\{\bG(t)\}_{t\ge 0}$ is again a symmetric Brownian motion, except that it is initialized 
at $\bG(\beta^2) = \beta\bW$.

As for similar observation processes derived in the previous pages, 
$\bY(t)$ satisfies an SDE, namely 
\begin{align}
\de\bY(t) &= \bM(\bY(t);t)\de t + \de\bB(t)\, ,\\
\bM(\bY;t) & :=
\E\Big\{\frac{1}{n} \btheta\btheta^{\sT}\Big|\frac{t}{n}\btheta \btheta^{\sT} + \bG(t) = \bY\Big\} \, ,
\end{align} 
since we incorporated the observation $\bX$ in $\bY(t)$, 
this SDE has to be solved with initialization at $t=\beta^2$:
\begin{align}
\bY(\beta^2)= \beta \bX\, .
\end{align}

\begin{algorithm}
\caption{Diffusion-based sampling for spiked models}\label{alg:general-sampling-2}
\textbf{Input: }Data $\bX$, parameters $(\beta, K_{\AMP}, L, \Delta)$;	
\begin{algorithmic}[1]
\State Set $\bhY_0 = \beta \bX\in\reals^{n\times n}$;
	\For{$\ell = 0,1,\cdots, L - 1$}
		\State Draw $\bW_{\ell + 1} \sim \GOE(n)$ independent of everything so far;
		\State Let $\hbM(\bhY_\ell,\ell\Delta + \beta^2)$ be the Bayes AMP estimate of $\bM(\bhY_\ell,\ell\Delta + \beta^2)$
		(see main text);
		\State Update $\bhY_{\ell + 1} = \bhY_\ell + \hbM( \bhY_\ell, \ell\Delta + \beta^2) \Delta + \sqrt{\Delta} \bW_{\ell + 1}$;	
	\EndFor
\State Compute $\bX^{\salg}  = \hbM( \bhY_L, L\Delta + \beta^2)$,
and let $\lambda_1(\bX^{\salg})$, $\bv_1(\bX^{\salg})$ be its top eigenvalue/eigenvector;
   \If { $\pi_\Theta$ is symmetric}\;
   {        
 Draw      $s\sim\Unif(\{+1,-1\})$\;
 }
        \Else \;
    {
    Compute    $s=\tcA(\bv_1(\bX^{\salg}))$\;
        }
        \EndIf\;
\State \Return $\btheta^{\salg}  = s\sqrt{n\lambda_1(\bX^{\salg})}\bv_1(\bX^{\salg})$; 
\end{algorithmic}
\end{algorithm}

The resulting sampling procedure is outlined as \cref{alg:general-sampling-2}.
 Two components are unspecified: $(i)$~An algorithm $\tcA:\reals^n\to\{+1,-1\}$ 
 such that, with high probability, $\tcA(\bv_1) = \sign \<\bv_1,\btheta\>$
 when $\bv_1=\bv_1(\bX^{\salg})$;
$(ii)$~An algorithm to compute an approximation  
 $\hbM( \bY; t)$ for the conditional expectation $\bM(\bY;t)$. 
 The first algorithm is completely analogous to $\cA$ introduced in Lemma \ref{lemma:sign}.
 
 Finally, in order to compute an approximation of 
 $\bM(\bY;t)$, we use once more the AMP algorithm, whereby we replace $\bX$ by $\bY$:
\begin{align}
 & \bz_t^0 = \bnu_t\, , \\
	& \hat{\bm}_t^k = \sF(\bz_t^k; \tilde\alpha^k_t)\, , \\
	& \bz_t^{k + 1} = \bY \hat{\bm}_t^k- \tilde b_t^k \hat{\bm}_t^{k - 1}\, ,
\end{align}
where $\bnu_t$ is a randomly selected top eigenvector of $\bY$, normalized such that $\|\bnu_t\|_2 = \sqrt{nt(t - 1)}$, 
\begin{align} 
	& \tilde{\alpha}_t^0 = t - 1, \\
	& \tilde{\alpha}_t^{k + 1} = t \E[\E[\Theta \mid \tilde{\alpha}_t^{k} \Theta + (\tilde{\alpha}_t^{k})^{1/2}G]^2], \\
	& \tilde{b}_t^k = t \E[(\Theta - \E[\Theta \mid \tilde{\alpha}_t^{k} \Theta + (\tilde{\alpha}_t^{k})^{1/2}G])^2].
\end{align}
We then let $\hbm(\bY;t)= \hbm^{K_{\AMP}}_t$ and $\hbM(\bY;t) = \hbm(\bY;t)\hbm(\bY;t)^{\sT}$.

\section{Numerical experiments}
\label{sec:Numerical}

In this section we present numerical experiments in which we use
\cref{alg:Spiked-Sampling-AMP} to sample from the Bayes posterior of the 
spiked model. In  these experiments we focus on the $\ZZ_2$-synchronization task,
i.e., let $\pi_{\Theta} = (\delta_{+1}+\delta_{-1})/2$.  
Since the prior is discrete, several standard techniques  (e.g., Langevin or Hamiltonian Monte Carlo)
do not apply. 
Additionally, there are no guarantees for Gibbs sampling, also known as the 
``Glauber dynamics'' in this context.
We recall that, in this case, the algorithmic threshold for Bayes-optimal 
estimation coincides with the information-theoretic one
$\beta_{\salg}(\pi_{\Theta})=\beta_{\sIT}(\pi_{\Theta}) = 1$.
%

\begin{figure}[ht]
     \centering
     \begin{subfigure}{\textwidth}
         \centering       \includegraphics[width=\textwidth]{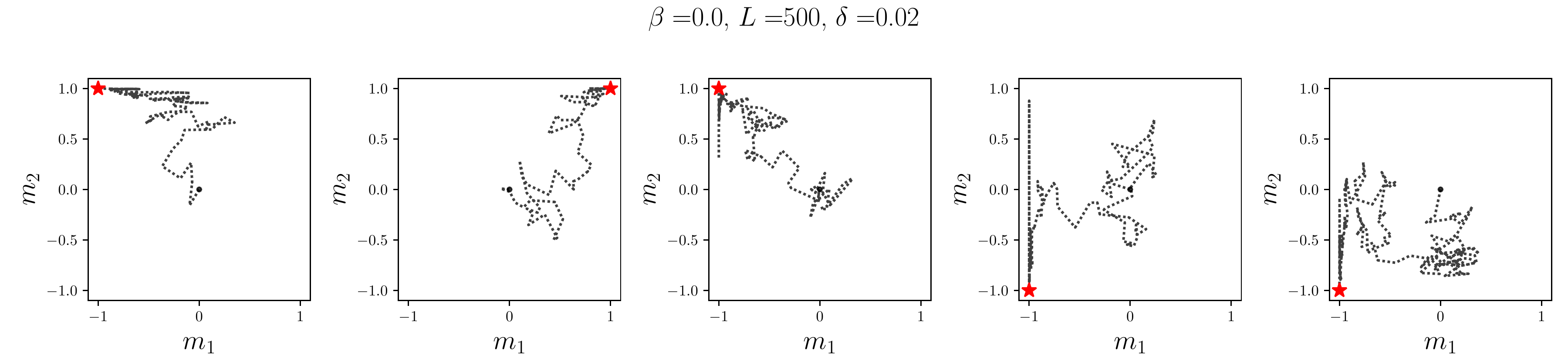}
     \end{subfigure} \\
     \vfill
     \begin{subfigure}{\textwidth}
         \centering      \includegraphics[width=\textwidth]{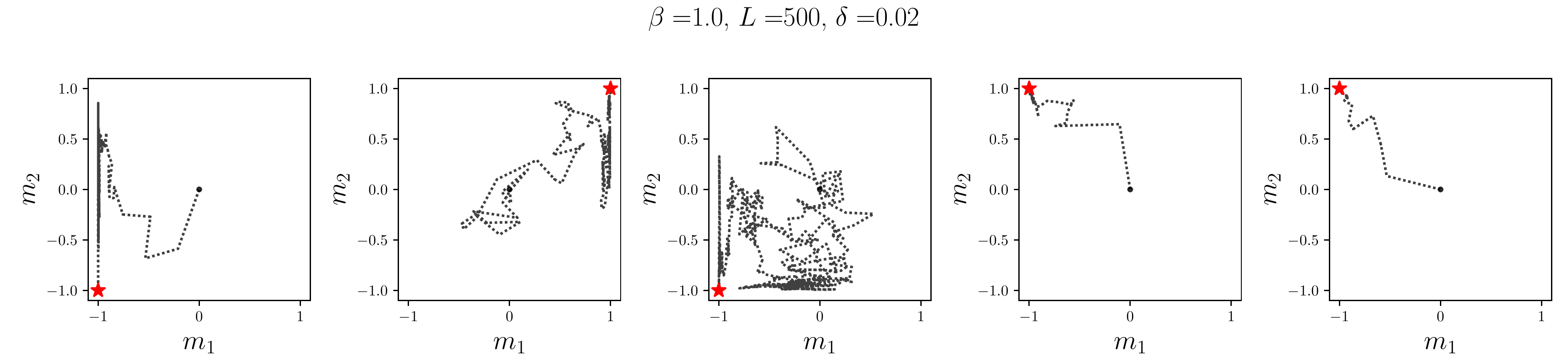}
     \end{subfigure} \\
     \vfill
     \begin{subfigure}{\textwidth}
         \centering        \includegraphics[width=\textwidth]{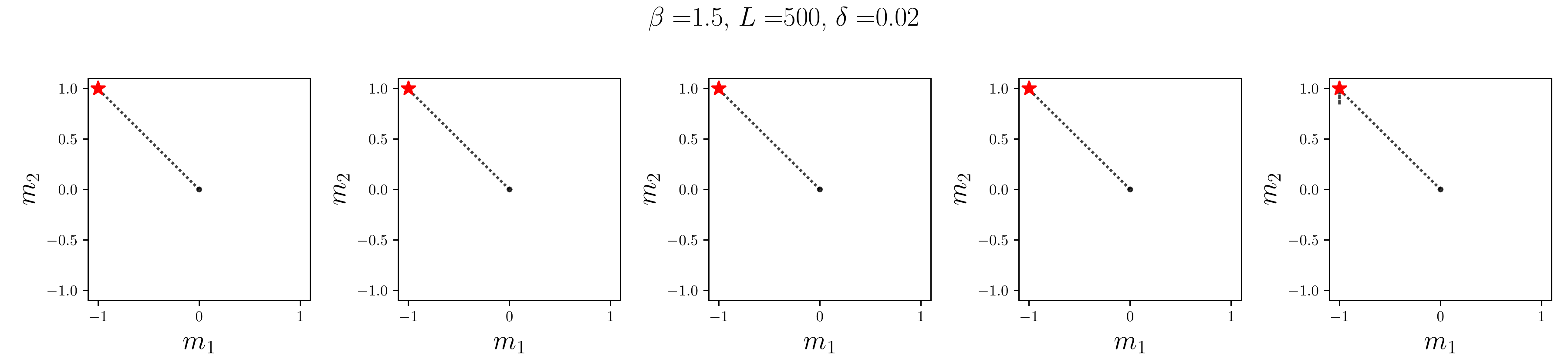}
     \end{subfigure}
        \caption{Trajectories of the first and second coordinates of the estimated 
        mean vectors computed by  \cref{alg:Spiked-Sampling-AMP},
        in the case of $\ZZ_2$-synchronization.
        For this experiment we set $n = 1000$, $L = 500$, and $\Delta = 0.02$. }
        \label{fig:trajectory}
\end{figure}

 In our first experiment, we set $L = 500$, $\Delta = 0.02$, and $n = 1000$. 
 In Figure \ref{fig:trajectory}, we plot the trajectories of the first and second coordinates of 
the mean vectors generated by the algorithm:
 $\hbm(\hby_{\ell},\ell\Delta)$, $\ell\in\{0,\dots,L\}$.
 For each realization of the data, we run five independent experiments and
 plot the resulting
 trajectories.  
 
 We see that both $\hm_1$ and $\hm_2$ converge to either $+1$ or $-1$ as $t \to \infty$, 
 regardless of the value of $\beta$. When the signal-to-noise ratio is below the
 recovery  $\beta \leq \beta_{\sIT}(\pi_{\Theta}) = 1$, the trajectory 
  appears to converge to an arbitrary corner. On the contrary, when the signal-to-noise ratio 
  is above the information-theoretic threshold ($\beta > 1$), most trajectories that correspond
   to the same data  $\bX$ consistently converge to the same corner, which 
is correlated with the actual signal $\btheta$.

In our second experiment, we consider the inner product 
$\< \btheta_{\leq 10}, \btheta^{\salg}_{\leq 10} \>$, 
where $\btheta_{\leq 10} \in \RR^{10}$ is the vector comprising of  
the first ten coordinates of $\btheta$ and similarly, 
$\btheta^{\salg}_{\leq 10} \in \RR^{10}$ contains the first ten coordinates of $\btheta^{\salg}$.
This inner product  takes values in the set 
 $\{-10, -8, \cdots, 8, 10\}$. Let $\btheta'$ be a sample from the target posterior distribution 
 $\mu_{\bX}$. 
 If $\btheta^{\salg}$ has distribution close to the target posterior $\mu_{\bX}$,
  we expect the distribution of
  $\< \btheta_{\leq 10}, \btheta^{\salg}_{\leq 10} \>$ to be  close to the 
  one of  $\< \btheta_{\leq 10}, \btheta_{\leq 10}' \>$.
  Let us emphasize that this is not a consequence of  \cref{thm:main2},
  since $\btheta'\mapsto  \< \btheta_{\leq 10}, \btheta_{\leq 10}' \>$
  is not $O(1/\sqrt{n})$-Lipschitz.
  
We can compute an asymptotically exact
prediction for the distribution of $\< \btheta_{\leq 10}, \btheta_{\leq 10}' \>$
as follows.
For large $n$, projecting $\btheta'\sim \mu_{\bX}$ onto an $O(1)$ subset of coordinates
leads to approximately independent entries (modulo an overall sign), 
with marginal distributions given by the AMP state evolution \cite{gerschenfeld2007reconstruction,deshpande2017asymptotic}. 
More explicitly, the distribution of $\< \btheta_{\leq 10}, \btheta_{\leq 10}' \>$
is expected to be approximately the same as the sum of $10$ independent Rademacher
random variables $Z_1,\dots, Z_{10}$, with $\E\{Z_i\} = \hbm_i(\bfzero,0)^2$.
This prediction can be easily evaluated numerically.

\begin{figure}[H]
  \begin{subfigure}{0.33\textwidth}
    \centering
    \includegraphics[width=\linewidth]{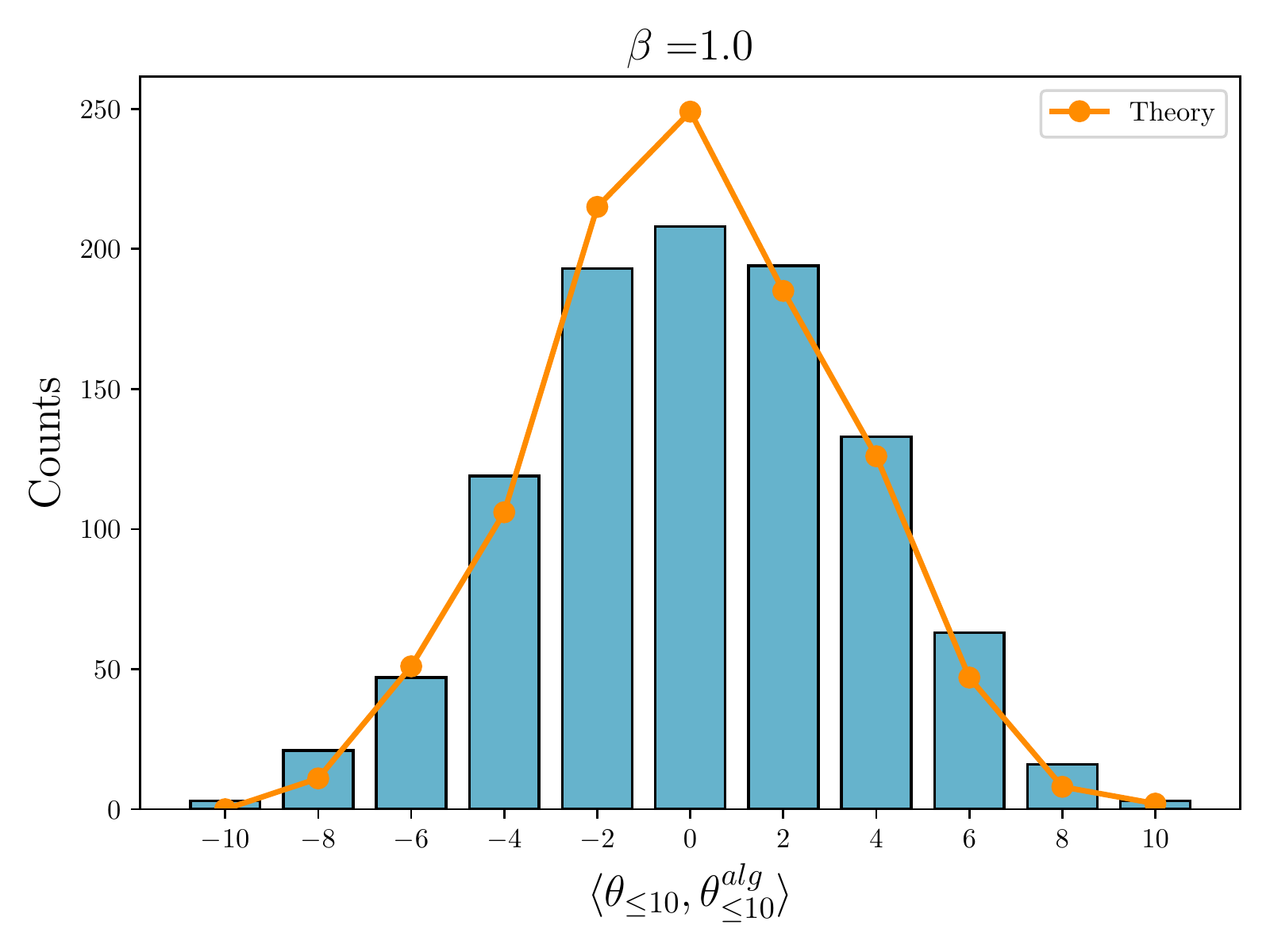}
  \end{subfigure}%
  \begin{subfigure}{0.33\textwidth}
    \centering
    \includegraphics[width=\linewidth]{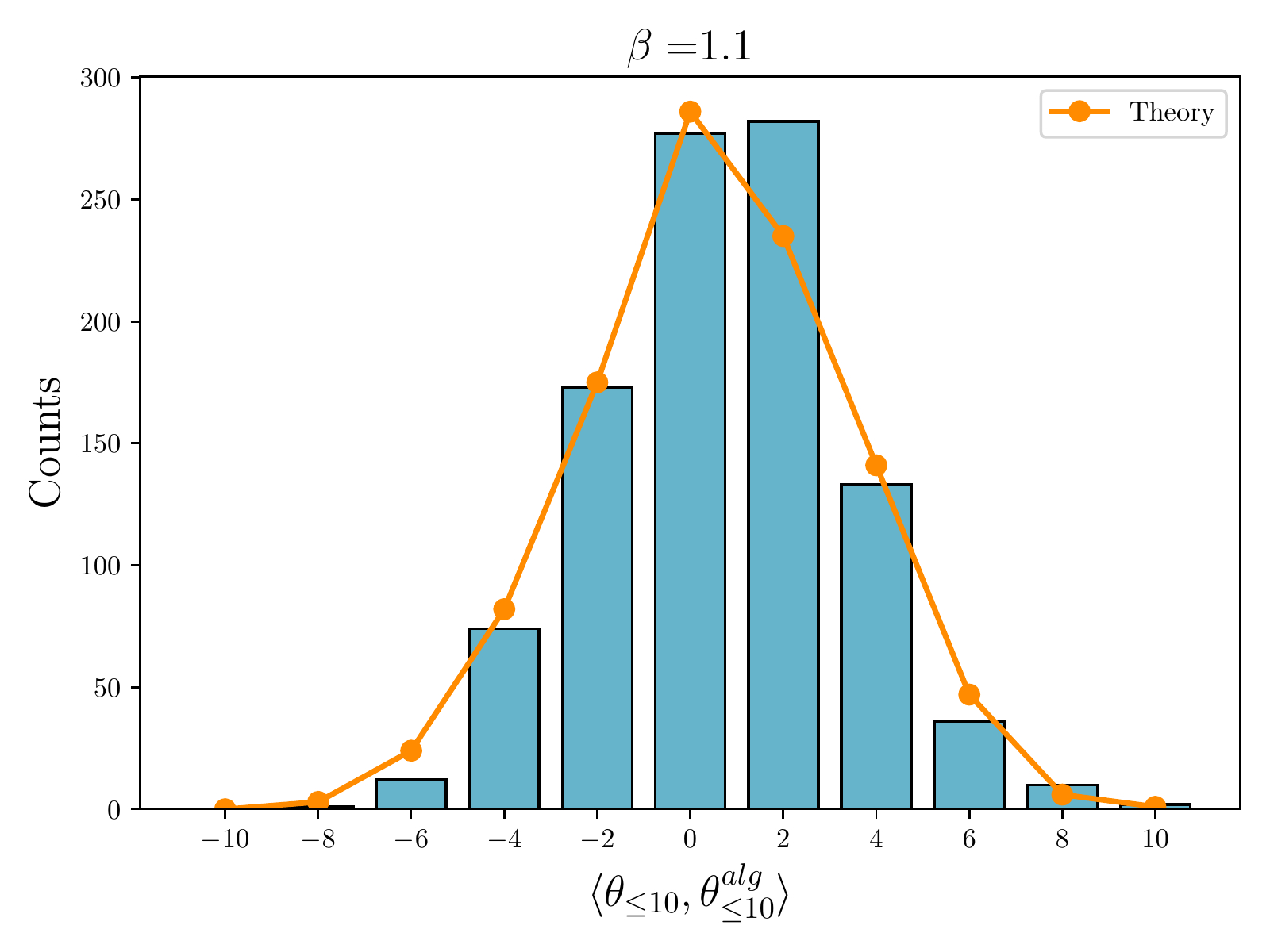}
  \end{subfigure}
  \begin{subfigure}{0.33\textwidth}\quad
    \centering
    \includegraphics[width=\linewidth]{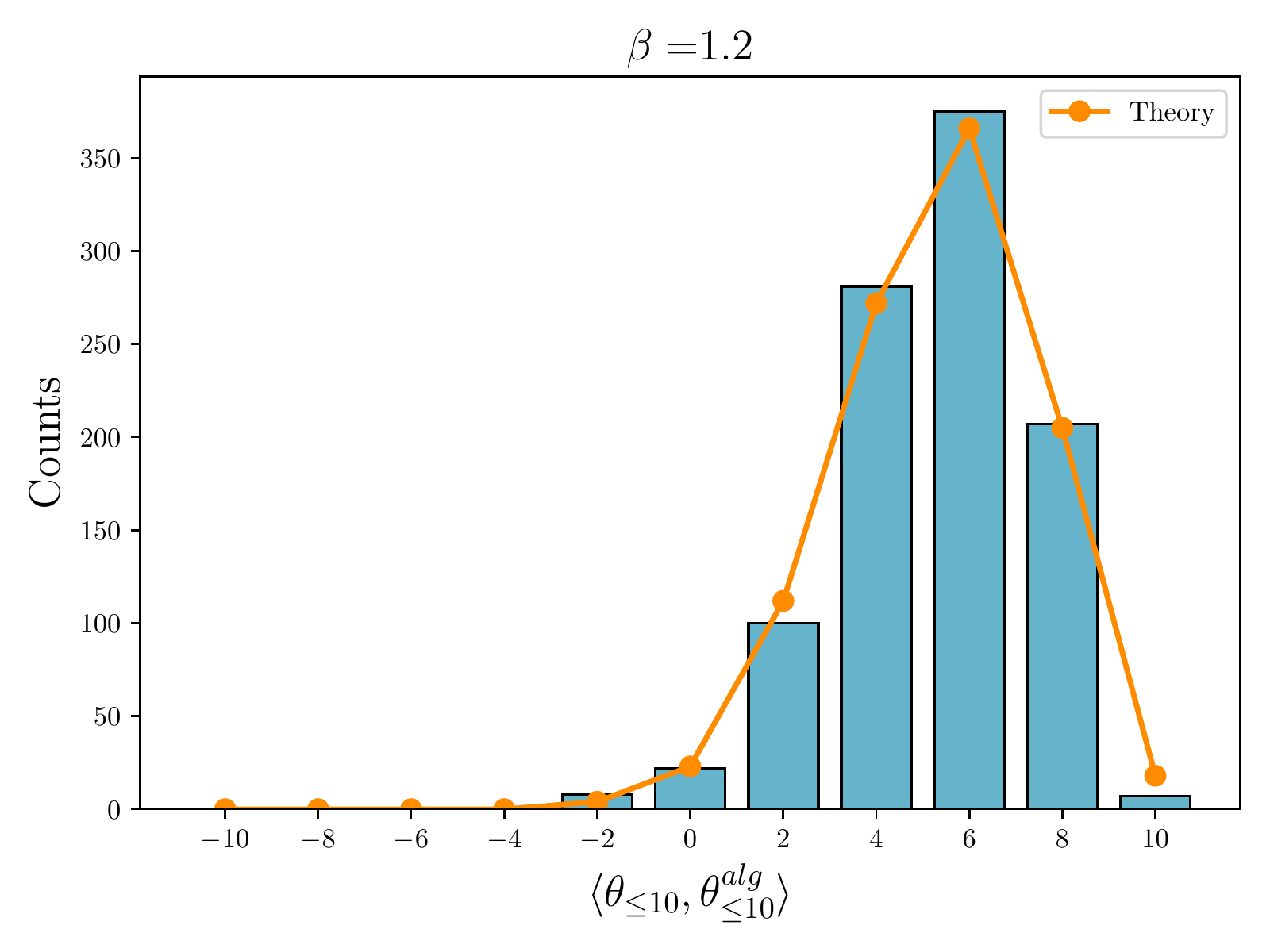}
  \end{subfigure}
  \medskip

  \begin{subfigure}{0.33\textwidth}
    \centering
    \includegraphics[width=\linewidth]{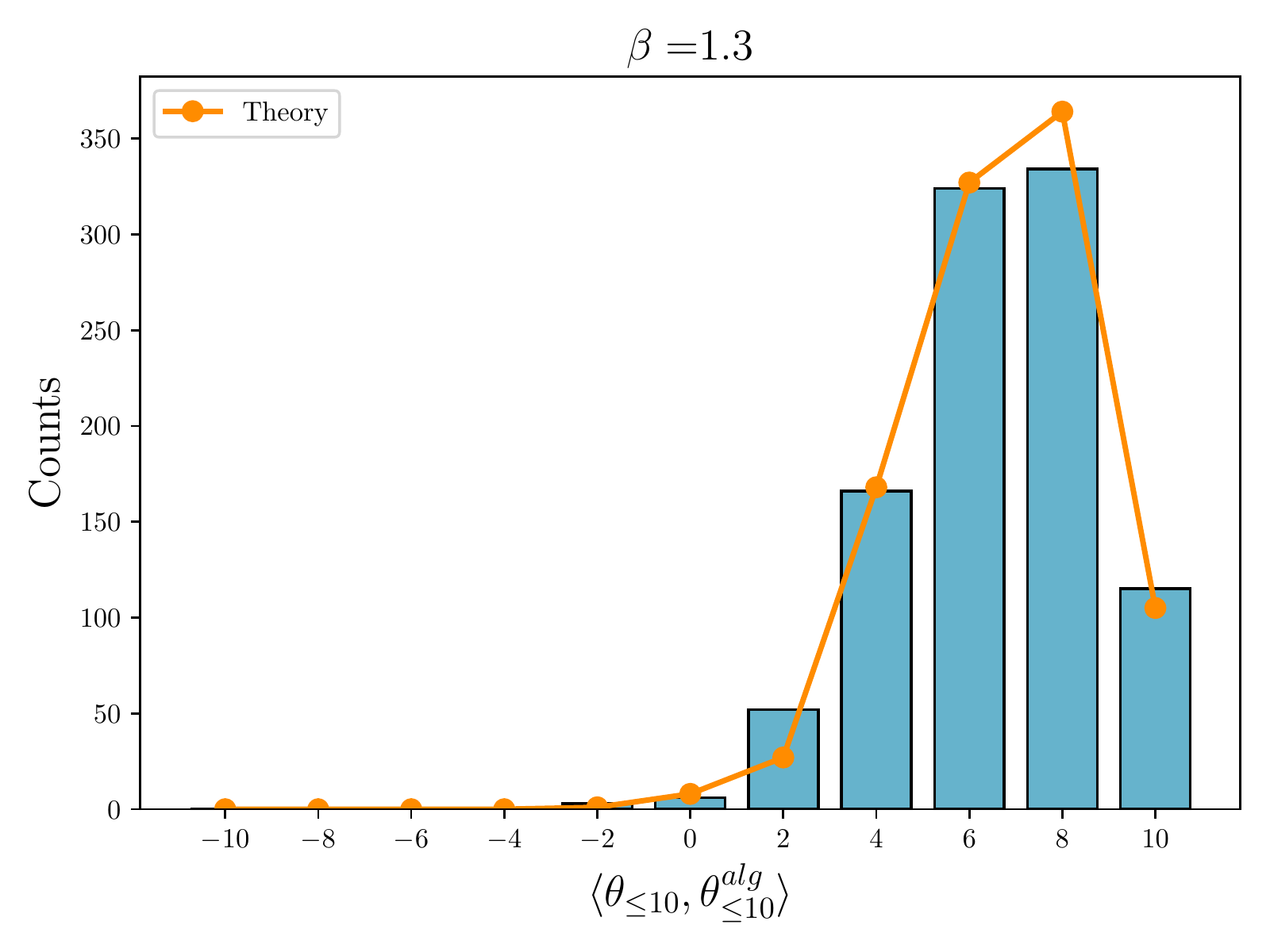}
  \end{subfigure}
  \begin{subfigure}{0.33\textwidth}
    \centering
    \includegraphics[width=\linewidth]{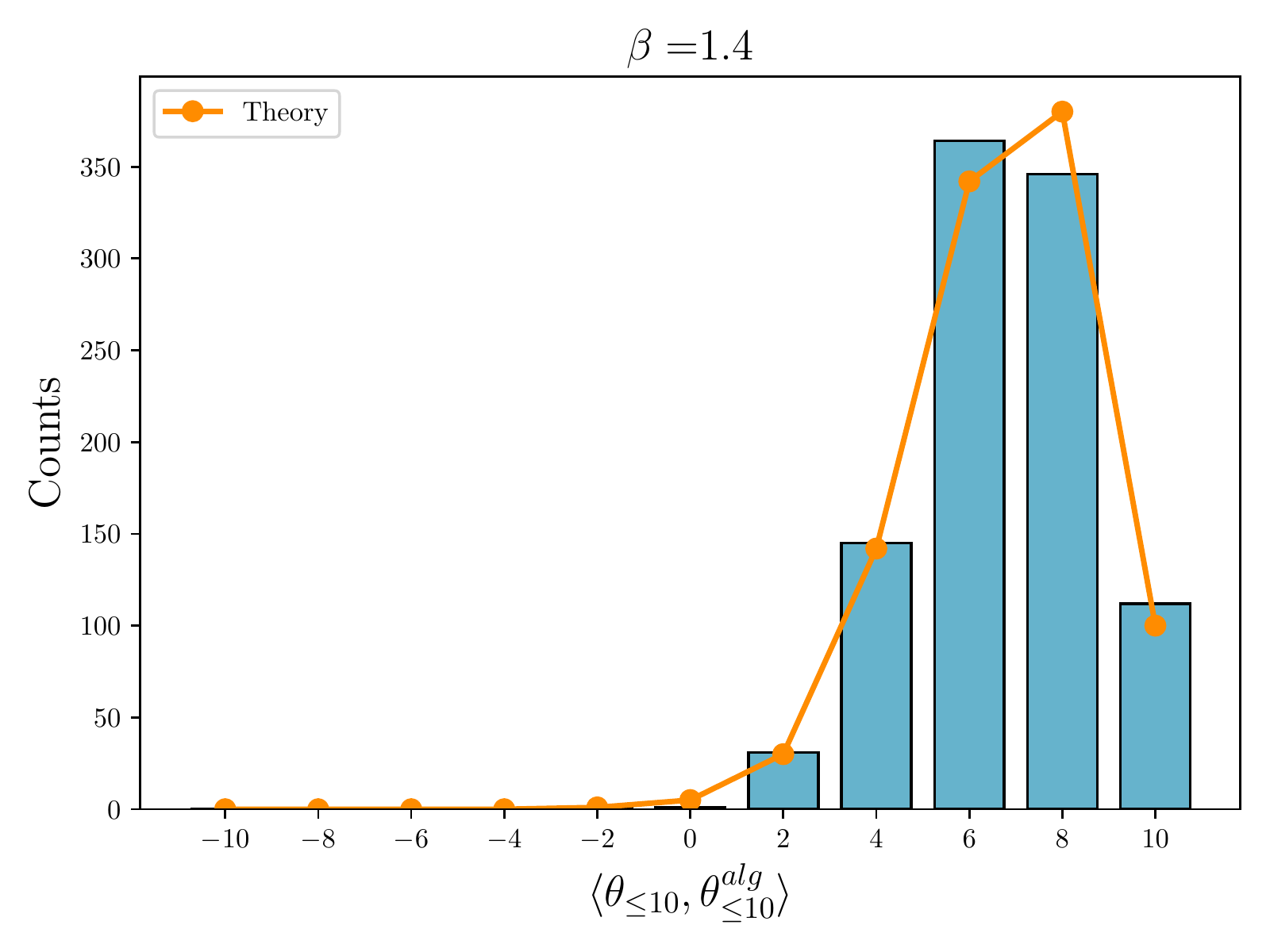}
  \end{subfigure}
  \begin{subfigure}{0.33\textwidth}
    \centering
    \includegraphics[width=\linewidth]{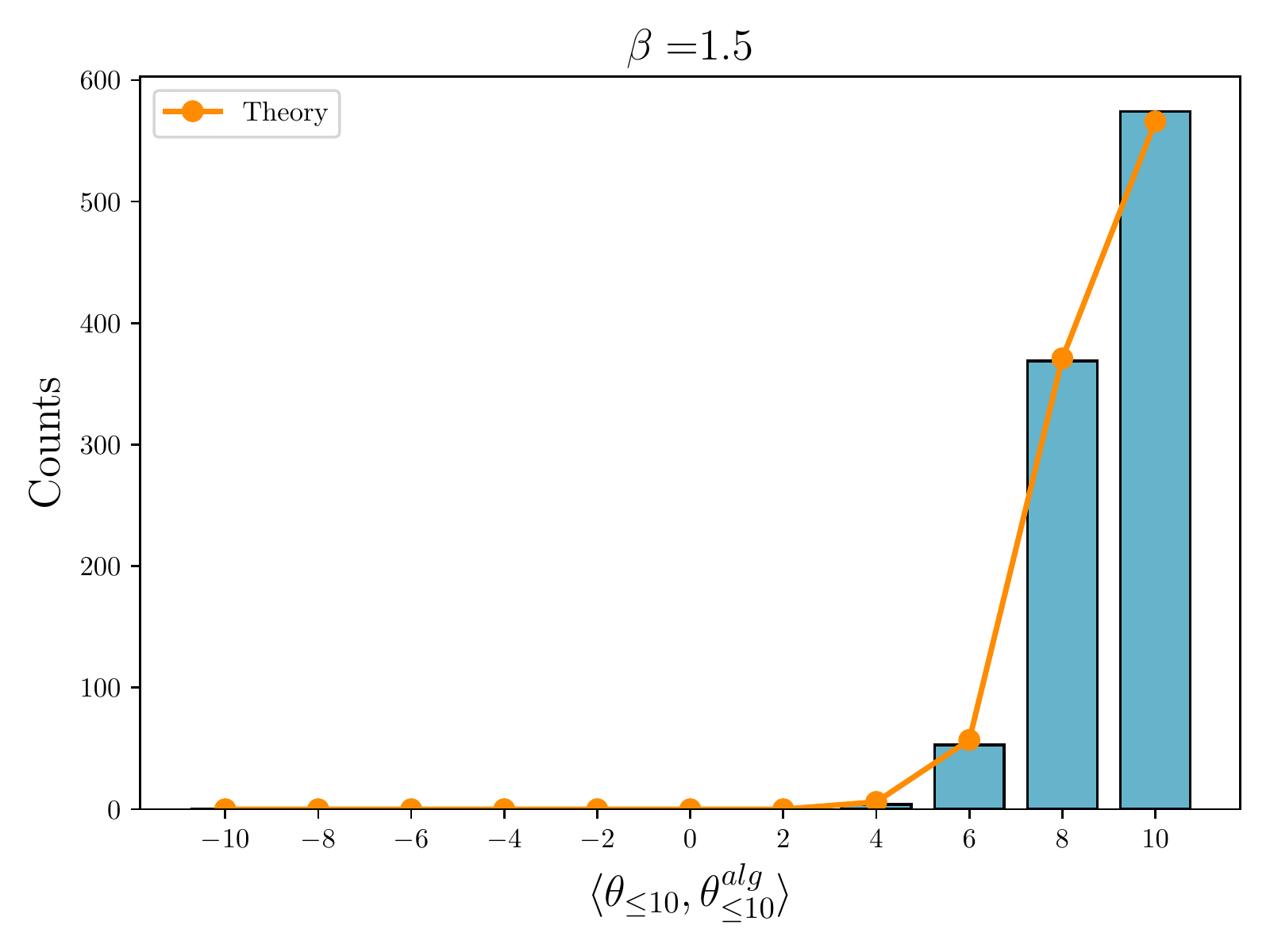}
  \end{subfigure}
  \caption{Histograms: empirical distributions of $\< \btheta_{\leq 10}, \btheta_{\leq 10}^{\salg} \>$ 
  for samples generated by  \cref{alg:Spiked-Sampling-AMP}, based on a single realization
  of the data $(\bX,\btheta)$ at each value of $\beta$. 
  Continuous line: theoretical prediction approximating the distribution of $\langle \btheta_{\leq 10}, \btheta_{\leq 10}^{\salg} \rangle$
  with the true posterior. }
  \label{fig:distribution}
\end{figure}

In Figure \ref{fig:distribution} we compare the empirical distributions of
 $\< \btheta_{\leq 10}, \btheta_{\leq 10}^{\salg} \>$ with
  the theoretical predictions just described. 
  Here, we take
   $n = 1000$, $L = 500$, $\Delta = 0.02$ and multiple values of $\beta$. 
   For each value of $\beta$, we
   draw a single realization $(\bX, \btheta)$, and 
   generate 1000 samples via  \cref{alg:Spiked-Sampling-AMP} for those data. 
    The   experimental outcomes match well with the 
    theoretical predictions, witnessing that the algorithm behaves better 
    than what is guaranteed by our theory.

We just mentioned that, under the posterior, the joint distribution 
of a small subset of the coordinates of $\btheta$ is expected to be well
approximated by a product form. Let us emphasize that this does not mean that the
distribution of the whole vector $\btheta$ is close in $W_{2,n}$ to a product distribution.
In order to highlight the nontrivial correlations in $\mu^{\salg}_{\bX}$,
 we consider the normalized log-likelihood:
\begin{align*}
	\cL(\bX, \btheta^{\salg}) := \frac{\beta}{2n}\langle \btheta^{\salg}, \bX \btheta^{\salg} \rangle. 
\end{align*}
Our theory implies that the distribution of $\btheta^{\salg}$ 
is close to the true posterior. For $\btheta\sim\mu_{\bX}(\,\cdot\,)$, we have 
\begin{align}
\plim_{n\to\infty}	\cL(\bX, \btheta) &= \lim_{n\to\infty}	\E\left\{\int \cL(\bX, \btheta)
\mu_{\bX}(\de\btheta)\right\}\nonumber\\
&= \lim_{n\to\infty}	\frac{\beta}{2n} \E\left\{
\E\big\{\<\btheta, \bX \btheta\>|\btheta\big\}\right\}\nonumber\\
&=\frac{\beta^2}{2}\, \, .\label{eq:logl-limit}
\end{align}
Note that the function $\btheta\mapsto \cL(\bX, \btheta)$ is Lipschitz
(over the domain of interest $\|\btheta\|_2\le \sqrt{n}$) with Lipschitz
constant $(\beta/n)\sup_{\|\btheta\|_2\le\sqrt{n}}\|\bX\btheta\|_{2}\le C(\beta)/\sqrt{n}$
(with high probability with respect to the choice of $\bX$).
Here, $C(\beta) > 0$ is a constant that depends only on $\beta$. 
Hence, \cref{thm:main2} implies that \cref{alg:Spiked-Sampling-AMP} will 
produce samples with the correct expectation for the value of this function.

\begin{figure}[ht]
	\centering
	\includegraphics[width=0.7\linewidth]{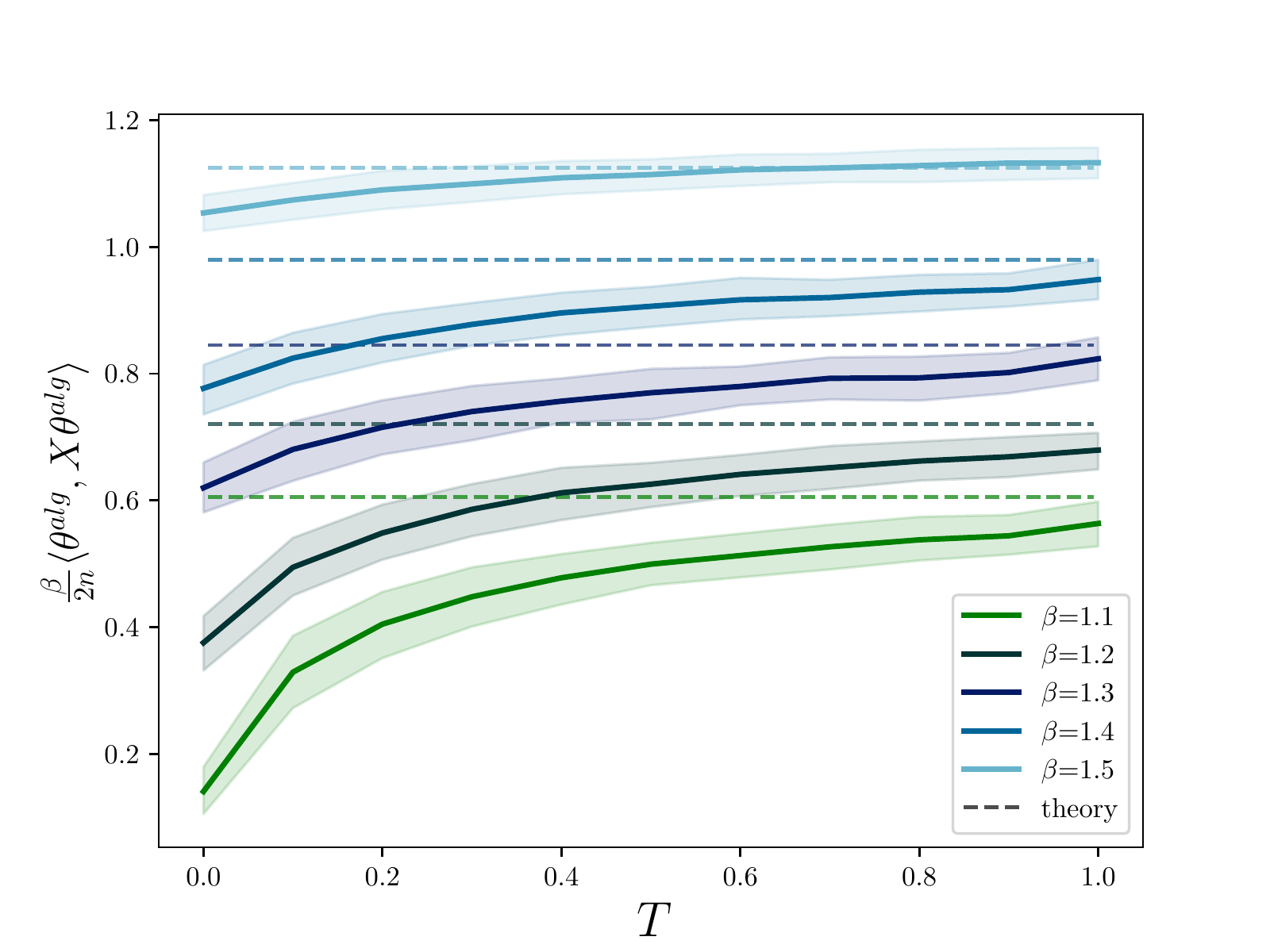}
	\caption{Bands: normalized log-likelihood achieved by \cref{alg:general-sampling-2}.
	 Dashed lines: theoretical predictions.}
	\label{fig:logl}
\end{figure} 

The simple calculation in the last display also shows that the posterior $\mu_{\bX}$ cannot be
approximated in $W_{2,n}$ by a product measure. Indeed, it is possible to show that, 
at least for certain values of $\beta$,
\begin{align}
\plim_{n\to\infty}\frac{\beta}{2n} \<\bm(\bzero,0),\bX\bm(\bzero,0)\> < \frac{\beta^2}{2}\, ,
\end{align}
with strict inequality\footnote{Consider for instance $\beta= 1+\eps$. Then,
 the results of \cite{deshpande2017asymptotic} imply
$\plim_{n\to\infty} \|\bm(\bzero,0)\|_2^2/n\le C\eps$, and therefore 
$\plim_{n\to\infty}(\beta/2n) \<\bm(\bzero,0),\bX\bm(\bzero,0)\>\le C'\eps$ for numerical constants $C, C'$.}.

In Figure \ref{fig:logl}, we compare the theoretical  prediction given in Eq.~\eqref{eq:logl-limit}
with numerical realizations for $\cL(\bX, \btheta^{\salg})$. 
 We fix $\Delta = 0.01$, $n=1000$ and consider  several 
values of the signal-to-noise ratio 
$\beta$ and the number of steps $L$. For each value of the parameters, we repeat the experiment
 independently for 300 times, and display our results. The bands represent the 
 10\% and 90\% quantiles. From the figure, we see that the likelihood increases with $T$, 
 which agrees with the expectation that at larger $T$, $\mu_{\bX}^{\salg}$ matches better with the actual posterior.
 The agreement with the theoretical prediction is again excellent.

\section*{Acknowlegements}

This work was supported by the NSF through award DMS-2031883, the Simons Foundation through
Award 814639 for the Collaboration on the Theoretical Foundations of Deep Learning, the NSF
grant CCF-2006489 and the ONR grant N00014-18-1-2729, and a grant from Eric and Wendy
Schmidt at the Institute for Advanced Studies. Part of this work was carried out while Andrea
Montanari was on partial leave from Stanford and a Chief Scientist at Ndata Inc dba Project N.
The present research is unrelated to AM’s activity while on leave. 

\newpage
\bibliographystyle{alpha}
\bibliography{bib}

%
%

\newpage

\begin{appendices}

\section{Technical preliminaries}

This section summarize some technical facts that will be useful in the proof.

\subsection{Supporting lemmas}

\begin{remark}\label{lemma:nishimori}
Let $\btheta,\bX$ be a couple of random variables (vectors) whose joint distribution is given
by the general Bayesian model \eqref{eq:FirstBayes}. Let $\btheta^{(1)},\dots, \btheta^{(k)}$
be i.i.d. samples from the posterior $\mu_{\bX}(\;\cdot\;) x:= \prob(\;\cdot \; |\bX)$,
independent of $\btheta$. Then
\begin{align}
\bX,\btheta^{(1)},\dots,\btheta^{(k)} \ed \bX,\btheta,\btheta^{(1)},\dots,\btheta^{(k-1)}\, .
\end{align}  
(Here $\ed$ denotes equality in distribution.)

This fact is immediate (just write the joint distribution) and is known in
physics as the ``Nishimori identity.'' 
\end{remark}

\begin{lemma}[Lemma 3.2 in \cite{panchenko2013sherrington}]\label{lemma:diff-of-derivative}
	If $f$ and $g$ are two differentiable convex functions, then for any $b > 0$, 
	\begin{align*}
		|f'(a) - g'(a)| \leq g'(a + b) - g'(a - b) + \frac{d}{b},
	\end{align*}
	where $d = |f(a + b) - g(a + b)| + |f(a - b) - g(a - b)| + |f(a) - g(a)|$. 
\end{lemma}

\begin{lemma}[Lemma 4.15 in \cite{alaoui2022sampling}]\label{lemma:rounding}
	Suppose probability distributions $\mu_1, \mu_2$ on $[-1,1]^n$ are given. Sample $\bm_1 \sim \mu_1$ and $\bm_2 \sim \mu_2$ and let $\btheta_1, \btheta_2 \in \{-1,+1\}^n$ be standard randomized roundings, respectively of $\bm_1$ and $\bm_2$. (Namely, the coordinates of $\btheta_i$ are conditionally independent given $\bm_i$, with $\E[\btheta_i \mid \bm_i] = \bm_i$.) Then
	\begin{align*}
		W_{2,n}(\mathcal{L}(\btheta_1), \mathcal{L}(\btheta_2)) \leq 2 \sqrt{W_{2,n}(\mu_1, \mu_2)}. 
	\end{align*}
\end{lemma}

\begin{lemma}[Proposition 2.1.2 in \cite{vershynin2018high}]\label{lemma:gaussian-tail}
	Let $g \sim \normal(0,1)$. Then for all $t > 0$, it holds that 
	\begin{align*}
		\left(\frac{1}{t} - \frac{1}{t^3}\right) \cdot \frac{1}{\sqrt{2\pi}} e^{-t^2 / 2} \leq \P(g \geq t) \leq \frac{1}{t} \cdot \frac{1}{\sqrt{2\pi}} e^{-t^2 / 2}.
	\end{align*}
\end{lemma}

\begin{lemma}
\label{lemma:a-ratio}
	For any $v \in \RR\backslash \{0\}$ and $\bar a >  a > 0$, it holds that 
	\begin{align*}
		\frac{1 - e^{-v\bar a}}{\bar a} \leq \frac{1 - e^{-v a}}{ a}. 
	\end{align*}
\end{lemma}
\begin{proof}[Proof of \cref{lemma:a-ratio}]
	Define $f(x) = x^{-1} (1 - e^{-v x})$ for $x > 0$. Then, $f'(x) = x^{-2} (e^{-vx} + vx e^{-vx} - 1)$. If $v \neq 0$, then $f'(x) \leq 0$ for all $x > 0$. The proof is complete. 
\end{proof}

\begin{lemma}
\label{lemma:strong-convexity-property}
	Let $f: \RR^n \to \RR$ be smooth and $\mu$-strongly convex for $\mu > 0$,
	and  denote by $\bx_{\ast}$ a local (hence global) minimizer of $f$. Then it holds that 
	\begin{align*}
		\frac{1}{2\mu} \|\nabla f(\bx)\|_2^2 \geq f(\bx) - f(\bx_{\ast}) \geq \frac{\mu}{2} \cdot \|\bx - \bx_{\ast}\|_2^2. 
	\end{align*}
\end{lemma}

\begin{lemma}
\label{lemma:strong-convexity-min}
	Let $f: \RR^n \to \RR$ be $\mu$-strongly convex for $\mu > 0$,
	and assume $\|\nabla f(\bx_0)\|\le c_0$. Then $f$ has a unique local minimum
	$\bx_*$, satisfying $\|\bx_*-\bx_0\|\le c_0/\mu$. 
\end{lemma}



\section{Proofs for the general sampling scheme}
\subsection{Proof of Theorem \ref{thm:gen}}
\label{app:ProofGen}

	We  couple $\{\bw_\ell\}_{1 \leq \ell \leq L}$ and $\{\bB(t)\}_{0 \leq t \leq T}$ 
	by letting $\bw_{\ell}= \bB(\ell\Delta)-\bB((\ell-1)\Delta)$.
	We also define $A_\ell = \|\hby_{\ell} - \by(\ell\Delta)\|_2 / \sqrt{N}$ for all $\ell \in \{0\} \cup [L]$, and write $t_{\ell}:=\ell\Delta$.
	
	Let $\Omega$ be the intersection of the events listed in points \hyperlink{A1}{$\mathsf{(A1)}$}, \hyperlink{A2}{$\mathsf{(A2)}$} and \hyperlink{A3}{$\mathsf{(A3)}$}. 
	By taking a union bound, we obtain that $\prob(\Omega)\ge 1-5\eta$. 
	We will prove by induction that, on $\Omega$, the following holds for all $\ell\le L$:
	\begin{align}
	\hby_{\ell}\in B(\ell)\, ,\;\;\mbox{ and }\;\;\;
		A_\ell \leq \frac{C_1\sqrt{\Delta} + \eps_1 + \eps_2 + \eps_3}{C_2}\cdot \big(e^{C_2\ell\Delta}-1\big)\, .
		\label{eq:BoundIndAell}
	\end{align}
	By definition, we see that 
	$A_0 = 0$ and $\hby_{0}=\by(0)=\bfzero\in B(0)$. Next, assume that the induction hypothesis holds up to step  $\ell-1$.
On the event $\Omega$:
	\begin{align*}
		 A_{\ell} - A_{\ell - 1} 
		 \leq & \;\frac{1}{\sqrt{N}} \int_{t_{\ell - 1}}^{t_{\ell}}
		 \|\hbm(\hby_{\ell - 1},t_{\ell - 1}) - \bm(\by(t), t)\|_2 \dd t \\
		\leq 
		&\;
		 \frac{\Delta}{\sqrt{N}}\|\hbm(\by(t_{\ell - 1}),t_{\ell - 1}) - \bm(\by(t_{\ell - 1}),t_{\ell - 1})\|_2 \\
		& + \sup_{t \in [t_{\ell - 1}, t_{\ell}]} \frac{\Delta}{\sqrt{N}} \|\bm(\by(t), t) - \bm(\by(t_{\ell - 1}), t_{\ell - 1})\|_2 \\
		&\;+\frac{\Delta}{\sqrt{N}} \|\hbm(\hby_{\ell - 1},t_{\ell - 1}) - \hbm(\by(t_{\ell - 1}),
		t_{\ell - 1})\|_2 \\
		\leq & \; \Delta  \cdot \left(\eps_1+C_1\sqrt{\Delta}  + \eps_2 + C_2 A_{\ell - 1} + \eps_3 \right). 
	\end{align*} 
	Substituting in the induction hypothesis, we obtain 
	$A_\ell \leq \frac{C_1\sqrt{\Delta} + \eps_1 + \eps_2 + \eps_3}{C_2}\cdot \big(e^{C_2\ell\Delta}-1\big)$
	as desired. The claim $\hby_{\ell}\in B(\ell)$ follows from the stated condition on
	$r_{\ell}$. 
	This completes the induction proof for \cref{eq:BoundIndAell}.
	
Applying the bound given in \cref{eq:BoundIndAell} with $\ell=L$
and using once more assumptions  \hyperlink{A1}{$\mathsf{(A1)}$} and \hyperlink{A3}{$\mathsf{(A3)}$},
we see that on $\Omega$, it holds that (recall $T = L \Delta$)
	\begin{align}\label{eq:m-AT}
		\frac{1}{\sqrt{n}}\|\bm_{\btheta}(\by(T), T)  - \hbm_{\btheta}(\hby_L,T)\|_2 \leq & \frac{1}{\sqrt{n}} 
		\|\bm_{\btheta}(\by(T), T) - \hbm_{\btheta}(\by(T),T)\|_2\nonumber\\
		& + \frac{1}{\sqrt{n}} \|\hbm_{\btheta}(\by(T),T) - \hbm_{\btheta}(\hby_L,T)\|_2 \nonumber \\
		\leq & \eps_1 + C_2A_L + \eps_3 \nonumber \\
		\leq  & \eps_1 + \eps_3 + (C_1 \sqrt{\Delta} + \eps_1 + \eps_2 + \eps_3)\cdot e^{C_2T}=:\Delta_0\, .
	\end{align}
	The above upper bound further implies that (denoting by $\sP_{\oR}$ the projection
	onto $\Ball^n(\bfzero,\oR \sqrt{n})$)
	\begin{align}
		 & W_{2, n}(\mu_{\data}, \mu_{\data}^{\salg}) \\
		\leq & W_{2, n}(\mu_{\data}, \Law( \bm_{\btheta}(\by(T), T))) + W_{2, n}(\Law( \bm_{\btheta}(\by(T), T)), \mu_{\data}^{\salg}) \nonumber\\
		\leq & \frac{1}{\sqrt{n}} \EE[\|\hbm_{\btheta}(\hat{\by}_L,L\delta) - \bm_{\btheta}(\by(T), T)\|_2^2]^{1/2}  +  W_{2, n}(\mu_{\data}, \Law( \bm_{\btheta}(\by(T), T))) \nonumber \\
		\leq & \frac{1}{\sqrt{n}} \big( \EE[\|\sP_{\oR}(\hbm_{\btheta}(\hat{\by}_L,L\delta)) - \sP_{\oR}(\bm_{\btheta}(\by(T), T))\|_2^2]^{1/2} + \EE[\|\bm_{\btheta}(\by(T), T) - \sP_{\oR}(\bm_{\btheta}(\by(T), T))\|_2^2]^{1/2}\big) \nonumber\\
		& + W_{2, n}(\mu_{\data}, \Law( \bm_{\btheta}(\by(T), T)))  \\ 
		\leq & \Delta_0 +\E\Big[\frac{1}{n}\|\bm_{\btheta}(\by(T), T)\|_2^2\bbone_{\|\bm_{\btheta}(\by(T), T)\|_2\ge \overline{R}\sqrt{n}}\Big]^{1/2}+
		10 \overline{R}\eta +
		 W_{2, n}\big(\mu_{\data}, \Law( \bm_{\btheta}(\by(T), T))\big) \, .\nonumber 
		\end{align}
	This implies  \cref{eq:ZeroGeneral}.   \cref{eq:SecondGeneral} follows
	by using the moment assumption to bound the 
	expectation on the right-hand side and
	optimizing over
	$\overline{R}$, and finally applying Lemma \ref{lemma:W2-distance}.

%
%
\subsection{Concentration of the stochastic localization process}

The next lemma is a slight generalization of analogous results in
\cite{eldan2020taming,alaoui2022sampling}.
\begin{lemma}\label{lemma:W2-distance}
Let $\mu\in \cuP_2(\reals^n)$ be a probability measure with finite second moment, and
$\by(t) = t\bH\btheta +\sqrt{t}\bg$ for $(\btheta,\bg)\sim\mu\otimes 
\normal(\bzero,\id_n)$.
Further, let $\mu_t(\,\cdot\,):= \prob(\btheta\in\,\cdot\,|\by(t))$, and
$\bm(\by(t); t) :=\E[\btheta|\by(t)]$.
Finally, denote by $\proj_{\ker(\bH)}$ the projector onto the null space of $\bH$.

Then the following inequalities hold for all  $t > 0$:
	\begin{align}
	\E\Cov(\mu_t) & \preceq \proj_{\ker(\bH)}\Cov(\mu)\proj_{\ker(\bH)}+
	\frac{1}{t}\bH^{+}(\bH^{+})^{\sT}\, ,\label{eq:FirstGeneralBound}\\
		W_{2,n}(\mu, \Law(\bm(\by(t); t)))^2 &\leq \frac{1}{n}
		\Tr(\proj_{\ker(\bH)}\Cov(\mu))+
		 \frac{1}{nt} \Tr(\bH^{+}(\bH^{+})^{\sT}
		\big)\, .\label{eq:SecondGeneralBound} 
	\end{align}
\end{lemma}
\begin{proof}
By rescaling $\bH$, we can assume without loss of generality $t=1$.
We can also center $\mu$ so that $\E(\btheta) = \bfzero$.
We will write, for simplicity, $\by=\by(1)$.
 Note that
\begin{align}
\E\Cov(\mu_t) & = \E\Big\{\big(\btheta-\E(\btheta|\by)\big)\big(\btheta-\E(\btheta|\by)\big)^{\sT}\Big\}\\
& \le \E\Big\{\big(\btheta-\bH^{+}\by\big)\big(\btheta-\bH^{+}\by\big)^{\sT}\Big\}\\
&= \proj_{\ker(\bH)}\Cov(\mu)\proj_{\ker(\bH)}+\bH^{+}(\bH^{+})^{\sT}\, ,
\end{align}
where the inequality follows by the optimality of posterior expectation
under quadratic losses.
This proves the first claim \eqref{eq:FirstGeneralBound}. 

In order to prove the second one, denote
by $B(n,t)$ the right-hand side of Eq.~\eqref{eq:SecondGeneralBound}.
By taking the trace of the former inequality, we obtain 
\begin{align}
\E \big\{W_{2,n}(\mu_t,\delta_{\bm(\by(t),t)})^2\big\} = \frac{1}{n}\E\Tr\Cov(\mu_t)\le B(n,t)\, .
\end{align}
Since $(\mu,\nu)\mapsto W_2(\mu,\nu)^2$ is jointly convex in $(\mu,\nu)$,
Jensen's inequality implies
\begin{align}
\E \big\{W_{2,n}(\mu_t,\delta_{\bm(\by(t),t)})^2\big\}\ge 
W_{2,n}(\E\mu_t,\E\delta_{\bm(\by(t),t)})^2 = W_{2,n}(\mu, \Law(\bm(\by(t); t)))^2 \, ,
\end{align}
which completes our proof.
\end{proof}
%
%

\section{Proofs for \cref{alg:Spiked-Sampling-AMP}}
\label{sec:proof-alg:Spiked-Sampling-AMP}

\subsection{Proof of \cref{lemma:sign}}\label{sec:proof-of-lemma:sign}

The proof of \cref{lemma:sign} will be based on the following lemma, which is a straightforward
consequence of the fact that the characteristic function uniquely identifies the corresponding 
probability measure.
\begin{lemma}\label{lemma:symmetric-characteristic}
Let $\rP$ be a probability measure on $\RR$.
Then $\rP$ is symmetric (i.e. $\rP(A) =\rP(-A)$ for every Borel set $A$) if
and only if its characteristic function $\varphi_{\rP}$ is real-valued or,
equivalently, if and only if $\varphi_{\rP}(t) = \varphi_{\rP}(-t)$
for every $t\in\RR$.
\end{lemma}

Recall that $\bnu$ is a top eigenvector of $\bX$ with norm
$\|\bnu\|_{2}^2 = n\beta^2(\beta^2-1)$. We denote by $\lambda_1$
the corresponding eigenvalue, and note that this is almost surely non-degenerate 
(because the law of $\bX$ is absolutely continuous with respect to Lebesgue).
Let 
\begin{align}
\bnu_+ = s\bnu\, ,\;\;\;\;\; s:= \sign \< \bnu, \btheta \>\, .
\end{align}
Note that we can assume $s$ independent of $\btheta,\bW$ 
(because we can define $\bnu$ to be taken uniformly at random among the 
two eigenvectors with given norm.)

For any $\bOmega \in \Orth(n)$ satisfying $\bOmega \btheta = \btheta$ and 
is independent of $\bW$, we have $\bOmega \bW \bOmega^{\sT} \ed \bW$. 
Moreover, if we replace $\bW$ by $\bOmega \bW \bOmega^{\sT}$, then $\lambda_1$ 
is the top eigenvalue of $\bOmega \bX \bOmega^{\sT} = \beta \btheta \btheta^{\sT} 
/ n + \bOmega \bW \bOmega^{\sT}$ and $\bOmega \bnu$ is the corresponding eigenvector.
 As a result, we can conclude that the following two conditional distributions 
 are equal:
\begin{align*}
	\bOmega \bW \bOmega^{\sT}, \bOmega \bnu_+ \mid \btheta, \bOmega \ed 
	 \bW,  \bnu_+ \mid \btheta, \bOmega\, .
\end{align*}
Let $\bP^{\perp}_{\btheta}$ be the projector orthogonal to $\btheta$.
The above invariance implies that, conditioning on 
$\btheta$, $\<\btheta,\bnu_+\>_+/\|\btheta\|_2 = \rho_{\|}$ and 
$\|\bP_{\btheta}^{\perp}\bnu_+\|_2=\rho_{\perp}$, we have
(on $\btheta\neq \bzero$):
\begin{align*}
\bnu_+ = \rho_{\|}\frac{\btheta}{\|\btheta\|_2}+ 
\rho_{\perp} \bu\, ,
\end{align*}
where $\bu$ is a uniformly random unit vector orthogonal to $\btheta$.
Equivalently,
\begin{align*}
\bnu_+ = \rho_{\|}\frac{\btheta}{\|\btheta\|_2}+ 
\rho_{\perp} \frac{\bP^{\perp}_{\btheta}\bg}{\|\bP^{\perp}_{\btheta}\bg\|_2}\, ,
\end{align*}
where $\bg\sim\normal(\bzero,\id_n)$  independent of $\btheta$.

By \cite{benaych2011eigenvalues} 
$\rho_{\|}/\sqrt{n} = (\beta^2-1)+o_{n,P}(1)$ and therefore, using the normalization of
$\bnu_+$, $\rho_{\perp}/\sqrt{n} = (\beta^2-1)^{1/2}+o_{n,P}(1)$.
Using the fact that $\|\btheta\|_2= \sqrt{n} +O_{n,P}(1)$  and $\|\bP^{\perp}_{\btheta}\bg\|_2 =
\sqrt{n} +O_{n,P}(1)$ (both hold by the law of large numbers,
since $\int \theta^2\, \pi_{\Theta}(\de\theta)=1$), we obtain
\begin{align*}\label{eq:bnu}
	\plim_{n\to\infty}\frac{1}{n} \big\|\bnu_+ - (\beta^2 - 1) \btheta - \sqrt{\beta^2 - 1} \bg  \big\|_2^2 = 0\, .
\end{align*}

By \cref{lemma:symmetric-characteristic} we can  assume without loss of generality 
$\Im\varphi_{\pi_\Theta}((\beta^2-1)t_0)>\delta_0>0$ for some $\delta_0,t_0>0$ (we use $\Im$ to denote the imaginary part of a complex number). 
The proof for the other sign can be handled symmetrically.
We then define 
\begin{align}
\cA(\bnu):= \sign T_n(\bnu) \, , \;\;\;\;\;
T_n(\bnu) := \frac{1}{n}\sum_{i=1}^n\sin(t_0\nu_i)\, .
\end{align}
Note that $T_n(\bnu) = s\, T_n(\bnu_+)$ and therefore the proof is completed by showing that,
with high probability $T_n(\bnu_+)>0$. Indeed, let $\delta_1:= \exp(-(\beta^2-1)t_0/2)\delta_0$.
Then, letting $\overline\bnu_+:=(\beta^2 - 1) \btheta - \sqrt{\beta^2 - 1} \bg$, 
by Eq~\eqref{eq:bnu}, we have, with high probability
\begin{align}
\big|T_n(\overline\bnu_+)-T_n(\bnu_+)\big|\le \frac{\delta_1}{2}\, .
\end{align}
On the other hand, by the law of large numbers (for $(\Theta,G)\sim \pi_{\Theta}\otimes\normal(0,1)$)
\begin{align}
\plim_{n\to\infty}T_n(\overline\bnu_+) &= \E\sin(t_0(\beta^2 - 1) \Theta- \sqrt{t_0(\beta^2 - 1)} G)\\
& =  \exp\big(-(\beta^2-1)t_0/2\big) \, \Im\varphi_{\pi_\Theta}((\beta^2-1)t_0)\ge \delta_1\, .
\end{align}
Together with the previous display, we complete the proof.

\subsection{Proof of \cref{thm:main2}}\label{sec:appendix-spiked}

This section is devoted to analyzing \cref{alg:Spiked-Sampling-AMP}. In particular, we will outline the proof of
 \cref{thm:main2}, while delaying most of the technical details to Appendix \ref{sec:supporting-lemmas}.

\subsubsection{The tilted measure}
\label{sec:Tilted}

As stated in \cref{alg:Spiked-Sampling-AMP}, after $\ell$ steps, the state of the algorithm is given by vectors
$\hby_{\ell}$ and $\hbm(\hby_{\ell},\ell\Delta)$. In particular, $\hbm(\hby_{\ell},\ell\Delta)$
is  interpreted as an estimate the posterior mean of $\btheta$ given $\by(\ell \delta) = \hby_{\ell}$ and $\bX$. 
  
This interpretation has to be slightly modified if $\pi_{\Theta}$ is
symmetric.
To be explicit, for any $\by \in \RR^n$ and $t > 0$, we define the `tilted measure' $\mu_{\bX, \by, t}$ as follows
(in the formulas below $Z(\bX, \by, t)$  are normalizing constants defined by $\int \mu_{\bX, \by, t}(\de\btheta) = 1$):
\begin{itemize}
	\item \textbf{If $\pi_{\Theta}$ is not a symmetric distribution,} then 
	$\mu_{\bX, \by, t}$ is the posterior distribution of $\btheta$ given  $\bX$ and $\by(t)=\by$:
\begin{align*}
	& \mu_{\bX, \by, t}(\de\btheta) := \frac{1}{Z(\bX, \by, t)} \exp \left( \frac{\beta}{2} \< \btheta, \bX \btheta \>
	-\frac{\beta^2}{4n}\|\btheta\|_2^4
	+ \langle \by, \btheta \rangle - \frac{t}{2}\|\btheta\|_2^2 \right) \pi_{\Theta}^{\otimes n}(\de\btheta) \, .
\end{align*}  
\item \textbf{If $\pi_{\Theta}$ is a symmetric distribution,} 
then we let $\bnu(\bX)$ be a randomly selected leading eigenvector of
 $\bX$\footnote{Almost surely, the leading eigenvalue 
of $\bX$ is non-degenerate, and therefore there are two choices $\{+\bv_1,-\bv_1\}$
for the normalized leading eigenvector. 
We let $\bnu(\bX)\sim\Unif(\{+c\bv_1,-c\bv_1\})$, $c:=\sqrt{n\beta^2(\beta^2-1)}$ independently
of $\bX$, $\by(t)$.}. Then we break the symmetry by conditioning on the sign of
$\<\btheta, \bnu(\bX)\>$. Namely, we let 
\begin{align*}
	& \mu_{\bX, \by, t}(\de\btheta) := \frac{1}{Z(\bX, \by, t)}
	 \exp \left( \frac{\beta}{2} \langle \btheta, \bX \btheta \rangle -\frac{\beta^2}{4n}\|\btheta\|_2^4
	 + \langle \by, \btheta \rangle - \frac{t}{2}\|\btheta\|_2^2 \right)
	 \bbone_{ \langle \btheta, \bnu(\bX) \rangle \geq 0}\, \pi_{\Theta}^{\otimes n}(\de\btheta).
\end{align*}
\end{itemize}
For the symmetric prior case, 
 the relation between $\mu_{\bX, \by, t}$ and the actual posterior distribution is
 given by 
 $$\prob(\btheta\in A|\bX,\by(t)=\by) =\frac{1}{2}\,\mu_{\bX, \by, t}(A)+\frac{1}{2}\,\mu_{\bX, \by, t}(-A).$$ 
As a consequence, if we can approximately draw samples from the distribution $\mu_{\bX, \bfzero, 0}$, then we can also sample from the target posterior 
 $\prob(\btheta\in \,\cdot\, |\bX)$ with the same approximation guarantees in $W_2$
 distance. 
 Therefore, it is sufficient to generate $\btheta \sim \mu_{\bX, \bfzero, 0}$
 and then flip its sign with probability $1/2$, which is what we do in \cref{alg:Spiked-Sampling-AMP}.
 Hereafter, we will focus on sampling from $\mu_{\bX, \bfzero, 0}$.
 We note that in both cases $\mu_{\bX, \by, t}$ is the tilted measure for $\mu_{\bX, \bfzero, 0}$. 

Throughout this proof, with a slight abuse of notations we use $\bm(\by, t)$ to denote the mean of the tilted measure
$\bm(\by, t) := \int \btheta \,  \mu_{\bX, \by, t}(\dd\btheta)$.

\subsubsection{Proof outline}

The proof consists of checking the assumptions of Theorem
\ref{thm:gen}. 
Namely, in Section \ref{sec:AMP-optimal}  we prove that the AMP estimate is close to the 
expectation, thus verifying \hyperlink{A1}{$\mathsf{(A1)}$};
 in Section \ref{sec:path-regularity}
we verify the path-regularity assumption \hyperlink{A2}{$\mathsf{(A2)}$};
finally, in Section  \ref{sec:AMP-stable} we establish the Lipschitz continuity of the AMP estimate,
verifying  assumption \hyperlink{A3}{$\mathsf{(A3)}$}.

We conclude the proof in Section \ref{sec:proof}.
Throughout the proof, the prior distribution $\pi_{\Theta}$ is fixed and hence we do not detail the dependency of various quantities on $\pi_{\Theta}$. 
We will on the other hand track dependencies on other objects by emphasizing them inside parentheses.
In this proof, we use $\hbm^k(\by,t)$ to denote the estimate produced by the AMP algorithm
of Eq.~\eqref{eq:general-AMP} after $k$ iterations, with inputs $\bX,\by$ and $\bnu$ that is computed using the approach given in \cref{alg:Spiked-Sampling-AMP}.

\subsubsection{AMP achieves Bayes optimality}\label{sec:AMP-optimal}

The analysis of the Bayes AMP algorithm uses the following characterization in terms of 
state evolution, which is adapted from \cite{montanari2021estimation}.
Here, we refer to a function $\psi:\reals^m\to\reals$ as \emph{pseudo-Lipschitz}
if $|\psi(\bx_1)-\psi(\bx_2)|\le C(1+\|\bx_1\|_2+\|\bx_2\|_2)\|\bx_1-\bx_2\|_2$ holds for
 all $\bx_1, \bx_2 \in \RR^m$.
 
 We note that the analysis in this section requires only moment conditions on $\pi_{\Theta}$. Namely, we do not assume $\pi_{\Theta}$ is discrete or continuous.  
 
%
\begin{proposition}[\cite{montanari2021estimation}]\label{propo:SE-Basic}
Consider the Bayes AMP algorithm with spectral initialization, defined in Eq.~\eqref{eq:general-AMP},
and the state evolution recursion of Eq.~\eqref{eq:AMPSE}. Assume $\pi_{\Theta}$ has unit second moment.
Then, for any fixed $k \in \NN_+$, $t\ge 0$, and any pseudo-Lipschitz
test function $\psi:\reals^2\to\reals$, we have
\begin{align}
\plim_{n\to\infty}\frac{1}{n}\sum_{i=1}^{n}\psi(\theta_i,z_{t,i}^k)=
\E\big[\psi(\Theta,\gamma_t^k \Theta + ({\gamma_t^k})^{1/2} G)\big]\, ,
\end{align}
where expectation is with respect to $(\Theta,G)\sim\pi_{\Theta}\otimes\normal(0,1)$.
Here, we follow the procedure stated in \cref{alg:Spiked-Sampling-AMP} to take the spectral initialization for Bayes AMP.  
\end{proposition}

Building on state evolution, we prove Bayes optimality when the signal-to-noise
ratio is above a suitable constant threshold $\beta_0$.
	Here it is understood that the spectral initialization is selected
	 following the procedure stated in \cref{alg:Spiked-Sampling-AMP}. 
	 Further, we recall that 
	 $\hbm^{k}(\by, t)= \hat{\bm}^k_t$ is the output of AMP algorithm 
	\eqref{eq:general-AMP} at time $t$ after $k$ iterations, and $\bm(\by, t)$ is the mean 
	vector corresponding to the tilted measure as defined in Section \ref{sec:Tilted}.
%
\begin{lemma}\label{lemma:AMP-optimal}
Assume $\pi_{\Theta}$ has unit second moment and bounded eighth moment. 
Then there exists a  constant $\beta_0(\pi_{\Theta})$ that depends uniquely on $\pi_{\Theta}$, such that 
	the following  holds.
For any $\beta >\beta_0(\pi_{\Theta})$, the first stationary point of $\gamma \mapsto \Phi(\gamma, \beta, t)$  
on $(t,\infty)$ is also the unique global minimum of the same function over $\gamma\in (t, \infty)$ (recall this is defined in \cref{eq:Phi}).
In addition, for any $\eps, T>0$, 
there exists $K(\beta,T,\eps) \in \NN_{>0}$ that depends only on $(\beta, T, \eps)$, 
 such that for any $t \in [0, T]$:
	\begin{align*}
\plimsup\limits_{n \to \infty}\frac{1}{\sqrt{n}} \|\bm(\by(t), t) - \hbm^{K(\beta,T,\eps)}(\by(t), t)\|_2 \leq \eps,
	\end{align*}
Further, under the same conditions, we have 
\begin{align*}
	q_{\beta, t} = \lim_{n \to \infty}\frac{1}{n} 
	\E\left[ \|\E[\btheta \mid \bX, \by(t)]\|_2^2 \right] , 
\end{align*}
where $q_{\beta, t} = \E[\E[\Theta \mid \gamma_{\beta, t} \Theta + \sqrt{\gamma_{\beta, t}} G]^2]$,
and  $\gamma_{\beta, t}$ is the unique global minimum of $\gamma \mapsto \Phi(\gamma, \beta, t)$ on $(t, \infty)$.

\end{lemma}
The proof of \cref{lemma:AMP-optimal} is deferred to \cref{sec:proof-of-lemma:AMP-optimal}. 

\begin{remark}\label{rmk:Bayes}
Denote by $\cuE^{(1)}_{\beta, L, \Delta, \eps, n}$ the event 
that AMP returns an accurate approximation of the posterior mean 
for all $t\in\{0,\Delta,\dots, L\Delta\}$. Namely, we define:
\begin{align}\label{eq:E1}
\cuE^{(1)}_{\beta, L, \Delta, \eps, n} := \left\{  \frac{1}{\sqrt{n}} \big\| \bm(\by(\ell\Delta), \ell \Delta ) 
- \hbm^{K(\beta,T,\eps)}(\by(\ell\Delta), \ell \Delta) \big\|_2 \leq \eps\;\;
\forall\ell \in \{0,1,\cdots, L - 1\}  \right\}.
\end{align} 
By \cref{lemma:AMP-optimal}, we have $\prob(\cuE^{(1)}_{\beta, L, \Delta, \eps, n})=1-o_n(1)$.
\end{remark}

\subsubsection{Path regularity}\label{sec:path-regularity}

We next consider 
assumption \hyperlink{A2}{$\mathsf{(A2)}$} of Theorem \ref{thm:gen}:
Namely, we show that the path $t \mapsto \bm(\by(t), t)$ is
regular.  
\begin{lemma}\label{lemma:path-regularity}
Assume $\pi_{\Theta}$ has unit second moment and bounded eighth moment.
Then there exists a  constant
	 $\beta_0(\pi_{\Theta})$ that depends only on $\pi_{\Theta}$, 
	 such that the following holds:
For fixed $\beta \geq \beta_0(\pi_{\Theta})$ and $T \in \RR_{> 0}$,  
	   there exists a constant $C_{\sreg}= C_{\sreg}(\beta)>1 $ that depends only on $\beta$, 
	   such that for all $0 \leq t_1 < t_2 \leq T$, 
	\begin{align*}
		\plim\limits_{n \to \infty} \sup\limits_{t \in [t_1, t_2]} \frac{1}{n} 
		\big\|\bm(\by(t), t) - \bm(\by(t_1), t_1)\big\|_2^2 & = \plim\limits_{n \to \infty}\frac{1}{n} \big\|\bm( \by(t_1), t_1) - \bm(\by(t_2), t_2)\big\|_2^2 \\
		& \leq \frac{C_{\sreg}}{2}\cdot |t_1 - t_2|. 
	\end{align*}
\end{lemma}
The proof of \cref{lemma:path-regularity} is deferred to Appendix \ref{sec:proof-of-lemma:path-regularity}. 

\begin{remark}\label{rmk:Path}
Define $\cuE_{\beta, L, \Delta, n}^{(2)}$ as the following event: 
%
\begin{align}\label{eq:E3}
\cuE_{\beta, L, \Delta, n}^{(2)} := \left\{ \sup\limits_{t \in [ \ell\Delta, (\ell + 1)\Delta]} 
\frac{1}{\sqrt{n}} \big\|\bm( \by(t),t) - \bm(\by(\ell\Delta), \ell\Delta)\big\|_2 \leq C_{\sreg} \sqrt{\Delta} 
\;\; \forall\ell \in \{0\} \cup [L - 1]\right\}.
\end{align}
\cref{lemma:path-regularity} implies that $\P(\cuE_{\beta, L, \Delta, n}^{(2)}) 
=1-o_n(1)$. 
\end{remark}

\subsubsection{AMP is Lipschitz continuous} \label{sec:AMP-stable}

The crucial technical step is to prove that Bayes AMP is Lipschitz continuous in a neighborhood
of $\by(t)$, thus establishing Assumption \hyperlink{A3}{$\mathsf{(A3)}$} of Theorem \ref{thm:gen}.
 To be specific, we will use a change of variables technique introduced in \cite{celentano2021local}
 to prove that Bayes AMP is contractive in the new variables, for $\beta$  above certain threshold.
 
To define the change of variables technique, 
for $\gamma > 0$, we define $\Gamma_\gamma,\Psi_\gamma: \RR \to \RR$ by
\begin{align}\label{eq:transforms}
\begin{split}
	& \Gamma_{\gamma}(h) := \int_0^h \Var[\Theta \mid \gamma \Theta + \sqrt{\gamma} G = s]^{1/2} \dd s, \\
	& \Psi_{\gamma}(p) := \E[\Theta \mid \gamma \Theta + \sqrt{\gamma} G = \Gamma_{\gamma}^{-1}(p)], 
\end{split}
\end{align}
where $(\Theta, G) \sim \pi_{\Theta} \otimes \normal(0,1)$. 
Note that $\Gamma_{\gamma}, \Psi_{\gamma}$ are both strictly increasing. 
The mappings $\Gamma_{\gamma},\Psi_{\gamma}$ are specifically designed such that if we let $p = \Gamma_{\gamma}(h)$ and $m =  \Psi_{\gamma}(p)$, then the following factorization equality holds:
\begin{align*}
	\Var[\Theta \mid \gamma \Theta + \sqrt{\gamma} G = h]^{1/2} = \frac{\dd m}{\dd p} = \frac{\dd p}{\dd h}.
\end{align*}
One can verify that both $\Gamma_{\gamma}$ and $\Psi_{\gamma}$ are strictly increasing.
Furthermore, both $\Gamma_{\gamma}$ and $\Psi_{\gamma}$ are 
$M_{\Theta}$-Lipschitz continuous, where
 $M_{\Theta} = \|\pi_{\Theta}\|_{\infty}:=\sup(|\theta|:\; \theta\in\supp(\pi_{\Theta}))$.

Recall the posterior expectation function $\sF$ is defined in \cref{eq:F}.
For $t \in \RR_{> 0}$, $k \in \NN_+$, we define the AMP mapping 
 $T_{\AMP}^{(t,k)}: \RR^{n}\times\RR^n\times\RR^n \to \RR^n$ via
%
\begin{align*}
	T_{\AMP}^{(t,k)}(\bm, \overline{\bm}, \by) := \sF(\beta \bX \bm + \by - b_t^k 
	\overline{\bm},  \gamma_t^{k + 1}) \, .
\end{align*}
The AMP iteration \cref{eq:general-AMP} can therefore be rewritten
as
\begin{align*}
\hbm_t^{k + 1} = T_{\AMP}^{(t,k)}(\hbm_t^k, \hbm_t^{k - 1}, \by(t))\, .
\end{align*}
Let $\hat{\bp}_t^k := \Psi_{\alpha_t^k}^{-1}(\hat{\bm}_t^k)$ and define the AMP
 mapping in $\bp$-domain by
\begin{align*}
	\tilde T_{\AMP}^{(t,k)}(\bp, \overline{\bp}, \by) = \Psi^{-1}_{\gamma_t^{k + 1}}(
	\sF(\beta \bX \Psi_{\gamma_t^k}(\bp) + \by - b_t^k \Psi_{\gamma_t^{k - 1 }}(
	\overline{\bp}),  \gamma_t^{k + 1})) .
\end{align*} 
We immediately see that the vectors $\hat{\bp}_t^k$ satisfy the recursion
\begin{align*}
\hat\bp_t^{k + 1} = \tilde T_{\AMP}^{(t,k)}(\hat\bp_t^k, \hat\bp_t^{k - 1}, \by(t)).
\end{align*}
Note that the range of $\Psi_{\gamma}$ is $(a_{\Theta},b_{\Theta})$, where
$a_{\Theta} = \inf\,{\rm supp}(\pi_{\Theta})$, $b_{\Theta} = \sup\, {\rm supp}(\pi_{\Theta})$.
 For $\bm \in (a_{\Theta}, b_{\Theta})^n$, we define 
 $\bD_{\gamma}(\bm) := \diag\{\Var[\, \Theta \mid \gamma \Theta + \sqrt{\gamma}G =
  \Gamma_{\gamma}^{-1}(\Psi_{\gamma}^{-1}(\bm))]^{1/2}\} \in \RR^{n \times n}$, where,
  by convention,  the conditional variance operator applies on vectors entrywise. 
  
  We compute the Jacobian matrices of the AMP  mappings $T_{\AMP}^{(t,k)}$ and 
  $\tilde T_{\AMP}^{(t,k)}$:
\begin{align}
	& \frac{\dd T_{\AMP}^{(t,k)}(\bm, \bm^-, \by)}{\dd (\bm, \bm^-, \by)} =
	 \left( \beta \bD_{ \gamma_t^{k + 1}}(\bm^+)^2 \bX;\; - b_t^k \bD_{\gamma_t^{k + 1}}(\bm^+)^2;\; \bD_{\gamma_t^{k + 1}}(\bm^+)^2 \right), \label{eq:general-jacob-m} \\
	 &\frac{\dd \tilde T_{\AMP}^{(t,k)}(\bp, \obp, \by)}{\dd (\bp, \obp, \by)} = \left( \beta \bD_{ \gamma_t^{k + 1}}(\bm^+) \bX \bD_{\gamma_t^{k}}(\bm);\; -b_t^k \bD_{\gamma_t^{k + 1}}(\bm^+) \bD_{\gamma_t^{k - 1}}(\bm^-), \bD_{\gamma_t^{k + 1}}(\bm^{+}) \right). \label{eq:general-jacob-p}
\end{align}
(Here it is understood that $\bm^{+} = T_{\AMP}^{(t,k)}(\bm, \bm^{-}, \by)$,
$\bm=\Psi_{\gamma_t^k}(\bp)$, and $\bm^-=\Psi_{\gamma_t^{k-1}}(\bp^-)$.) 

Roughly speaking, we will show that, if the signal strength $\beta$ is large, and 
after a large number of iterations $k$, then 
 most elements $\bD_{\gamma_t^{k}}(\hbm^t_k)$, 
$\bD_{\gamma_t^{k + 1}}(\hbm^{k+1}_t)$, $\bD_{\gamma_t^{k - 1}}(\hat\bm_t^{k-1})$ become small. 
This in turn will imply that the operator norms of the Jacobian matrices 
in \cref{eq:general-jacob-m,eq:general-jacob-p} are small.
Finally, this can be used to prove that the AMP  mapping is contractive. 

The next two lemmas formalize this argument. In the first lemma, we provide an upper
 bound on $\|\bD_{ \gamma_{l \delta }^{k}}(\hat{\bm}_{l\delta}^k)\|_F^2 / n$ 
 with a function of $\beta$. In the same lemma, we also show that AMP is 
  with high probability Lipschitz continuous if we allow the Lipschitz constant to depend on $(\beta, \pi_{\Theta})$ and the number of iterations. 
\begin{lemma}\label{lemma:general-AMP-hits-strong-signal-region}
Assume $\pi_{\Theta}$ has unit second moment and is discrete.
Let $K(\beta,T,\eps)$ be the constant of Lemma \ref{lemma:AMP-optimal}.
Then there exist constants $\beta_0, C_{\sconv} > 0$ that 
 depend uniquely on $\pi_{\Theta}$, such that the following  hold:
 For all $\beta \geq \beta_0$, there exist 
 $k_0(\beta)  \in \NN_+$, $\Lipc(\beta) \in \RR_{> 0}$ which are  functions of $(\pi_{\Theta},\beta)$ only, 
 such that for all 
 $\eps>0$, $ t \in [0,T]$, the following hold with probability $1-o_n(1)$:
 %
 \begin{enumerate}
 \item For all  $k_0(\beta) \leq k \leq K(\beta, T, \eps)$, 
	\begin{align}
		& \frac{1}{n}\|\bD_{\gamma_{t }^{k}}(\hbm^k(\by(t),t))\|_F^2 \leq 
		C_{\sconv}^{-1}\exp(-C_{\sconv} \beta^2), \label{eq:general-m}\\
		& b_t^k \leq  C_{\sconv}^{-1}\exp(-C_{\sconv} \beta^2)\, . \label{eq:general-b}
	\end{align}
\item	For $k\in \{k_0(\beta) - 1, k_0(\beta), k_0(\beta) + 1\}$:
	\begin{align}
		&\sup_{\by_1\neq\by_2} \frac{\|\hbm^{k}(\by_1, t) - \hbm^{k}( \by_2, t)  \|_2}{\|\by_1 - \by_2\|_2} \leq
		\Lipc(\beta)\, , \label{eq:general-Lipschitz1} \\
		& \sup_{\by_1\neq\by_2}\frac{\|\hbp^{k}(\by_1, t) - \hbp^{k}( \by_2, t)  \|_2}{\|\by_1 - \by_2\|_2} \leq
		\Lipc(\beta)\, , \label{eq:general-Lipschitz2}
	\end{align}
	where $\hbp^k(\by, t) := \Psi^{-1}_{\gamma_t^k}(\hbm^k(\by, t))$.
	\end{enumerate}
\end{lemma}
We postpone the proof of \cref{lemma:general-AMP-hits-strong-signal-region} to Appendix 
\ref{sec:proof-of-lemma:general-AMP-hits-strong-signal-region}. 

 By classical estimates on the norm of spiked random matrices \cite{benaych2011eigenvalues}, with probability
  $1 - o_n(1)$ we have 
 $\|\bX\|_{\op} \leq \beta + \beta^{-1} + 1$.
 We denote by $\cuE^{(3)}_{\beta, L, \Delta, \eps, n}$ the intersection
 of this event and the one of \cref{lemma:general-AMP-hits-strong-signal-region}.
 Namely, we define
\begin{align}
	\cuE^{(3)}_{\beta, L, \Delta, \eps, n} := & 
		\Big\{ \mbox{Eq.~\eqref{eq:general-m} holds for all $k_0(\beta) \leq k \leq K(\beta, L\Delta,  \eps)$
		and all $t / \Delta \in \{0\} \cup [L]$}, \nonumber
		 \\
 &  \mbox{ \cref{eq:general-Lipschitz1} and \cref{eq:general-Lipschitz2} hold for all $k
 \in \{k_0(\beta), k_0(\beta) \pm 1\}$ 
	and all $t / \Delta \in \{0\} \cup [L]$ }, \nonumber\\
	& \mbox{and }\|\bX\|_{\op} \leq 1+\beta + \beta^{-1}  \Big\}.\label{eq:E2}
\end{align}
By the last lemma and a union bound, we have $\P(\cuE^{(3)}_{\beta, L, \Delta, \eps, n})=1-o_n(1)$. In what follows, we will be mainly working on the set $\cuE^{(1)}_{\beta, L, \Delta, \eps, n} \cap \cuE_{\beta, L, \Delta, n}^{(2)} \cap \cuE^{(3)}_{\beta, L, \Delta, \eps, n}$, which occurs with probability $1 - o_n(1)$ by the lemmas we establish. 

The next lemma from \cite{celentano2021local} is useful for bounding the operator norms of the 
Jacobian matrices. 
\begin{lemma}[Lemma C.2. in \cite{celentano2021local}]\label{lemma:C2}
	For $\bt \in [0,1]^n$ and $\xi > 0$, denote by $S(\bt, \xi)$ the subset of indices $i \in \{1,\cdots, n\}$ for which $t_i \geq \xi$. Then there exist universal constants $C, C', c > 0$ such that for $\bW \sim \GOE(n)$, any $\xi > 0$ and $0 < q < 1$, 
	\begin{align}\label{eq:C2}
		\P\left( \sup_{
		\substack{\bt_1, \bt_2 \in [0,1]^n:\\ |S(\bt_1, \xi)|\vee |S(\bt_2, \xi)| \leq nq}} \|\diag(\bt_1) \bW \diag(\bt_2)\|_{\op} \geq C'(\xi + \sqrt{q \log(e / q)}) \right) \leq Ce^{-cqn}. 
	\end{align}
\end{lemma} 
\cref{lemma:general-AMP-hits-strong-signal-region,lemma:C2} together
imply that AMP is a contraction, whence Lipschitz, in a neighborhood of $\by(t)$. 
\begin{lemma}\label{lemma:MainLip}
We assume the assumptions of \cref{lemma:general-AMP-hits-strong-signal-region},
and let $k_0(\beta), K(\beta,T,\eps)$ be as defined there.
Then there exists  $\beta_0 > 0$ that 
 depends uniquely on $\pi_{\Theta}$, such that the following  hold: 
 For all $\beta \geq \beta_0$,
 there exists $r(\beta), \Lips(\beta)> 0$ depending uniquely on $(\pi_{\Theta},\beta)$ such that, for all $t\in[0,T]$, 
 the following holds with probability $1-o_n(1)$ for all $k_0(\beta)\le k\le K(\beta,T,\eps)$:
 \begin{align}
 \sup_{\by_1\neq\by_2\in \Ball^n(\by(t),r(\beta))}
 \frac{\|\hbm^k(\by_1,t)-\hbm^k(\by_2,t)\|_2}{\|\by_1-\by_2\|_2}\le 2\, \Lips(\beta)\, .
 \end{align}
 (Here, $\Ball^n(\bx_0;r):=\{\bx\in\RR^n:\;\|\bx-\bx_0\|\le r\}$.)
\end{lemma}
The proof of \cref{lemma:MainLip} can be found in Appendix \ref{sec:proof-of-lemma:main}.

\subsubsection{Completing the proof of \cref{thm:main2}}\label{sec:proof}

We are now in position to apply Theorem \ref{thm:gen} to prove \cref{thm:main2}.  
First of all, notice that in the present
case $\bH=\id_n$ and therefore $\bm(\by,t)=\bm_{\btheta}(\by,t)$. We 
set  $T$, $\eps$ and $\delta$ as follows
\begin{align}
T&=\frac{4}{\xi^2}\,,\\
 \eps &= \frac{r(\beta)\wedge \xi}{8}\, e^{-8\Lips(\beta)/\xi^2}\, ,\label{eq:ChoiceEps}\\
 \sqrt{\Delta}& = \frac{r(\beta)\wedge\xi}{8C_{\sreg}} \, e^{-8\Lips(\beta)/\xi^2}\, .
 \label{eq:ChoiceDelta}
\end{align}
and set $K_{\AMP} = K(\beta,T,\eps)$, where $K(\beta,T,\eps)$ is defined  by Lemma \ref{lemma:AMP-optimal} and $\Lips(\beta)$ is defined by \cref{lemma:MainLip}.

We next check that
assumptions \hyperlink{A1}{$\mathsf{(A1)}$}, \hyperlink{A2}{$\mathsf{(A2)}$} and 
\hyperlink{A3}{$\mathsf{(A3)}$} hold (with $\eta=o_n(1)$):
\begin{itemize}
\item[${\sf (A1)}$] By \cref{rmk:Bayes}, this assumption holds with $\eps_1=\eps$.
\item[${\sf (A2)}$] By \cref{rmk:Path}, this assumption holds with $C_1=C_{\sreg}$,
$\eps_2=0$.
\item[${\sf (A3)}$] By \cref{lemma:MainLip}, this assumption holds with $C_2=2\Lips(\beta)$,
$r_{\ell}=r(\beta)$ and $\eps_3 = 0$. We need to check the lower bound on $r_{\ell}$ that is required by
assumption \hyperlink{A3}{$\mathsf{(A3)}$}.
For that purpose, note that 
\begin{align}
(C_1\sqrt{\Delta}+\eps_1+\eps_2 + \eps_3)\frac{e^{C_2L\Delta}}{C_2}& \le 
\big(C_{\sreg}\sqrt{\Delta}+\eps\big) e^{8\Lips(\beta)/\xi^2}\\
& \le \frac{r(\beta)}{2} <r_{\ell}\, ,
\end{align}
where in the first step we used $L\Delta= T= 2/\xi$ and, without loss of generality,
$\Lips(\beta)\ge 1$. In the second inequality, we used 
the choices for $\eps$, $\Delta$ given in Eqs.~\eqref{eq:ChoiceEps} and
\eqref{eq:ChoiceDelta}.
\end{itemize}
Note that, since $\|\pi_{\Theta}\|_{\infty} < \infty$, we then have $\int (\|\btheta\|^2/n)^2 \mu_{\bX}(\de\btheta)\le R^4$, where $R > 0$ is a constant depending only on $\pi_{\Theta}$.  

Applying \cref{eq:SecondGeneral} from Theorem \ref{thm:gen}, we obtain that, 
for any $\eta>0$, the following holds with probability
$1-o_n(1)$ with respect to the choice of $\bX$
	\begin{align}
 W_{2, n}(\mu_{\bX}, \mu_{\bX}^{\salg})
		&\leq  \eps + \big(C_{\sreg}\sqrt{\Delta}+\eps\big) e^{8\Lips(\beta)/\xi^2}
		+C R\eta^{1/2}
		+\frac{1}{\sqrt T} \\
		& \le \frac{7\xi}{8}+C R\eta^{1/2} \le \frac{9}{10}\xi\, .
	\end{align}
where $C > 0$ is a numerical constant, and the last inequality follows by choosing a suitably small $\eta$.  

\section{Proofs for \cref{alg:Spiked-Sampling-continuous}}
\label{sec:proof-alg:Spiked-Sampling-continuous}

\subsection{Proof of \cref{thm:spiked-continuous}}
\label{sec:spiked-continious}


For the sake of simplicity, we consider the case of non-symmetric $\pi_{\Theta}$. 
The case of symmetric $\pi_{\Theta}$ can be handled by applying random 
sign flips at the end of the sample generation process, as discussed in \cref{sec:discrete-spiked-sampling}. 
We skip this discussion here to avoid duplication.

Throughout this section we use the shorthand  $\bm_{\sB} = \bm(\by(t),t)$
for the posterior expectation $\bm(\by(t),t):= \E[\btheta \mid \bX, \by(t)]$. 
In the next lemma, we show that $\bm_{\sB}$ is an approximate stationary point of the TAP free energy.
 Its proof is given in Appendix \ref{sec:proof-lemma:approximate-stationary}. 
\begin{lemma}
\label{lemma:approximate-stationary}
	Assume the conditions of Theorem \ref{thm:spiked-continuous}, and $\beta \geq \beta_0(\pi_{\Theta})$ where $\beta_0(\pi_{\Theta})$ is defined in \cref{lemma:large-beta-cond1}. 
	Then, $\bm_{\sB}$ is an approximate stationary point of the TAP free energy, in the sense that 
	\begin{align*}
		\plim_{n\to \infty}\frac{1}{n} \|\nabla \cuF_{\sTAP} (\bm_{\sB}; \bX, \by(t), \beta, t)\|_2^2 = 0\, . 
	\end{align*}
\end{lemma}
We then show that the TAP free energy is strongly convex in a neighborhood of $\bm_{\sB}$,
 if restricted to the sphere $\Spher_{\beta, t}$. We delay the proof of the
 next  lemma to Appendix \ref{sec:proof-lemma:projected-hessian}. 
\begin{lemma}
\label{lemma:projected-hessian}
	Assume the conditions of  Theorem \ref{thm:spiked-continuous}.
	Then there exists a constant $\beta_0(\pi_{\Theta})$  that depends only on $\pi_{\Theta}$,
	such that the following holds for all $\beta\ge \beta_0(\pi_{\Theta})$.
	There exists $r, \kappa > 0$ that depend only on $\beta$, such that for any $t\ge 0$, with probability $1 - o_n(1)$ the following happens. 
	For all $\bm \in \mathsf{B}(\bm_{\sB}, \sqrt{n} r)$, $\by \in \RR^n$ and $\langle \bx, \bm \rangle = 0$, 
	\begin{align*}
	\bx^{\top}  \nabla^2_{\bm} \cuF_{\sTAP} (\bm; \bX, \by, \beta, t)  \bx \succeq  \kappa \|\bx\|_2^2 \, . 
	\end{align*}
\end{lemma}
 
Next, we prove that for a sufficiently large $k$, the mapping  
$\bw \mapsto \cuF_{\sTAP}(\bphi_{\hat \bm_t^k}(\bw); \bX, \by, \beta, t)$ 
 is with high probability strongly convex in a small neighborhood of the origin, where we recall that $\bphi_{\bm}$ 
 is defined in \cref{eq:varphi-project}, and $\hat \bm_t^k$ is the $k$-th AMP iterate as defined in \cref{eq:general-AMP}. 
\begin{lemma}
\label{lemma:tangent-space-convextiy}
	Under the conditions of  Theorem \ref{thm:spiked-continuous},
	 there exists a constant $\beta_0(\pi_{\Theta})$ that depends only on $\pi_{\Theta}$,
	such that, for all $\beta\ge \beta_0(\pi_{\Theta})$, the following holds.
	There exist $k \in \NN_{+}$, $c_1, c_2 > 0$ and $R, R_y \in (0, \sqrt{q_{\beta, t}} / 10)$ with $R_y \leq c_1 R / 100$
	 that depend only on $(\beta, \pi_{\Theta})$, such that 
	  the following statement is true with probability $1 - o_n(1)$:
	  For all $\by\in\Ball^{n}(\by(t),R_y\sqrt{n})$,
	   the mapping $\bw \mapsto \cuF_{\sTAP}(\bphi_{\hat \bm_t^k}(\bw); \bX, \by, \beta, t)$
	    is $c_1$-strongly convex and $c_2$-smooth on $\bw \in \Ball^{n - 1} (\bfzero,R\sqrt{n})$.	    
	   In addition, we denote the minimizer of 
	$\bw \mapsto \cuF_{\sTAP}(\bphi_{\hat \bm_t^k}(\bw); \bX, \by, \beta, t)$ 
on	$\bw \in \Ball^{n - 1} (\bfzero,R\sqrt{n})$ by $\bw^{\ast}(\by)$. Then $\bw^{\ast}(\by)$ exists for all $\by \in \Ball^n(\by(t), R_y \sqrt{n})$ and satisfies
	$\|\bw^{\ast}(\by)\|_2 < R \sqrt{n} / 4$. 
	   %
	    	    
\end{lemma}
We prove \cref{lemma:tangent-space-convextiy} in Appendix 
\ref{sec:proof-lemma:tangent-space-convextiy}.

Finally, we apply \cref{thm:gen} to prove \cref{thm:spiked-continuous}. 
To this end, it suffices to verify \hyperlink{A1}{$\mathsf{(A1)}$} and  \hyperlink{A3}{$\mathsf{(A3)}$}, 
as \hyperlink{A2}{$\mathsf{(A2)}$} has already been checked in Remark \ref{rmk:Path}. 
In this case, $C_1 = C_{\sreg}$ and $\eps_2 = 0$. 
To verify \hyperlink{A1}{$\mathsf{(A1)}$}, we prove the following lemma. We defer the proof of the lemma to Appendix \ref{sec:proof-lemma:A1-4}. 
\begin{lemma}
\label{lemma:A1-4}
	We assume the conditions of \cref{thm:spiked-continuous}. Then for any $\eps > 0$, there exists a sufficiently large $k$, such that with probability $1 - o_n(1)$ 
	\begin{align*}
		\frac{1}{\sqrt{n}}  \|\bphi_{\hat \bm_t^k}(\bw^{\ast}(\by(t))) - \bm_{\sB}\|_2 \leq \eps\, . 
	\end{align*}
\end{lemma}


By \cref{lemma:tangent-space-convextiy} we know that the mapping $\bw \mapsto \cuF_{\sTAP}(\bphi_{\hat \bm_t^k}(\bw); \bX, \by(t), \beta, t)$ is with high probability  strongly convex and smooth
for	$\bw \in \Ball^{n - 1} (\bfzero,R\sqrt{n})$ and has a minimizer $\bw^{\ast}(\by(t))$. 
Therefore, in this case we can find an arbitrarily accurate estimate to $\bw_{\ast}(\by(t))$ by running gradient descent with an appropriate step size. 
\hyperlink{A1}{$\mathsf{(A1)}$} then follows from \cref{lemma:A1-4} with an arbitrarily small $\eps_1$. 
In the sequel, we verify \hyperlink{A3}{$\mathsf{(A3)}$}.

We finally prove \hyperlink{A3}{$\mathsf{(A3)}$}. 
For $\by_1, \by_2 \in \{\by \in \RR^n: \|\by - \by(t)\|_2 \leq R_y \sqrt{n}\}$, with probability $1 - o_n(1)$ it holds that 
\begin{align*}
	0 \geq & \cuF_{\sTAP}(\bphi_{\hat \bm_t^k}(\bw_{\ast}(\by_1));  \bX, \by_1, \beta, t) - \cuF_{\sTAP}(\bphi_{\hat \bm_t^k}(\bw_{\ast}(\by_2)); \bX, \by_1, \beta, t) \\
	\geq & - \| \nabla_{\bw} \cuF_{\sTAP}(\bphi_{\hat \bm_t^k}(\bw); t, \bX, \by_1, \beta) |_{\bw = \bw_{\ast}(\by_2)} \| \cdot \|\bw_{\ast}(\by_1) - \bw_{\ast}(\by_2)\|_2 + \frac{c_1}{8} \|\bw_{\ast}(\by_1) - \bw_{\ast}(\by_2)\|_2^2. 
\end{align*}
Note that 
\begin{align*}
	\nabla_{\bw} \cuF_{\sTAP}(\bphi_{\hat \bm_t^k}(\bw); \bX, \by_1, \beta, t) |_{\bw = \bw_{\ast}(\by_2)} = -\nabla_{\bw} \bphi_{\hat \bm_t^k}(\bw)  (\by_1 - \by_2)\mid_{\bw = \bw_{\ast}(\by_2)}. 
\end{align*}
%
%
%
Checking the proof of \cref{lemma:tangent-space-convextiy}, we see that $\bw \mapsto \bphi_{\bm}(\bw)$ is with high probability $5$-Lipschitz continuous. 
As a consequence, $\|\bphi_{\hat \bm_t^k}(\bw_{\ast}(\by_1)) - \bphi_{\hat \bm_t^k}(\bw_{\ast}(\by_2))\|_2 \leq 40 \|\by_1 - \by_2\|_2 / c_1$.
The proof is complete as we can employ gradient descent to find the approximate minimizers $\bw_{\ast}(\by)$, thanks to strong convexity. 
Namely, we have verified \hyperlink{A3}{$\mathsf{(A3)}$} with $\eps_3 = o_n(1)$ and $C_2 = 40 / c_1$.

\subsection{Proof of \cref{lemma:approximate-stationary}}
\label{sec:proof-lemma:approximate-stationary}

We let $g_m(\lambda) = \lambda m - \phi(\lambda, w)$, then $g_m'(\lambda) = m - \E[\Theta \mid \sqrt{\omega} \Theta + G = \lambda \omega^{-1/2}]$ and $g_m''(\lambda) = - \Var[\Theta \mid \sqrt{\omega} \Theta + G = \lambda \omega^{-1/2}] < 0$. Hence, $g_m$ is strongly concave, and has a unique maximizer. 
We denote by $\lambda(m)$ its maximizer, which satisfies 
\begin{align*}
	m = \E[\Theta \mid \sqrt{\omega} \Theta + G = \lambda (m) \omega^{-1/2}]. 
\end{align*} 
The above equality further implies that $\frac{\partial m}{\partial \lambda(m)} = \Var[\Theta \mid \omega \Theta + \omega^{1/2} G = \lambda(m)]$. 
By the envelope theorem, 
\begin{align*}
	\frac{\partial}{\partial m} h(m, \beta^2 q_{\beta, t})  = \lambda(m).  
\end{align*}
We then take the gradient of the TAP free energy, which gives  
\begin{align*}
	\nabla \cuF_{\sTAP} (\bm; \bX, \by(t), \beta, t) = - \beta \bX \bm + \beta^2 (1 - q_{\beta, t}) \bm - \by(t) + \lambda(\bm). 
\end{align*}
Recall that $\hat \bm_t^k$ is the $k$-th iterate of the Bayes AMP algorithm given in \cref{eq:general-AMP}. 
By Lemma \ref{lemma:AMP-optimal} it holds that for any $\eps > 0$, 
there exists $k_{\eps} \in \NN_+$, such that for any fixed $k\ge k_{\eps}$, with probability $1 - o_n(1)$ 
\begin{align}\label{eq:B-k}
	 \|\hat \bm_t^{k_{\eps}} - \bm_{\sB}\|_2^2 / n \leq \eps. 
\end{align}
Plugging this fact into the AMP iteration \eqref{eq:general-AMP}, we obtain that for $k\ge k_{\eps}+1$,
\begin{align}
\frac{1}{n} \|- \beta \bX \hat \bm_t^{k} + \beta^2 (1 - q_{\beta, t}) 
\hat \bm_t^{k} - \by(t) + \lambda(\hat \bm_t^{k})\|_2^2 \leq C
\eps + o_P(1). 
\end{align}
By random matrix theory, with probability $1 - o_n(1)$ we have $\|\bX\|_{\op} \leq C_{\delta}$ 
for some positive constant $C_{\delta}$ that is a function of $\delta$ only. 
Also, note that $\eps$ above is an arbitrary positive number. Hence, 
using again Eq.~\eqref{eq:B-k},
\begin{align*}
\plim_{n\to\infty}\	\frac{1}{n} \| \nabla \cuF_{\sTAP} (\bm_{\sB}; \bX, \by(t), \beta, t)\|_2^2 = 0. 
\end{align*}  
The proof is complete.

\subsection{Proof of \cref{lemma:projected-hessian}}
\label{sec:proof-lemma:projected-hessian}

We take the projected Hessian of $\cuF_{\sTAP}$: 
\begin{align*}
	\sP_{\bm}^{\perp} \nabla^2 \cuF_{\sTAP} (\bm;\bX, \by, \beta, t) 
	\sP_{\bm}^{\perp} = - \frac{\beta^2}{n} \mathsf{P}_{\bm}^{\perp} \btheta \btheta^{\top} \mathsf{P}_{\bm}^{\perp} 
	- \beta \mathsf{P}_{\bm}^{\perp} \bW \mathsf{P}_{\bm}^{\perp} + \beta^2 (1 - q_{\beta, t}) \mathsf{P}_{\bm}^{\perp} + \mathsf{P}_{\bm}^{\perp} \bD(\bm) \mathsf{P}_{\bm}^{\perp}, 
\end{align*}
where $\bD(\bm) = \diag ( \{ \Var[\Theta \mid \beta^2 q_{\beta, t} \Theta + \beta q_{\beta, t}^{1/2} G = \lambda(m_i)]^{-1} \}_{i \in [n]} )$.
We recall that $\lambda(m)$ is defined in \cref{sec:proof-lemma:approximate-stationary}. 
By the Brascamp–Lieb inequality, we see that the diagonal entries of $\bD(\bm)$ are no smaller than $\beta^2 q_{\beta, t} - \|U''\|_{\infty}$, which is no smaller than $\beta^2 / 4$ for a large enough $\beta$. 
 
Observe that 
\begin{align*}
	\plim_{n\to\infty}\frac{1}{n}\| \mathsf{P}_{\bm_{\sB}}^{\perp} \btheta  \|_2^2 =
	1 - q_{\beta, t} \leq \beta^{-2},    
\end{align*}
%
where the second inequality is by \cref{lemma:large-beta-cond1}. 
We can then choose $r$ small enough, 
such that with high probability, for all $\bm \in \mathsf{B}(\bm_{\sB}, \sqrt{n} r)$, it holds that 
\begin{align*}
	\frac{1}{n}\| \mathsf{P}_{\bm}^{\perp} \btheta\|_2^2  \leq  2\beta^{-2}. 
\end{align*}
Recall that $\bW \sim \GOE(n)$, hence with probability $1 - o_n(1)$ we have
 $\| \beta \mathsf{P}_{\bm}^{\perp} \bW \mathsf{P}_{\bm}^{\perp} \|_{\op} \leq 3 \beta$. 
Putting together the above arguments, we conclude that for all $\bx \in \RR^d$ that satisfies $\langle \bx, \bm \rangle = 0$, we have 
\begin{align*}
	\bx^{\top} \mathsf{P}_{\bm}^{\perp} \nabla \cuF_{\sTAP} (\bm; \bX, \by, \beta, t) \mathsf{P}_{\bm}^{\perp} \bx \geq (\beta^2 / 4 - 3 \beta - 2) \|\bx\|_2^2. 
\end{align*}
The proof is complete if we choose a sufficiently large $\beta_0(\pi_{\Theta})$.

\subsection{Proof of \cref{lemma:tangent-space-convextiy}}
\label{sec:proof-lemma:tangent-space-convextiy}

In this proof, we write $\bm = \hat \bm_t^k$ for the simplicity of presentation.
We emphasize that $\bm$ is random and implicitly depends on $k$. 

To prove $\bw \mapsto \cuF_{\sTAP}(\bphi_{\bm}(\bw);  \bX, \by, \beta, t)$ is strongly convex, we take the Hessian of this function: 
\begin{align}
	\label{eq:DecHess}& \nabla_{\bw}^2 \cuF_{\sTAP} (\bphi_{\bm}(\bw); \bX, \by, \beta, t) \nonumber\\
	 = & \underbrace{\bD_{\bw} \bphi_{\bm}(\bw)^{\top} \nabla_{\bm}^2 \cuF_{\sTAP} (\bphi_{\bm}(\bw); \bX, \by, \beta, t) \bD_{\bw} \bphi_{\bm}(\bw)}_{(i)} \\
&  +  \underbrace{\bD_{\bm} \cuF_{\sTAP} (\bphi_{\bm}(\bw); \bX, \by, \beta, t) \nabla_{\bw}^2 \bphi_{\bm}(\bw)}_{(ii)}. \nonumber 
\end{align}
In the above display, the Jacobian and the Hessian of $\bphi_{\bm}$ takes the following form: 
%
\begin{align*}
	\bD_{\bw} \bphi_{\bm}(\bw) &=  \frac{\sqrt{n q_{\beta, t}} \cdot \bT_{\bm}}{\|\bm + \bT_{\bm} \bw\|_2} -   \frac{\sqrt{n q_{\beta, t}} \cdot(\bm + \bT_{\bm} \bw)\otimes \bw^{\top} \bT_{\bm}^{\top} \bT_{\bm}}{\|\bm + \bT_{\bm} \bw\|_2^3} \in \RR^{n \times (n - 1)}, \\
	 \nabla_{\bw}^2 \bphi_{\bm}(\bw) &=
	 - \frac{\sqrt{n q_{\beta, t}}}{\|\bm + \bT_{\bm} \bw\|_2^3}
	 \big((\bT_{\bm} \otimes \bT_{\bm}^{\top} \bT_{\bm} \bw) + (\bT_{\bm} \otimes \bT_{\bm}^{\top} \bT_{\bm} \bw)^{(23)}\big)\\	 
	& -   \frac{\sqrt{n q_{\beta, t}} \cdot (\bm + \bT_{\bm} \bw) \otimes  \bT_{\bm}^{\top} \bT_{\bm}}{\|\bm + \bT_{\bm} \bw\|_2^3}  \\
	& + \frac{3\sqrt{n q_{\beta, t}} (\bm + \bT_{\bm} \bw) \otimes  \bT_{\bm}^{\top} \bT_{\bm} \bw \otimes  \bT_{\bm}^{\top} \bT_{\bm} \bw}{\|\bm + \bT_{\bm} \bw\|_2^5} \in \RR^{n \times (n - 1)^2}. 
\end{align*}
Here, the superscript $(23)$ denotes transposition with respect to second and third dimension.
For a sufficiently small $R$ and a sufficiently large $k$, we conclude that there exists $R_y > 0$, such that with probability 
$1 - o_n(1)$, it holds that 
%
\begin{align}
\|\bD_{\bw}^2 \bphi_{\bm}(\bw)\|_{\op} \leq  \frac{2}{\sqrt{n}}\, \label{eq:CurvatureSphere}
\end{align}
for all $\|\bw\|_2 \leq \sqrt{n}R$ and $\|\by - \by(t)\|_2 \leq \sqrt{n} R_y$. 
We recall that $\|\cdot\|_{\op}$ denotes the spectral norm of a tensor that is defined in \cref{sec:notation}. 

We next show that for a sufficiently large $k$ and a sufficiently small $R$, there exists $R_y > 0$ such that the following statement holds: 
There exists $\kappa_1 > 0$ that depends only on $\beta$, such that 
with probability $1-o_n(1)$, for all  $\bw \in \mathsf{B}(\mathbf{0}_{n - 1}, R \sqrt{n})$ and $\by \in \mathsf{B}(\by(t), R_y \sqrt{n})$,
\begin{align}
\label{eq:150}
	\nabla_{\bw}^2 \cuF_{\sTAP} (\bphi_{\bm}(\bw); \bX, \by, \beta, t)
	 \succeq \kappa_1 \id_{n-1}. 
\end{align}
We will separately analyze terms $(i)$ and $(ii)$ of the decomposition \eqref{eq:DecHess}. 
We will show that term $(i)$ dominates and is positive definite, 
while term $(ii)$ is negligible

We shall make use of \cref{lemma:projected-hessian} to lower bound $(i)$.
 By \cref{lemma:AMP-optimal}, 
 there exist a large enough $k$ and small enough 
 $R$, such that the following statement is true: There exists $R_y > 0$, such that  with probability $1 - o_n(1)$ we have 
$\|\bD_{\bw} \bphi_{\bm}(\bw)\|_{\op} \in [1/2, 2]$ and 
 $\|\bphi_{\bm}(\bw) - \bm_{\sB}\|_2 \leq \sqrt{n}r / 2$ (we recall that $r$ appears in 
 \cref{lemma:projected-hessian}) for all $\|\bw\|_2 \leq R \sqrt{n}$ and $\|\by - \by(t)\|_2 \leq R_y \sqrt{n}$. 
Then applying \cref{lemma:projected-hessian}, we see that for a sufficiently large $k$ and sufficiently small $R$, there exists $R_y > 0$, such that with probability $1 - o_n(1)$ the following statement is true:
For all $\|\bw\|_2 \leq R \sqrt{n}$ and $\|\by - \by(t)\|_2 \leq R_y \sqrt{n}$, it holds that 
\begin{align}
\label{eq:First-Part-Hessian}
\bD_{\bw} \bphi_{\bm}(\bw)^{\top} \nabla_{\bm}^2 \cuF_{\sTAP} (\bphi_{\bm}(\bw); \bX, \by, \beta, t) \bD_{\bw} \bphi_{\bm}(\bw) \succeq \kappa \bI_{n-1} 
\end{align}
for some $\kappa > 0$ that depends only on $\beta$. 
This lower bounds term $(i)$. 

We next proceed to upper bound term $(ii)$. To get such a result, we need to establish the 
following lemma, whose proof is deferred to Appendix \ref{sec:proof-lemma:lower-bound-var}.  
\begin{lemma}
\label{lemma:lower-bound-var}
	Under the conditions of \cref{lemma:tangent-space-convextiy}, there exists a positive constant $c_{\Theta, t}$ that depends only on $(\pi_{\Theta}, t)$, such that 
	\begin{align*}
		\Var[\Theta \mid t \Theta + \sqrt{t} G = x] \geq c_{\Theta, t}
	\end{align*}
	for all $x \in \RR^d$. In the above display, $(\Theta, G) \sim \pi_{\Theta} \otimes \normal(0,1)$. 
\end{lemma} 
We then upper bound term $(ii)$. Recall that  
\begin{align}
\label{eq:upper-bound-U-beta-t}
	\nabla_{\bm}^2 \cuF_{\sTAP} (\bm; \bX, \by, \beta, t) = - \frac{\beta^2}{n} \btheta \btheta^{\top}  - \beta \bW + \beta^2 (1 - q_{\beta, t}) \id_n +  \bD(\bm).  
\end{align}
Using \cref{eq:upper-bound-U-beta-t} and \cref{lemma:lower-bound-var}, we conclude that there exists a constant $U_{\beta, t} > 0$ that depends only on $(\beta, t, \pi_{\Theta})$, such that with probability $1 - o_n(1)$, for all $\by \in \RR^n$ it holds that 
\begin{align}
\label{eq:hessian-U-beta-t}
	\|\nabla^2_{\bm} \cuF_{\sTAP} (\bm; \bX, \by, \beta, t)\|_{\op} \leq U_{\beta, t}. 
\end{align}
Combining the above upper bound and \cref{lemma:approximate-stationary}, we conclude that for a sufficiently large $k$ and sufficiently small $R$, the following statement is true: 
There exists $R_y > 0$, 
such that with probability $1 - o_n(1)$, for all $\|\bw\|_2 \leq \sqrt{n}R$ and $\|\by - \by(t)\|_2 \leq \sqrt{n} R_y$, it holds that 
\begin{align}
\label{eq:kappa-100}
	\frac{1}{n} \|\nabla_{\bm} \cuF_{\sTAP} (\bphi_{\bm}(\bw); \bX, \by, \beta, t)\|_2^2 \leq \frac{\kappa}{100}.  
\end{align}
Using this fact together with Eq.~\eqref{eq:CurvatureSphere} and 
Eq.~\eqref{eq:First-Part-Hessian}, we conclude that under the current conditions, with probability $1 - o_n(1)$ 
Eq.~\eqref{eq:150} holds.

We then establish the smoothness property.
Recall that we have just proved that for a large enough $k$ and a small enough $R$, there exists $R_y > 0$, such that with probability $1 - o_n(1)$, the following statement is true: $\|\bD_{\bw} \bphi_{\bm}(\bw)\|_{\op} \leq 2$, \cref{eq:CurvatureSphere,eq:hessian-U-beta-t,eq:kappa-100} hold.
Putting these results together complete the proof for the smoothness property.

Finally, we prove the claim associated with $\bw_{\ast}(\by)$.
For any $\eps > 0$, by choosing a large enough $k$ and small enough $R_y$, we have $$\|\nabla_{\bw} \cuF_{\sTAP} (\bphi_{\bm}(\bw);  \bX, \by, \beta, t) |_{\bw = 0}\|_2^2 / n \leq \eps + o_P(1)$$ 
for all $\|\by - \by(t)\|_2 \leq \sqrt{n} R_y$. 
Therefore, by \cref{lemma:strong-convexity-property}, with probability $1 - o_n(1)$ we have $\| \bw_{\ast}(\by) \|_2 / \sqrt{n} \leq 20 \sqrt{\eps} / \kappa $ ($\kappa$ appears in \cref{eq:First-Part-Hessian}). The proof is complete via choosing a sufficiently small $\eps$.

\subsection{Proof of \cref{lemma:A1-4}}
\label{sec:proof-lemma:A1-4}

By \cref{lemma:AMP-optimal} we can choose a large enough $k$ such that with probability $1 - o_n(1)$ 
\begin{align*}
	\frac{1}{\sqrt{n}} \| \bphi_{\hat \bm_t^k}(\mathbf{0}_{n - 1}) - \bm_{\sB} \|_2 \leq \eps / 2. 
\end{align*} 
With a sufficiently large $k$ and sufficiently small $R$, it holds that with probability $1 - o_n(1)$ we have $\|\bD_{\bw}\bphi_{\hat \bm_t^k}(\bw)\|_{\op} \leq 2$ for all $\bw \in \Ball(\mathbf{0}_{n - 1}, \sqrt{n} R)$. 
By \cref{lemma:tangent-space-convextiy}, with probability $1 - o_n(1)$  
we have $\|\bw^{\ast}(\by(t))\|_2 < R \sqrt{n} / 4$. 
Choosing a small enough $R$, we see that with probability $1 - o_n(1)$ 
\begin{align*}
	\frac{1}{\sqrt{n}} \| \bphi_{\hat \bm_t^k}(\mathbf{0}_{n - 1}) - \bphi_{\hat \bm_t^k}(\bw^{\ast}(\by(t))) \|_2 \leq \eps / 2. 
\end{align*}
By triangle inequality, with probability $1 - o_n(1)$ it holds that $\| \bm_{\sB} - \bphi_{\hat \bm_t^k}(\bw^{\ast}(\by(t)))\|_2 / \sqrt{n} \leq \eps$.

\section{Proofs for \cref{alg:LR-high-SNR}}

\subsection{Proof of \cref{lemma:PosteriorAMP_LR}}
\label{sec:proof-lemma:PosteriorAMP_LR}

To prove \cref{lemma:PosteriorAMP_LR},
it suffices to verify that the fixed point equation $E = \mmse(\delta(\sigma^2 + E)^{-1})$ has a 
unique positive solution for all $\delta > \delta_{\salg}$ and noise level no larger than $\sigma_0$. 
If the solution is unique, then Bayes optimality 
(in the sense of \cref{eq:lemma5.1-eq}) follows by a standard argument,
see e.g. \cite{deshpande2014information} or the discussions that follows
 \cite[Theorem 2.5]{celentano2023mean}.

The next lemma establishes that the fixed point is unique for all 
$\delta$ large enough and all $\delta$ small enough. 
(The latter will be used to prove sampling guarantees at small $\delta$, cf. Appendix
\ref{sec:low-SNR}.)
We refer to 
Appendix \ref{sec:proof-lemma:unique-stationary} for its proof.
\begin{lemma}
\label{lemma:unique-stationary}
	We assume the conditions of \cref{lemma:PosteriorAMP_LR}. 
	Then the following statements hold: 
	\begin{enumerate}
		%
	\item There exists a constant $c_{0} \in (0, \infty)$ that depends only on $\pi_{\Theta}$, such that for all $\frac{\delta^{1/2}}{\sigma^3} < c_{0}$, the  fixed point equation below for all $t \geq 0$ has a unique positive solution.
		\begin{align}
		\label{eq:fixed-point-equation}
			E = \mmse\left(\delta (\sigma^2 + E)^{-1} + t\right). 
		\end{align}
	\item  There exists a constant $C_{0} \in (0, \infty)$ that depends only on $\pi_{\Theta}$, such that for all $\delta > C_{0}$, the fixed point equation \eqref{eq:fixed-point-equation} has a unique positive solution for all $t \geq 0$. 
	\item Assume the conditions of claim 1 or claim 2, and denote the unique solution to fixed point equation \eqref{eq:fixed-point-equation} by $E_{\ast} (\delta, \sigma, t)$. Then for all $k \in \NN$,
	\begin{align*}
		|E_k - E_{\ast}(\delta, \sigma, t)| \leq \frac{|1 - E_{\ast}(\delta, \sigma, t)|}{2^{k + 1}}.  
	\end{align*} 
	\end{enumerate} 
\end{lemma}

\subsection{Proof of \cref{thm:LR-main}}
\label{sec:proof-thm:LR-main}

As discussed in Section \ref{sec:general-setting},
it is sufficient to prove that the posterior mean estimator 
$\hbm_{\btheta}(\bX,\by_0,\sigma^2)$ is Lipschitz continuous for all 
$\delta\ge \delta_0(\pi_{\Theta},\sigma_0^2)$ and all $\sigma^2 \le \sigma^2_0$. 
We will prove  this fact  by using the local convexity of TAP free energy proven in
\cite{celentano2023mean}. 

We next replicate Theorem 3.1 and 3.4 of \cite{celentano2023mean} below for readers' convenience,
 which states that under the assumptions of \cref{thm:LR-main}, $\cuF_{\sTAP}$ has a unique local minimizer
  around the posterior first and second moment vectors. 
as before, we use the shorthands $\bm_{\sB}$ and $\bs_{\sB}$ 
for the posterior expectation and second moment:
\begin{align*}
	\bm_{\sB} = \E[\btheta \mid \bX, \by_0], \qquad \bs_{\sB} = \E[\btheta^2 \mid \bX, \by_0]. 
\end{align*}
\begin{theorem}[Theorem 3.1 and 3.4 of \cite{celentano2023mean}]
\label{lemma:LR-local-convexity}
	Let Assumptions \ref{assumption:LR}, \ref{assumption:three-points} hold and
	 further assume the conditions of Theorem \ref{thm:AMP-LinReg}.
	 Then there exist $\eps, \kappa > 0$ such that with probability approaching 1 as $n, p \to \infty$, $\cuF_{\sTAP}(\bm, \bs; \by_0, \bX, \sigma^2)$ has a unique local minimizer $(\bm_{\ast}, \bs_{\ast})$ in $\mathsf{B}((\bm_{\sB}, \bs_{\sB}), \eps \sqrt{p}) = \{(\bm, \bs) \in \Gamma^p: \|(\bm, \bs) - (\bm_{\sB}, \bs_{\sB})\|_2 < \eps \sqrt{p}\}$, and 
	\begin{align*}
		\inf_{(\bm, \bs) \in \mathsf{B}((\bm_{\sB}, \bs_{\sB}), \eps \sqrt{p})} \lambda_{\min} \left( \nabla^2 \cuF_{\sTAP}(\bm, \bs; \by_0, \bX, \sigma^2) \right) \geq \kappa. 
	\end{align*}
	Furthermore, 
	\begin{align*}
		\frac{1}{p} \left[ \|\bm_{\ast} - \bm_{\sB}\|_2^2 + \|\bs_{\ast} - \bs_{\sB}\|_2^2 \right] \toP 0. 
	\end{align*}
\end{theorem}

\begin{remark}
	If we assume either the conditions of the first point in \cref{lemma:unique-stationary},
	 or the conditions of the second point in \cref{lemma:unique-stationary}, then
	 the conditions of  Theorem \ref{thm:AMP-LinReg} follow. 
\end{remark}

Given local convexity, \cite{celentano2023mean} then design a two-stage algorithm 
that finds a local minimizer of the TAP free energy and at the same time, 
outputs consistent estimates of the posterior expectation. The algorithm implements
 Bayes AMP in the first stage, and runs natural gradient descent (NGD) in the second.
Apart from $(\bX, \by_0)$, several additional inputs are required to implement this algorithm: the number of AMP steps $K_{\AMP}$, the number of NGD steps $K_{\mathsf{NGD}}$, NGD step size $\eta$, and the noise level $\sigma^2$.  
For the compactness of presentation, we denote this algorithm by 
$\cA(\bX, \by_0,  K_{\AMP}, K_{\mathsf{NGD}}, \eta, \sigma^2)$, where we highlight 
the dependence on the number of AMP steps $K_{\AMP}$, the number of NGD steps $K_{\mathsf{NGD}}$,
and the NGD step size $\eta$.

We state below the theoretical guarantee for this procedure.
%
\begin{lemma}
\label{thm:LR-TAP}
Denote the output of $\cA(\bX, \by_0, K_{\AMP}, K_{\mathsf{NGD}}, \eta, \sigma^2)$ by 
$\hbm_{\btheta}(\bX, \by_0, \sigma^2)$. Further, recall that $\bm_{\sB}=\bm(\bX, \by_0, \sigma^2) = 
\EE[\btheta \mid \bX, \by] \in \RR^p$. 
We assume  $\delta > C_{0}$, where $C_{0}$ is given in the second point of
 \cref{lemma:unique-stationary}.
	\begin{enumerate}
			\item 
			 Then, fixing $(\delta, \pi_{\Theta})$, for all $\sigma^2 > 0$
			  the following statement is true: For any $\eps$, there exists 
			  $(K_{\AMP}, K_{\mathsf{NGD}}, \eta)$, such that with probability 
			  $1 - o_n(1)$, it holds that 
		\begin{align*}
			\frac{1}{\sqrt{p}} \|\hat \bm_{\btheta}(\bX, \by_0, \sigma^2) - \bm_{\btheta}(\bX, \by_0, \sigma^2) \|_2 \leq \eps. 
		\end{align*}
\item 
 There exists $r > 0$, such that for any $\by_1, \by_2 \in 
 \Ball^p(\by_0, r \sqrt{p})$, it holds that 
		\begin{align*}
			 \|\hbm_{\btheta} (\bX, \by_1, \sigma^2) - \hbm (\bX, \by_2, \sigma^2) \|_2 \leq C_2  \|\by_1 - \by_2\|_2 + \eps \sqrt{p}. 
		\end{align*}
		for a positive constant $C_2$ that is independent of $\sigma^2$. 
	\end{enumerate}
\end{lemma}
We defer the proof of \cref{thm:LR-TAP} to Appendix \ref{proof:thm:LR-TAP}.

Invoking the first and the second points of \cref{thm:LR-TAP}, we are able to verify
 \hyperlink{A1}{$\mathsf{(A1)}$} and \hyperlink{A3}{$\mathsf{(A3)}$}  of \cref{thm:gen}. 
Verification of Assumption \hyperlink{A2}{$\mathsf{(A2)}$} under the
 linear model is similar to that of the spiked model, and we skip it here for the sake of simplicity. 
The proof then follows as $\eps$ is arbitrarily small.

\subsection{Proof of \cref{lemma:unique-stationary}}
\label{sec:proof-lemma:unique-stationary}

We let $G(E) := \mmse({\delta (\sigma^2 + E)^{-1} } + t)$. Observe that $G$ is non-decreasing on
 $[0, 1]$, with $G(1) \leq 1$ and $G(0) \geq 0$. Hence to prove
  $E = G(E)$ has a unique fixed point in $[0,1]$, it suffices to prove
   $\sup_{E \in [0, 1]}|G'(E)| \leq 1 / 2$. 


%
%
%
%
%
%

\paragraph{Proof of point 1.}
By a direct calculation (see also the proof of Proposition \ref{prop:1}) we see that
\begin{align}
\label{eq:mmse-p}
	\mmse'(\gamma) = \gamma^{-1/2} \E\left[ (\Theta - \E[\Theta  \mid \gamma \Theta  + \sqrt{\gamma} G]) \cdot \Var[\Theta  \mid \gamma \Theta  + \sqrt{\gamma} G] \cdot G \right]
\end{align}
We denote by $m_8$ the eighth moment of $\pi_{\Theta}$. 
\cref{eq:mmse-p} also leads to the following upper bound: 
\begin{align*}
	|G'(E)| \leq \frac{\delta^{1/2}m_8^{3 / 8}}{ \sigma^3}, 
\end{align*}
which holds for all $E \in [0, 1]$ and $t > 0$. As a result, for a sufficiently small $c_{0}$ we have $\sup_{E \in [0, 1]} |G'(E)| \leq 1 / 2$.

\paragraph{Proof of point 2.}
Observe that $\Theta  - \E[\Theta  \mid \gamma \Theta  + \sqrt{\gamma} G] = - \gamma^{-1/2} G + \gamma^{-1/2} \E[G \mid \gamma \Theta  + \sqrt{\gamma} G]$. Therefore, 
\begin{align*}
	|G'(E)| \leq & \frac{1}{ \delta } \cdot \left| \E\left[ (-G + \E[G \mid \gamma \Theta  + \sqrt{\gamma} G]) \cdot \E[(-G + \E[G \mid \gamma \Theta  + \sqrt{\gamma} G])^2 \mid \gamma \Theta  + \sqrt{\gamma} G] G \right] \right| \\
	\leq & \frac{C_0}{\delta},
\end{align*}
where $C_0$ is a positive numerical constant. Therefore, for a large enough $\delta$, we have $\sup_{E \in [0, 1]}|G'(E)| \leq 1 / 2$. 

\paragraph{Proof of point 3.} 
Under the conditions at previous points $E\mapsto G(E)$ is a contraction, whence the claim follows.

\subsection{Proof of \cref{thm:LR-TAP}}
\label{proof:thm:LR-TAP}

\subsubsection*{Proof of the second point}

%
%

The consistency proof is a straightforward consequence of \cref{thm:AMP-LinReg}.  
We then proceed to the Lipschitz continuity proof. 
To this end, we define a modified TAP free energy: 
\begin{align*}
	\bar \cuF_{\sTAP}(\bm, \bs; \by_0, \bX, \sigma^2) = \cuF_{\sTAP}(\bm, \bs; \by_0, \bX, \sigma^2) - \frac{1}{2\sigma^2} \|\by_0\|_2^2. 
\end{align*}
Modifying the TAP free energy in this way does not change the local and global minimizers. 


Invoking \cref{lemma:LR-local-convexity}, we obtain that with high probability, $\bar \cuF_{\sTAP}$ is locally strongly convex in $\mathsf{B}((\bm_{\sB}, \bs_{\sB}), r_1 \sqrt{p})$ for some $r_1 > 0$. 
We shall choose $K_{\AMP}$ large enough, such that with probability $1 - o_n(1)$, $\|(\hat \bm_{\AMP}^{K_{\AMP}}(\bX, \by_0, \sigma^2), \hat \bs_{\AMP}^{K_{\AMP}}(\bX, \by_0, \sigma^2)) - (\bm_{\sB}, \bs_{\sB})\|_2 \leq r_1 \sqrt{p} / 10$. 
Here, we use $\hat \bm_{\AMP}^{K_{\AMP}}(\bX, \by_0, \sigma^2)$ and $\hat \bs_{\AMP}^{K_{\AMP}}(\bX, \by_0, \sigma^2)$ to represent the estimates of the posterior first and second moments, obtained via Bayes AMP with $K_{\AMP}$ iterations and inputs $(\by, \bX)$ (recall this is given in Section \ref{sec:general-setting}). 
We then choose $r$ small enough, such that with probability $1 - o_n(1)$, for all $\by' \in \mathsf{B}(\by_0, r\sqrt{n})$, we have $\|(\hat \bm_{\AMP}^{K_{\AMP}}(\bX, \by', \sigma^2), \hat \bs_{\AMP}^{K_{\AMP}}(\bX, \by', \sigma^2)) - (\bm_{\sB}, \bs_{\sB})\|_2 \leq r_1 \sqrt{p} / 5$. This is attainable as the non-linearities for Bayes AMP are all Lipschitz continuous due to the bounded support assumption, and the operator norm of the design matrix is with high probability upper bounded.
We then run NGD initialized at the AMP final iterates.

Observe that changing the response vector $\by_0$ does not alter the Hessian matrix of the TAP free energy. We then draw the following conclusion: With probability $1 - o_n(1)$, for all $\by' \in \mathsf{B}(\by_0, r)$, the TAP free energy $\bar \cuF_{\sTAP}(\bm, \bs; \by', \bX, \sigma^2)$ is strongly convex in $\mathsf{B}((\bm_{\sB}, \bs_{\sB}), r_1 \sqrt{p})$ and $\|(\hat \bm_{\AMP}^{K_{\AMP}}(\bX, \by', \sigma^2), \hat \bs_{\AMP}^{K_{\AMP}}(\bX, \by', \sigma^2)) - (\bm_{\sB}, \bs_{\sB})\|_2 \leq r_1 \sqrt{p} / 5$. We denote the local minimizer of the modified TAP free energy in this region by $(\bm_{\ast} (\by'), \bs_{\ast}(\by'))$ to emphasize the dependency on $\by'$.

We then upper bound $\|\bm_{\ast} (\by_0) - \bm_{\ast}(\by')\|_2$ for $\by' \in \mathsf{B}(\by_0, \sqrt{n} r)$. Note that by \cref{lemma:LR-local-convexity}, $p^{-1/2}\|(\bm_{\ast}(\by_0), \bs_{\ast}(\by_0)) - (\bm_{\sB}, \bs_{\sB})\|_2 = o_P(1)$, hence $(\bm_{\ast}(\by_0), \bs_{\ast}(\by_0))$ is inside the interior of the set $\mathsf{B}((\bm_{\sB}, \bs_{\sB}), \sqrt{n} r)$. Therefore, 
\begin{align*}
	& \nabla \bar\cuF_{\sTAP}(\bm_{\ast}(\by_0), \bs_{\ast}(\by_0); \by_0, \bX, \sigma^2) = 0,  \\
	& \nabla \bar\cuF_{\sTAP}(\bm_{\ast}(\by_0), \bs_{\ast}(\by_0); \by', \bX, \sigma^2) = - \frac{\bX^{\top} \by'}{\sigma^2} + \frac{\bX^{\top} \by_0}{\sigma^2}. 
\end{align*}
The second equation above is implied by the first. 
We define 
\begin{align*}
	\Gamma_{\leq }^p = \{(\bm, \bs) \in \Gamma^p \cap \mathsf{B}((\bm_{\sB}, \bs_{\sB}), r_1 \sqrt{p}): \bar \cuF_{\sTAP}(\bm, \bs; \by', \bX, \sigma^2) \leq \bar \cuF_{\sTAP} (\bm_{\ast}(\by_0), \bs_{\ast}(\by_0); \by', \bX, \sigma^2) \}. 
\end{align*}
Then, for any $(\bm, \bs) \in \Gamma_{\leq }^p$, we have 
\begin{align*}
	0 \geq & \bar \cuF_{\sTAP}(\bm, \bs; \by', \bX, \sigma^2) - \bar \cuF_{\sTAP} (\bm_{\ast}(\by_0), \bs_{\ast}(\by_0); \by', \bX, \sigma^2) \\
	\geq & - \sigma^{-2} \|\bX\|_{\op} \|\by_0 - \by'\|_2 \|(\bm, \bs) - (\bm_{\ast}(\by_0), \bs_{\ast}(\by_0)) \|_2 + \frac{\kappa}{2} \|(\bm, \bs) - (\bm_{\ast}(\by_0), \bs_{\ast}(\by_0))\|_2^2, 
\end{align*}
which further implies that 
\begin{align}
\label{eq:1488}
	\|(\bm_{\ast}(\by'), \bs_{\ast}(\by')) - (\bm_{\ast}(\by_0), \bs_{\ast}(\by_0))\|_2 \leq \frac{2}{\kappa \sigma^2} \|\bX\|_{\op} \|\by_0 - \by'\|_2. 
\end{align}
Standard application of random matrix theory implies that there exists a positive constant $C_{\delta}$ depending only on $\delta$, such that  with probability $1 - o_n(1)$ we have $\|\bX\|_{\op} \leq C_{\delta}$. Using this and \cref{eq:1488} above, we deduce that we can choose a small enough $r$, such that with probability $1 - o_n(1)$, for all $\by' \in \mathsf{B}(\by_0, r)$, 
\begin{align*}
	\|(\bm_{\ast}(\by_0), \bs_{\ast}(\by_0)) - (\bm_{\ast}(\by'), \bs_{\ast}(\by'))\|_2 \leq \frac{2 C_{\delta}}{\kappa \sigma^2}  \|\by_0 - \by'\|_2 \leq \sqrt{p} r_1  / 10,  
\end{align*}
which tells us that $(\bm_{\ast}(\by_0), \bs_{\ast}(\by_0))$ is within the interior of $\mathsf{B}((\bm_{\sB}, \bs_{\sB}), r_1 \sqrt{p})$. 
Repeating the above procedure with $(\by_0, \by')$ replaced by $(\by_1, \by_2)$ for $\by_1, \by_2 \in \mathsf{B}(\by_0, r\sqrt{n})$, we conclude that with probability $1 - o_n(1)$, for all $\by_1, \by_2 \in \mathsf{B}(\by_0, r\sqrt{n})$
\begin{align*}
	\|(\bm_{\ast}(\by_1), \bs_{\ast}(\by_1)) - (\bm_{\ast}(\by_2), \bs_{\ast}(\by_2))\|_2 \leq \frac{2 C_{\delta}}{\kappa \sigma^2}  \|\by_1 - \by_2\|_2. 
\end{align*}
The proof is then complete, as the TAP free energy is strongly convex, and we can choose an appropriate $(\eta, K_{\mathsf{NGD}})$ and run NGD initialized at the final AMP iterate to consistently approximate $(\bm_{\ast}, \bs_{\ast})$. (Follows from \cite[Lemma H.1]{celentano2023mean})

\section{High-dimensional regression at low SNR}
\label{sec:low-SNR}

Alternatively, to sample from the posterior distribution of the linear coefficients, we may track the stochastic localization process $\dd \bz(t) = \bm(\bz(t), t) \dd t + \dd \bB(t)$ that is initialized at the origin, where $\bB(t)$ is a standard Brownian motion in $\RR^p$. For the linear regression setting that we consider, $\bm(\bz(t), t) = \E[\btheta \mid \bX, \by_0, \bz(t)]$. Here, $\bz(t) = t \btheta + \bG(t) \in \RR^p$, and $(\bG(t))_{t \geq 0}, (\bB(t))_{t \geq 0}$ are $p$-dimensional standard Brownian motions. Approximating $\bz(T)$ for a sufficiently large $T$ then gives approximate samples from the target posterior distribution.

We define the posterior first and second moment vectors as 
\begin{align*}
	\bm_{\sB, t} = \E[\btheta \mid \bX, \by_0, \bz(t)], \qquad \bs_{\sB, t} = \E[\btheta^2 \mid \bX, \by_0, \bz(t)], 
\end{align*}
both are in $\RR^p$. 
We consider the TAP free energy, defined over the domain $(\bm, \bs) \in \Gamma^p \subseteq \RR^p \times \RR^p$, 
\begin{align}
\label{eq:TAP-low-SNR}
	\cuF_{\sTAP}(\bm, \bs; t, \by_0, \bX, \bz(t), \sigma^2) = & \frac{n}{2} \log 2\pi \sigma^2 + D_0(\bm, \bs) + \frac{1}{2\sigma^2} \|\by_0 - \bX\bm\|_2^2 + \frac{n}{2} \log \left( 1 + \frac{S(\bs) - Q(\bm)}{\sigma^2} \right) \\
	& + \frac{1}{2t} \|\bz(t) - t \bm\|_2^2 + \frac{tp}{2}(S(\bs) - Q(\bm)). \nonumber
\end{align}
We ignore the dependency on $(\by_0, \bX, \bz(t), \sigma^2)$ when there is no confusion.

In the sequel, we shall prove the existence of a local minimizer of the TAP free energy $\cuF_{\sTAP}(\bm, \bs; t)$ over $(\bm, \bs) \in \Gamma^p$, near the true posterior first and second moments.
To this end, we establish the local convexity of the TAP free energy.

\subsection{Convexity of the TAP free energy}


We show that the TAP free energy is globally strongly convex in the low SNR regime.
 We formally state this result as \cref{lemma:low-SNR}, which  is a direct consequence of \cite[Proposition 3.3]{celentano2023mean}. The proof of the lemma is deferred to Appendix \ref{proof:lemma:low-SNR}. 

\begin{lemma}
\label{lemma:low-SNR}
	We assume Assumption \ref{assumption:LR}. Then there exist constants $\kappa$, $c_{0} > 0$ that depend only on $\pi_{\Theta}$, such that if $\delta / \sigma^2 < c_{0}$, then for a sufficiently large $(n, p)$,  $\nabla^2 \cuF_{\sTAP} (\bm, \bs; t) \succeq \kappa I_{2p}$ for all $(\bm, \bs) \in \Gamma^p$.  
\end{lemma}

\subsection{Sampling algorithm and guarantees}

We are then ready to state our sampling algorithms and the associated theoretical guarantees. 
We first give the algorithm for the low SNR regimes as \cref{alg:LR-low-SNR}. 
Motivated by \cref{lemma:low-SNR}, we propose to run NGD without the AMP first stage. 
We denote the output of this procedure by $\hat\bm(\bz, t)$ as an estimate of $\bm(\bz, t)$.  

\begin{algorithm}[ht]
\caption{Diffusion-based sampling for linear models in a low SNR regime}\label{alg:LR-low-SNR}
\textbf{Input: }Data $(\bX, \by_0)$, parameters $( K_{\mathsf{NGD}}, \eta, \delta, \sigma^2, L, \Delta)$;	
\begin{algorithmic}[1]
\State Set $\hat{\bz}_0 = \bfzero_p$;
\For{$\ell = 0,1,\cdots, L - 1$}
	\State Draw $\bw_{\ell + 1} \sim \normal(0,\id_p)$ independent of everything so far;
	\State Denote $\hat{\bm}(\hat{\bz}_\ell,\ell\Delta)$ the output of NGD on \eqref{eq:TAP-low-SNR} with $K_{\mathsf{NGD}}$ iterations with step size $\eta$
		\State Update $\hat{\bz}_{\ell + 1} = \hat{\bz}_\ell + \hat{\bm}( \hat{\bz}_\ell, \ell\Delta) \Delta + \sqrt{\Delta} \bw_{\ell + 1}$;	
	\EndFor
\State \Return $\btheta^{\salg}  = \hat{\bm}(\hat{\bz}_L, L\Delta)$;
\end{algorithmic}
\end{algorithm}

We next show that Bayes AMP consistently estimates the posterior mean for a sufficiently low SNR. We present this result as the lemma below: 
\begin{lemma}
\label{lemma:F2}
	We assume $\delta$ is sufficiently small such that the conditions of \cref{lemma:low-SNR} and the first point of \cref{lemma:unique-stationary} are satisfied. We also assume Assumption \ref{assumption:LR}. Then, for any $\eps, t> 0$, there exist $( K_{\mathsf{NGD}}, \eta)$, such that with probability $1 - o_n(1)$, it holds that 
		\begin{align*}
			\frac{1}{\sqrt{p}} \| \hat \bm(\bz(t), t) - \bm(\bz(t), t) \|_2 \leq \eps. 
		\end{align*}
		In addition, there exists a constant $C_1 > 0$, such that for any $\bz_1, \bz_2 \in \RR^p$,  
		\begin{align*}
			\|\hat\bm(\bz_1, t) - \hat\bm(\bz_2, t)\|_2 \leq C_1 \|\bz_1 - \bz_2\|_2 + \eps \sqrt{p}. 
		\end{align*}
\end{lemma}

We prove \cref{lemma:F2} in Appendix \ref{sec:proof-lemma:F2}. 
We present guarantee for \cref{alg:LR-low-SNR} in \cref{thm:LR-main-low-SNR}. 

\begin{theorem}
\label{thm:LR-main-low-SNR}
	 We denote by $\mu_{\bX, \by_0}^{\salg}$ the law of $\btheta^{\salg}$, obtained using \cref{alg:LR-low-SNR}, and denote by $\mu_{\bX, \by_0}$ the true posterior distribution. We also assume Assumptions \ref{assumption:LR} and \ref{assumption:three-points}. Then the following statements are true: 
		%
		We assume $\delta$ is small enough as in the second point of \cref{lemma:unique-stationary}. 
		Then for any $\xi > 0$, there exist $(K_{\AMP}, K_{\mathsf{NGD}}, \eta, L, \Delta)$ that depend only on $(\delta, \sigma^2, \pi_{\Theta})$, such that if \cref{alg:LR-high-SNR} takes as input $(\bX, \by_0, K_{\AMP}, K_{\mathsf{NGD}}, \eta, \delta, \sigma^2, L, \Delta)$, then with probability $1 - o_n(1)$ over the randomness of $(\bX, \by_0)$, 
		\begin{align*}
			W_{2, p}(\mu_{\bX, \by_0}, \mu_{\bX, \by_0}^{\salg}) \leq \xi. 
		\end{align*}
\end{theorem} 
\begin{proof}[Proof of \cref{thm:LR-main-low-SNR}]
	The proof is similar to that for \cref{thm:main2,thm:spiked-continuous,thm:LR-main}, we skip it to avoid redundancy. 
\end{proof}

\subsection{Proof of \cref{lemma:low-SNR}}
\label{proof:lemma:low-SNR}

According to Proposition 3.3 of \cite{celentano2023mean}, the lemma holds for $t = 0$. The proof is complete by observing that $\nabla^2 \cuF_{\sTAP} (\bm, \bs; t) = \nabla^2 \cuF_{\sTAP} (\bm, \bs; 0)$ for all $(\bm, \bs) \in \Gamma^p$ and $t > 0$.

\subsection{Proof of \cref{lemma:F2}}
\label{sec:proof-lemma:F2}

In this proof, we find the following lemma useful: 
\begin{lemma}
	\label{lemma:strong-convex-functions}
	If $f: \RR^n \to \RR$ is smooth and $\mu$-strongly convex for $\mu > 0$, then for $x_{\ast}$ that minimizes $f$, we have 
	\begin{align*}
		& \frac{1}{2\mu} \|\nabla f(x)\|_2^2 \geq f(x) - f(x_{\ast}) \geq \frac{\mu}{2} \|x - x_{\ast}\|_2^2, \\
		& f(x) - f(x_{\ast}) \leq \|\nabla f(x)\|_2 \cdot \|x - x_{\ast}\|_2.  
	\end{align*}
	For all $x, y \in \RR^n$, we have 
	\begin{align*}
		f(y) \geq f(x) + \langle \nabla f(x), y - x \rangle + \frac{\mu}{2}\|y - x\|_2^2. 
	\end{align*}
\end{lemma}

\subsubsection*{Bayes AMP with side information}

To prove the lemma, we first state the Bayes AMP algorithm with side information.  
Given inputs $(\bX, \by_0, \bz(t))$,  the AMP algorithm with side information takes the following form:
\begin{align}
\label{eq:linear-AMP-side}
	\left\{
	\begin{array}{l}
		\bb^{k + 1} = \bX^{\top} f_k(\ba^k, \by_0) - \xi_k g_k(\bb^k, \bz(t)), \\
		\ba^k = \bX \bg_k(\bb^k, \bz(t)) - \eta_k f_{k - 1}(\ba^{k - 1}, \by_0), 
	\end{array} \right.
\end{align}
where at initialization we set $\ba^{-1} = \mathbf{0}_n$. After $k$ iterations, algorithm \eqref{eq:linear-AMP-side} outputs $g_k(\bb^k, \bz(t)) \in \RR^p$ as an estimate of $\btheta$.  
In the above iteration, we assume $f_k, g_k: \RR^2 \to \RR$ are Lipschitz continuous and apply on matrices row-wisely (note that $f_{-1}$ maps from $\RR$ to $\RR$). The scalar sequences $(\xi_k)_{k \geq 0}$ and $(\eta_k)_{k \geq 0}$ are determined by the following state evolution recursion: 
\begin{align*}
	& \xi_k = \delta \E[\partial_1 f_k(\bar{G}_k, \bar{G}_{\ast} + \eps )], \qquad \eta_k =  \E[\partial_1 g_k(\mu_t \Theta + G_k, U)], \qquad \sigma_k^2 = \delta \E[f_{k - 1}(\bar{G}_{k - 1}, \bar{G}_{\ast} + \eps)^2] \\
	&  \bar\sigma_k^2 =  \E[g_t(\mu_k \Theta + G_k, U)^2], \qquad \bar\sigma_{k, \ast} =  \E[\Theta g_k(\mu_k \Theta + G_k, U)], \qquad \mu_k = \delta \E[\partial_{2} f_{k - 1}(\bar{G}_{k - 1}, \bar{G}_{\ast} + \eps) ], 
\end{align*} 
where $\Theta \sim \pi_{\Theta}$, $\bar{G}_{\ast} \sim \normal(0,1)$, $\eps \sim \normal(0, \sigma^2)$, $G_k \sim \normal(0, \sigma_t^2)$, $\bar{G}_k \sim \normal(0, \bar\sigma_k^2)$, and $U \sim t \Theta + \sqrt{t}G$. Here, $\Theta \perp G_k$, $\eps \perp (\bar{G}_{\ast}, \bar{G}_k)$, and $(\bar{G}_{\ast}, \bar{G}_k)$ are jointly Gaussian.
We also assume $G \sim \normal(0, 1)$ and is independent of anything else.  

Choosing appropriate non-linearities, we get the Bayes AMP: 
\begin{align*}
	& f_k(\bar{G}_k, \bar{G}_{\ast} + \eps) = \E[\eps \mid \bar{G}_k,  \bar{G}_{\ast} + \eps], \qquad g_k(\mu_k \Theta + G_k, U) = \E[\Theta \mid \mu_k \Theta + G_k, U]. 
\end{align*}
In this case, straightforward computation gives  $\bar\sigma_k^2 = \bar\sigma_{k, \ast}$ and $\sigma^2 \mu_k = \sigma_k^2$. In addition, the state evolution quantities admit the following expressions:  
\begin{align*}
	& f_k(x, y) = \frac{\sigma^2(y - x)}{\sigma^2 + 1 - \bar\sigma_k^2},  \\ 
	& g_k (x, u) = \E[\Theta \mid (\sigma^{-2}\mu_k \Theta + t) + \sqrt{\sigma^{-2}\mu_k + t} G = \sigma^{-2}x + u], \\
	& \mu_k = \frac{\delta\sigma^2}{\sigma^2 + 1 - \bar\sigma_{k - 1}^2}, \qquad \bar\sigma_k^2 =  \E\left[ \E[\Theta \mid (\sigma^{-2}\mu_k + t)  \Theta + \sqrt{\sigma^{-2}\mu_k + t} G ]^2 \right], \\
	& \xi_k = -\frac{\delta\sigma^2}{\sigma^2 + 1 - \bar\sigma_k^2}, \qquad \eta_k = \sigma^{-2}\E\left[ \Var[\Theta \mid (\sigma^{-2}\mu_k + t)  \Theta + \sqrt{\sigma^{-2}\mu_k + t} G] \right]. 
\end{align*}
We define $E_{k} := 1 - \bar\sigma_{k, \ast}$.
As initialization, we get $E_{-1} = 1$.
For $k \geq 0$ we have
\begin{align*}
	& E_{k} = \mmse \left( {\sigma^{-2} \mu_k} + t\right) = \mmse \left( { \delta (\sigma^2 + E_{k - 1})^{-1} } + t \right), \\
	& \mmse(\gamma) := \E\left[ \left( \Theta - \E[\Theta \mid \gamma \Theta + \sqrt{\gamma} G] \right)^2 \right]. 
\end{align*}

\subsubsection*{Proof of the lemma}


By assumption, fixed point equation \eqref{eq:fixed-point-equation} has a unique positive solution. 
Therefore, Bayes AMP consistently approximates the posterior expectation. 
More precisely, for any $\eps > 0$ and large enough $k \in \NN_+$, with probability $1 - o_n(1)$ we have  
\begin{align}
\label{eq:145}
	\frac{1}{\sqrt{p}} \|\bm(\bz(t), t) - g_{k}(\bb^{k}, \bz(t))  \|_2 \leq \eps,
\end{align}
where we recall $\bb^k$ is the $k$-th AMP iterate.  

Define $s_k (x, u) = \E[\Theta^2 \mid (\sigma^{-2}\mu_k \Theta + t) + \sqrt{\sigma^{-2}\mu_k + t} G = \sigma^{-2}x + u]$. 
Taking the gradient of the TAP free energy, we get 
\begin{align*}
	\nabla \cuF_{\sTAP}(\bm, \bs; t, \by_0, \bX, \bz(t), \sigma^2) = \left( \begin{array}{c}
			\blambda - \frac{1}{\sigma^2} \bX^{\top}(\by_0 - \bX \bm) - \frac{n / d}{\sigma^2 + S(\bs) - Q(\bm)} \bm - \bz(t) \\
			-\frac{1}{2} \bgamma + \frac{n / d}{2(\sigma^2 + S(\bs) - Q(\bm))} \mathbf{1} + \frac{t}{2} \mathbf{1}
		\end{array} \right). 
\end{align*}
Plugging the AMP iterates into the TAP free energy gradient formula, and using the state evolution of the AMP algorithm, we see that for any $\eps$, if we take a large enough $k$, then with probability $1 - o_n(1)$ we have 
\begin{align}
\label{eq:146}
	\frac{1}{n + p} \|\nabla \cuF_{\sTAP}(g_k(\bb^k, \bz(t)), s_k(\bb^k, \bz(t)); t, \by_0, \bX, \bz(t), \sigma^2)\|_2^2 \leq \eps. 
\end{align}
By \cref{lemma:low-SNR}, we see that for a large enough $(n, p)$, the TAP free energy is strongly convex with $\nabla^2 \cuF_{\sTAP} (\bm, \bs; t) \succeq \kappa I_{2p}$ for all $(\bm, \bs) \in \Gamma^p$. 
Further recall that $\Gamma$ is convex and open in $\RR^2$, hence $\Gamma^p$ is convex and open in $\RR^{2p}$. 
By strong convexity, we have 
\begin{align}
\label{eq:147}
	\|(\bm, \bs) - (\bm_{\ast}, \bs_{\ast})\|_2^2 \leq \frac{1}{\kappa^2} \|\nabla \cuF_{\sTAP}(\bm, \bs; t)\|_2^2. 
\end{align}
Putting together \cref{eq:145,eq:146,eq:147}, we conclude that for any $\eps > 0$, with probability $1 - o_n(1)$ we have $\| \bm_{\ast} - \bm(\bz(t), t) \|_2 \leq \eps \sqrt{p} / 2$. Once again employing the strong convexity of TAP free energy and applying \cite[Lemma H.1]{celentano2023mean}, we can choose an appropriate $(K_{\mathsf{NGD}}, \eta)$, such that $\|\hat\bm(\bz(t), t) - \bm_{\ast}\|_2 \leq \eps \sqrt{p} / 2$. By triangle inequality, 
\begin{align*}
	\|\hat\bm(\bz(t), t) - \bm(\bz(t), t)\|_2 \leq \eps \sqrt{p}. 
\end{align*}
Next, we write $\bm_{\ast}(\bz(t), t) = \bm_{\ast}$ to emphasize the dependency on $(\bz(t), t)$. 
We define 
\begin{align*}
	\bar \cuF_{\sTAP}(\bm, \bs; t, \by_0, \bX, \bz(t), \sigma^2) = & \frac{n}{2} \log 2\pi \sigma^2 + D_0(\bm, \bs) + \frac{1}{2\sigma^2} \|\by_0 - \bX\bm\|_2^2 + \frac{n}{2} \log \left( 1 + \frac{S(\bs) - Q(\bm)}{\sigma^2} \right) \\
	& - \langle \bm, \bz(t) \rangle + \frac{t}{2} \|\bm\|_2^2 + \frac{tp}{2}(S(\bs) - Q(\bm)). 
\end{align*}
Note that $(\bm_{\ast}(\bz, t), \bs_{\ast}(\bz, t))$ also minimizes $\bar \cuF_{\sTAP}(\bm, \bs; t, \by_0, \bX, \bz(t), \sigma^2)$ for all $\bz \in \RR^p$.  

For $\bz_1, \bz_2 \in \RR^p$, we then upper bound $\|\bm_{\ast}(\bz_1, t) - \bm_{\ast}(\bz_2, t)\|_2$. Straightforward computation gives 
\begin{align*}
	& \nabla \bar \cuF_{\sTAP} \big( \bm_{\ast}(\bz_1, t), \bs_{\ast}(\bz_1, t); t, \by_0, \bX, \bz_1, \sigma^2 \big) = 0, \\
	& \nabla \bar \cuF_{\sTAP} \big( \bm_{\ast}(\bz_1, t), \bs_{\ast}(\bz_1, t); t, \by_0, \bX, \bz_2, \sigma^2 \big) = - \bz_2 + \bz_1. 
\end{align*}
We define 
$$\Gamma_{\leq}^p = \{(\bm, \bs) \in \Gamma^p:\bar \cuF_{\sTAP} \big( \bm, \bs; t, \by_0, \bX, \bz_2, \sigma^2 \big) \leq \bar \cuF_{\sTAP} \big( \bm_{\ast}(\bz_1, t), \bs_{\ast}(\bz_1, t); t, \by_0, \bX, \bz_2, \sigma^2 \big)\}.$$
Then, for any $(\bm, \bs) \in \Gamma_{\leq}^p$, we have 
\begin{align*}
	0 \geq & \bar \cuF_{\sTAP} \big( \bm, \bs; t, \by_0, \bX, \bz_2, \sigma^2 \big) - \bar \cuF_{\sTAP} \big( \bm_{\ast}(\bz_1, t), \bs_{\ast}(\bz_1, t); t, \by_0, \bX, \bz_2, \sigma^2 \big) \\
	\geq & - \|\bz_1 - \bz_2\|_2 \cdot \|(\bm, \bs) - (\bm_{\ast}(\bz_1, t), \bs_{\ast}(\bz_1, t)) \|_2 + \frac{\kappa}{2} \|(\bm, \bs) - (\bm_{\ast}(\bz_1, t), \bs_{\ast}(\bz_1, t)) \|_2^2. 
\end{align*}
As a consequence,
\begin{align*}
	\|(\bm_{\ast}(\bz_1, t), \bs_{\ast}(\bz_1, t)) - (\bm_{\ast}(\bz_2, t), \bs_{\ast}(\bz_2, t))\|_2 \leq \frac{2}{\kappa} \|\bz_1 - \bz_2\|_2. 
\end{align*}
By \cref{lemma:strong-convex-functions}, 
\begin{align*}
	 & \bar \cuF_{\sTAP} \big( \bm_{\ast}(\bz_1, t), \bs_{\ast}(\bz_1, t); t, \by_0, \bX, \bz_2, \sigma^2 \big) - \bar \cuF_{\sTAP} \big( \bm_{\ast}(\bz_2, t), \bs_{\ast}(\bz_2, t); t, \by_0, \bX, \bz_2, \sigma^2 \big)  \\
	 \leq & \| \nabla \bar \cuF_{\sTAP} \big( \bm_{\ast}(\bz_1, t), \bs_{\ast}(\bz_1, t); t, \by_0, \bX, \bz_2, \sigma^2 \big) \|_2 \cdot \|(\bm_{\ast}(\bz_1, t), \bs_{\ast}(\bz_1, t)) - (\bm_{\ast}(\bz_2, t), \bs_{\ast}(\bz_2, t))\|_2 \\
	 \leq & \frac{2}{\kappa} \|\bz_1 - \bz_2\|_2^2. 
\end{align*}
Also by \cref{lemma:strong-convex-functions},
\begin{align*}
	& \bar \cuF_{\sTAP} \big( \bm_{\ast}(\bz_1, t), \bs_{\ast}(\bz_1, t); t, \by_0, \bX, \bz_2, \sigma^2 \big) - \bar \cuF_{\sTAP} \big( \bm_{\ast}(\bz_2, t), \bs_{\ast}(\bz_2, t); t, \by_0, \bX, \bz_2, \sigma^2 \big) \\
	\geq & \frac{\kappa}{2} \|(\bm_{\ast}(\bz_1, t), \bs_{\ast}(\bz_1, t)) - (\bm_{\ast}(\bz_2, t), \bs_{\ast}(\bz_2, t))\|_2^2.  
\end{align*}
Therefore, 
\begin{align*}
	\|(\bm_{\ast}(\bz_1, t), \bs_{\ast}(\bz_1, t)) - (\bm_{\ast}(\bz_2, t), \bs_{\ast}(\bz_2, t))\|_2 \leq \frac{2}{\kappa} \|\bz_1 - \bz_2\|_2.  
\end{align*}
The proof is complete.

\section{Proof of the supporting lemmas}
\label{sec:supporting-lemmas}

\subsection{Proof of \cref{lemma:AMP-optimal}}\label{sec:proof-of-lemma:AMP-optimal}

We will consider separately the case of non-symmetric and symmetric prior $\pi_{\Theta}$.
The latter case is more challenging and requires new ideas with respect to 
earlier work, e.g. \cite{alaoui2022sampling}.
We first prove the claim regarding the stationary point of $\Phi$. 

\subsubsection{Proof for the stationary point of $\Phi$}
\label{sec:proof-lemma:large-beta-cond1}

Recall that the state evolution recursion for rank-one estimation is given in
\cref{eq:AMPSE} and the corresponding free energy functional 
$\Phi$ is defined in \cref{eq:Phi}.
\begin{lemma}
\label{lemma:large-beta-cond1}
	Assume $\pi_{\Theta}$ has finite eighth moment and unit second moment. Then, there exists 
	$\beta_0(\pi_{\Theta})$ that depends only on $\pi_{\Theta}$, 
	such that for all $\beta \geq \beta_0(\pi_{\Theta})$, 
	there exists exactly one positive solution $\gamma_{\beta, t}\in (t,\infty)$ to 
	$\frac{\partial}{\partial \gamma} \Phi(\gamma, \beta, t) = 0$. Equivalently, there exist
	exactly one fixed point $\gamma_{\beta, t}>0$ for the state evolution recursion of \cref{eq:AMPSE}. This also corresponds to the global minimum of $\Phi$ on $(t, \infty)$. 
	In addition,  $\gamma_{\beta, t}\ge \beta^2-1$ and (recalling the definition in \cref{eq:Qdef})
	$q_{\beta,t} \ge 1-\beta^{-2}$.
\end{lemma}
\begin{proof}
Taking the derivative of $\gamma \mapsto \Phi(\gamma, \beta, t)$, we get 
\begin{align*}
	\frac{\partial}{\partial \gamma} \Phi(\gamma, \beta, t) = \frac{\gamma}{2\beta^2} - \frac{1}{2} + \frac{1}{2} \mmse(\gamma + t) = \frac{\gamma}{2\beta^2} - \frac{1}{2} \E\left[ \E[\Theta \mid (\gamma + t) \Theta + \sqrt{\gamma + t} G]^2 \right]. 
\end{align*}
Hence, we recover the well known fact that stationary points of
$\gamma\mapsto  \Phi(\gamma, \beta, t)$ coincide with fixed points of the 
state evolution recursion of \cref{eq:AMPSE}.
By comparing the optimal estimator to the optimal linear estimator, we get
 $\E\left[ \E[\Theta \mid (\gamma + t) \Theta + \sqrt{\gamma + t} G]^2 \right] \geq
  (\gamma + t)/(1+\gamma+t) \geq \gamma/(1+\gamma)$. 
  Therefore, for all sufficiently large $\beta$, in order for
   $\frac{\partial}{\partial \gamma} \Phi(\gamma, \beta, t) = 0$ to hold at some positive $\gamma$, 
   we must have $\gamma \geq \beta^2 -1\ge \beta^2/2$ (assuming $\beta^2\ge 2$). 
	
Standard computation implies that 
\begin{align*}
	\mmse'(\gamma) = \gamma^{-1/2} \E\left[ (\Theta - \E[\Theta \mid \gamma \Theta + \sqrt{\gamma} G]) \cdot \Var[\Theta \mid \gamma \Theta + \sqrt{\gamma} G] \cdot G \right],
\end{align*}
where $(\Theta, G) \sim \pi_{\Theta} \otimes \normal(0, 1)$.
Note that $\Theta - \E[\Theta \mid \gamma \Theta + \sqrt{\gamma} G] = - \gamma^{-1/2} G + \gamma^{-1/2} \E[G \mid \gamma \Theta + \sqrt{\gamma} G]$. 
Therefore, for all $\gamma \geq \beta^2 / 2$, it holds that $|\mmse'(\gamma + t)| \leq C_{0} \beta^{-3}$, where $C_0$ is a numerical constant. 
From these discussions, we see that for all $\gamma \geq \beta^2 / 2$, 
\begin{align*}
	\frac{\partial^2}{\partial \gamma^2} \Phi(\gamma, \beta, t) = \frac{1}{2\beta^2} + \frac{1}{2} \mmse'(\gamma + t) \geq \frac{1}{2\beta^2} - \frac{C_0}{2\beta^3},  
\end{align*}
which is positive for a sufficiently large $\beta$. 
We then deduce that there is at most one positive solution to 
$\frac{\partial}{\partial \gamma} \Phi(\gamma, \beta, t) = 0$. 
In addition, it is easy to check that $\frac{\partial}{\partial \gamma} \Phi(\gamma=0, \beta, t) < 0$,
which implies that the solution exists.
Finally, we note that $\frac{\partial}{\partial \gamma} \Phi(\gamma, \beta, t) > 0$ for a sufficiently large $\gamma$, implying the unique positive stationary point is also the global minimum. 
\end{proof}

\subsubsection{Proof of AMP consistency for non-symmetric $\pi_{\Theta}$}
\label{sec:ProofNonSymmetric}

In this section we assume $\pi_{\Theta}$ is non-symmetric together with the moment conditions.
Recall the state evolution sequence  $(\gamma_t^k)_{k\ge 0}$ is defined in
\cref{eq:AMPSE}, which we rewrite as $\gamma_t^{k+1} = \sM(\gamma_t^k;\beta,t)$.
By Lemma \ref{lemma:large-beta-cond1}, for $\beta>\beta_0(\pi_{\Theta})$ there is a unique 
$\gamma_{\beta, t} >0$ such that 
$\gamma_{\beta, t} = \sM(\gamma_{\beta, t};\beta,t)$. 
Observe that $\sM(0;\beta,t)> 0$ and $\gamma \mapsto \sM(\gamma;\beta,t)$ is strictly increasing, we then have
$\sM(\gamma;\beta,t)>\gamma$ for $\gamma<\gamma_{\beta, t}$.  
Therefore, $(\gamma_t^k)_{k\ge 0}$ is a increasing sequence, upper bounded by $\gamma_{\beta, t}$. In addition, 
and 
\begin{align}
\lim_{k\to\infty}\gamma_t^k=\gamma_{\beta, t}\, .
\end{align}

Our next proposition provides a more quantitative control of this convergence rate.
\begin{proposition}\label{prop:1}
Assume $\pi_{\Theta}$ has unit second moment and $\EE[ \Theta^8 ]^{1/4} \leq M$.
Then there exists $\beta_0=\beta_0(\pi_{\Theta}) > 0$ depending only on $\pi_{\Theta}$, 
such that for all 
	$\beta \geq \beta_0(\pi_{\Theta})$ and $t \geq 0$, $\gamma_{\beta, t}$ is the 
	unique positive solution to the fixed point equation \cref{eq:Phi}. 
	Furthermore, for all $k\ge 0$,
	\begin{align*}
		1 - 2^{-k} \leq \frac{\gamma_t^k}{\gamma_{\beta, t}} \leq 1. 
	\end{align*}
\end{proposition} 
\begin{proof}[Proof of Proposition \ref{prop:1}]
	We define the non-decreasing function 
	\begin{align}
	H(\gamma) &:= \beta^2(1 - \mathsf{mmse}(\gamma ))\\
	& = \beta^2\Big(1-\E\big[(\Theta-f_{\sB}(Y;\gamma))^2\big]\Big)\, ,
	\end{align}
	where $Y = \gamma \Theta+\sqrt{\gamma} G$ and 
	\begin{align*}
	f_{\sB}(y;\gamma) := \E[\Theta|\gamma \Theta+\sqrt{\gamma} G=y] = \sF\Big(\frac{y}{\sqrt{\gamma}};\gamma\Big)\,.
	\end{align*}
	We then have
	(here we repeatedly use the fact that $\E[(\Theta-f_{\sB}(Y;\gamma)) \, h(Y)]=0$ for
	 any function $h$ such that the expectation exists)
	\begin{align*}
	H'(\gamma) &= 2\beta^2\E\Big[\big(\Theta-f_{\sB}(Y;\gamma)\big) \partial_{\gamma}f_{\sB}(Y;\gamma)\Big]+
	 2\beta^2\E\Big[\big(\Theta-f_{\sB}(Y;\gamma)\Big) \partial_{Y}f_{\sB}(Y;\gamma) 
	 \big(\Theta + (G/2)\gamma^{-1/2}\big)\Big]\\
	 & = -\beta^2\gamma^{-1/2}\E\big[\big(\Theta-f_{\sB}(Y;\gamma)\big) \Var(\Theta|Y) G\big]\\
	 & \le \beta^2\gamma^{-1/2} \mmse(\gamma)^{1/2}
	  \E\big[\Var(\Theta|Y)^2G^2\big]^{1/2}\\
	  &\le 2\beta^2 \gamma^{-1/2}\mmse(\gamma)^{1/2}
	  \E\big[\Var(\Theta|Y)^4\big]^{1/4}\\
	  &\le 2 M^{3/2}\beta^2\gamma^{-1/2} \mmse(\gamma)^{3/4}
	 \, ,
	\end{align*} 
	where the last inequalities follow from Cauchy-Schwartz. 
	Recalling that $\mmse(\gamma)\le 1/\gamma$, we obtain
	\begin{align}
	0\le H'(\gamma+t)\le 2 M^{3/4}\frac{\beta^2}{(\gamma+t)^{5/4}}\, .
	\end{align}
	And therefore, for $\gamma\ge \gamma_t^0$, and all $\beta\ge \beta_0(\pi_{\Theta})$,
	we proved that
	$	0\le H'(\gamma+t)\le 1/2$
	
	This implies 
	\begin{align*}
		\left| \gamma_{\beta, t} - \gamma_t^{k + 1} \right| = \left| 
		H(t+\gamma_{\beta, t}) - H(t+\gamma_t^{k + 1}) \right| \leq \frac{1}{2} \left| \gamma_{\beta, t} - \gamma_t^{k + 1} \right|\, ,
	\end{align*}
	which concludes the proof of the proposition. 
	
	In addition, for $t_1, t_2 > 0$, it holds that 
	\begin{align*}
		|\gamma_{\beta, t_1} - \gamma_{\beta, t_2}| = & |H(t_1 + \gamma_{\beta, t_1}) - H(t_2 + \gamma_{\beta, t_2})| \\
		\leq & \frac{1}{2} |\gamma_{\beta, t_1} - \gamma_{\beta, t_2}| + \frac{1}{2} |t_1 - t_2|, 
	\end{align*}
	hence 
	\begin{align}
	\label{eq:75-gamma-star}
		|\gamma_{\beta, t_1} - \gamma_{\beta, t_2}| \leq |t_1 - t_2|. 
	\end{align}
\end{proof}

\begin{proposition}\label{prop:2}
	For any $\beta \geq \beta_0$ and $t \geq 0$, we have
	\begin{align*}
		\lim_{n \to \infty} \frac{1}{n} \E\left[\|\btheta - \bm(\by(t), t)\|_2^2 \right] = 1 - \frac{\gamma_{\beta, t}}{\beta^2}. 
	\end{align*}
\end{proposition}
\begin{proof}[Proof of Proposition \ref{prop:2}]
	This is a direct consequence of discussions in Section 2.4 of \cite{montanari2021estimation}. In particular, Proposition 2.2. 
\end{proof}

We finally notice that 
\begin{align}
\frac{1}{n} \E\Big\{\|\bm(\by(t), t) - \hbm^{k}(\by(t), t)\|_2^2\Big\}
= \frac{1}{n} \E\Big\{\|\btheta - \hbm^{k}(\by(t), t)\|_2^2\Big\}-\frac{1}{n} \E\Big\{\|\btheta - 
\bm(\by(t), t)\|_2^2\Big\}\, .
\end{align}
Therefore taking the limit $n\to\infty$, using Proposition \ref{prop:2},
 and noticing that all the expectations are 
of bounded random variables, we get
\begin{align}
\plim_{n\to\infty}\frac{1}{n} \|\bm(\by(t), t) - \hbm^{k}(\by(t), t)\|_2^2\le
\frac{\gamma_{\beta, t}}{\beta^2}-\frac{\gamma_t^k}{\beta^2}\, ,
\end{align}
which vanishes as $k \to \infty$.

\subsubsection{Proof of AMP consistency for symmetric $\pi_{\Theta}$}
\label{sec:ProofSymmetric}

We now consider the case of symmetric $\pi_{\Theta}$.
In this case, the posterior $\mu_{\bX,0}(\de\btheta)$ is symmetric
under flip $\btheta\to -\btheta$, and the original vector $\btheta$ is identifiable only 
up to a global sign.
We let $\bv_1= \bv_1(\bX)$ be a uniformly random eigenvector of $\bX$,
and denote by $\prob_0$ the joint distribution of $\btheta, \bX, \bv_1$.
(Since the top eigenvalue is almost surely non-degenerate, 
there are two possible choices for $\bv_1$
given $\bX$.)
We denote by $\prob_+$ the same distribution, conditioned to $\<\bv_1,\btheta\>>0$:
\begin{align}
\prob_+(\de\btheta,\de\bX,\de\bv_1)& = \frac{1}{\prob_0(\<\bv_1,\btheta\>>0)}
 \prob_0(\de\btheta,\de\bX,\de\bv_1) \,
 \bfone\{\langle \bv_1, \btheta \rangle \geq 0\}\\
 &= 2\, \prob_0(\de\btheta,\de\bX,\de\bv_1) \,
 \bfone\{\langle \bv_1, \btheta \rangle \geq 0\}\, .
\end{align}
Note that under $\prob_0$, $\bv_1$, and $\btheta$ are conditionally independent 
given $\bX$, while they are not under $\prob_+$. Also,
the marginal law  of $\bX,\bv_1$ is the same under the two distributions.

Conditionally on $\bX,\bv_1\sim \prob_{0}$, we let $\btheta^+_1,\btheta^+_2,\dots$ be
i.i.d. vectors with distribution  $\prob_+(\btheta\in\;\cdot\;|\bX,\bv_1)$
and $\btheta^0_1,\btheta^0_2,\dots$ be
i.i.d. vectors with distribution  $\prob_0(\btheta\in\;\cdot\;|\bX,\bv_1)$
(independent of the $\btheta^+_i$'s).
We will use the fact that, by an application of Remark \ref{lemma:nishimori},
if $\bX = \beta\btheta\btheta^{\sT}/n+\bW$,
\begin{align}
\Big(\bX,\bv_1,\btheta_1^0,\dots, \btheta_k^0,\btheta_1^+,\dots, \btheta_k^+\Big)
\ed \Big(\bX,\bv_1,\btheta,\btheta_1^0\dots, \btheta_{k-1}^0,\btheta_1^+,\dots, \btheta_k^+\Big)
\, .\label{eq:DistrSymmetry}
\end{align}
We first prove a concentration result  $\<\btheta^+_1, \btheta^+_2 \>$. 
\begin{lemma}\label{lemma:concentration-munew}
Let $\cuD(\beta)$ be the set of discontinuity points of $t\mapsto \gamma_{\beta, t}$, which we recall is the unique stationary point to \cref{eq:Phi}.
Then, for  all $t \in \reals_{\ge 0}\setminus\cuD(\beta)$, we have
	\begin{align}
&\lim_{n\to\infty}\frac{1}{n^2} \E[(\langle \btheta^+_1, \btheta^+_2 \rangle - \E[\langle \btheta_1^+, \btheta^+_2 \rangle])^2] = 0\, ,
\label{eq:ConcPlusPlus}\\
&\lim_{n \to \infty} 
\frac{1}{n} \E[\<\btheta^+_1, \btheta^+_2 \>] = 
\frac{\gamma_{\beta, t}  - t}{\beta^2}\, .\label{eq:PlusPlus}
	\end{align}
\end{lemma}
\begin{remark}
Note that $\cuD(\beta)$ is countable by monotonicity of the mapping $t \mapsto \gamma_{\beta, t}$, 
and further it is empty for all $\beta\ge \beta_0(\pi_{\Theta})$. Hence, in applying this lemma, we will disregard the set of
exceptional points $\cuD(\beta)$.
\end{remark}
\begin{proof}[Proof of \cref{lemma:concentration-munew}]
We divide the proof into two parts depending on the value of $t$.

\paragraph{Case I: ${\boldsymbol t > 0}$.}
We begin with a useful concentration result.
\begin{lemma}\label{lemma:proof-of-eq}
For all $0< t_1<t_2$, we have 
	\begin{align}\label{eq:orig-overlap-concentration}
	\lim_{n\to\infty}\frac{1}{n^2}\int_{t_1}^{t_2}
	 \E\Big[\Big(\langle \btheta^0_1, \btheta^0_2 \rangle - \E[\langle \btheta_1^0, \btheta^0_2 \rangle]\Big)^2\Big] \, \de t = 0\, .
	\end{align}
\end{lemma}
We present the proof of this fact in Appendix \ref{sec:proof-of-lemma:proof-of-eq}.
A similar statement is proven in \cite{lelarge2019fundamental}.

Let us next show  that this  implies the desired concentration result:
	\begin{align}\label{eq:dk}
 \E\left[ (\langle \btheta_1^+, \btheta^+_2 \rangle - \E[\langle \btheta_1^+, \btheta^+_2 \rangle])^2\right] \leq & 
  \E\left[ (\langle \btheta_1^+, \btheta_2^+ \rangle - \E[\langle\btheta_1^0, \btheta_2^0 \rangle] )^2 \right] \nonumber \\
= & \frac{\E\left[ (\langle \btheta^0_1, \btheta_2^0 \rangle - \E[\langle\btheta_1^0, \btheta_2^0 \rangle] )^2 
\bfone_{\< \bv_1, \btheta^1_0 \> \geq 0, \< \bv_1, \btheta_2^0 \> \geq 0)}\right]}
{\P(\< \bv_1, \btheta_1^0 \> \geq 0, \< \bv_1, \btheta_2^0 \> \geq 0)} \nonumber\\  	
\overset{(i)}{\leq} & 4
\E[(\< \btheta_1^0, \btheta_2^0 \> - \E[\< \btheta^0_1, \btheta^0_2 \>])^2] ,
	\end{align}
	Here in \emph{(i)} we made use of the fact that $\bv_1, \btheta_1^0 ,\btheta_2^0$
are conditionally independent given $\bX$, implying
	\begin{align}\label{eq:Porig-q}
	\P(\< \bv_1, \btheta_1^0 \> \geq 0, \< \bv_1, \btheta_2^0 \> \geq 0|\bX) = 
	 \P(\< \bv_1, \btheta_1^0 \> \geq 0|\bX)^2 = \frac{1}{4}\, ,
	\end{align}
	and therefore the same identity holds unconditionally. 
	
	
Recall that, 
by \cite{lelarge2019fundamental}, it holds that  (for any $t\ge 0$)
\begin{align}
\lim_{n \to \infty}\frac{1}{n^2}
\E[\<\btheta_1^0, \btheta_2^0\>^2] &= \lim_{n \to \infty}\frac{1}{n^2}
\E\Big\{\big\|\E\big[\btheta\btheta^{\sT}\big|\bX\big]\big\|_F^2\Big\}
=\frac{ (\gamma_{\beta, t} - t)^2}{\beta^4}\, .\label{eq:LimitSquare}
\end{align}
Using this, together with the concentration property \eqref{eq:orig-overlap-concentration},
we get, for all  $0< t_1<t_2$,
\begin{align*}
\lim_{n \to \infty}\frac{1}{n^2}\int_{t_1}^{t_2}
\Big(\E[\<\btheta_1^0, \btheta_2^0\>] - \frac{ \gamma_{\beta, t} - t}{\beta^2}\Big)^2 \de t = 0\, .
\end{align*}
Since $t\mapsto \E[\<\btheta_1^0, \btheta_2^0\>]=\E[\|\E[\btheta|\bX,\by(t)]\|^2]$
is non-decreasing (by Jensen),  the last limit holds pointwise.
Namely, at all continuity points of $t\mapsto \gamma_{\beta, t}$, 
\begin{align}
\lim_{n \to \infty}\frac{1}{n}\E[\<\btheta_1^0, \btheta_2^0\>] =
 \frac{ \gamma_{\beta, t} - t}{\beta^2}\, .\label{eq:Lim00}
 \end{align}
 Using this and \eqref{eq:LimitSquare}, we get, on $t\in (0,\infty)\setminus \cuD(\beta)$
 \begin{align}
 \lim_{n\to\infty}\frac{1}{n^2} \E[(\langle \btheta^0_1, \btheta^0_2 \rangle - \E[\langle \btheta_1^0, \btheta^0_2 \rangle])^2] = 0\, ,
 \end{align}
and therefore, using  \eqref{eq:dk}, we obtain the claim \eqref{eq:ConcPlusPlus}.

Finally, notice that the following  is a  consequence of \cref{eq:dk}:
\begin{align*}
\lim_{n \to \infty} \left(\frac{1}{n}\E[\langle \btheta_1^0, \btheta_2^0 \rangle] - \frac{1}{n}
\E[\langle \btheta_1^+, \btheta_2^+ \rangle] \right)^2 = 0.
\end{align*}
Putting the last limit together with Eq.~\eqref{eq:Lim00} implies Eq.~\eqref{eq:PlusPlus}.

\paragraph{Case II: ${\boldsymbol t = 0}$.} We begin with establishing the following lemma.

\begin{lemma}\label{lemma:inf}
Let $\beta_0(\pi_{\Theta})$ be as in the statement of Theorem \ref{thm:main2}.
Then for any $\beta\ge \beta_0(\pi_\Theta)$ and any $t\ge 0$, we have 
	\begin{align}
	\lim_{n\to\infty}\frac{1}{n^4}
	\E\Big[\Big(\langle\btheta^0_1, \btheta^0_2 \rangle^2 - \E[\langle\btheta^0_1, \btheta^0_2 \rangle^2]\Big)^2
	\Big] = 0\, .
	\end{align}
\end{lemma}

\begin{proof}
	Recall that $\Phi$ is defined in \cref{eq:Phi}. We let $\gamma_{\beta, t}$ be the first stationary point of $\gamma \mapsto \Phi(\gamma, \beta, t)$ on $(t, \infty)$. Following the notation of \cite{lelarge2019fundamental}, we let
	\begin{align*}
		D_t := \left\{ \beta > 0: \gamma \mapsto \Phi(\gamma, \beta, t) \mbox{ has a unique minimizer} \right\}.
	\end{align*}
	By the assumptions of Theorem \ref{thm:main2}, we know that 
	$[\beta_0(\pi_{\Theta}), \infty) \subseteq D_t$. Then, the claim of the lemma is
	 a direct consequence of \cite[Theorem 20]{lelarge2019fundamental}.  
\end{proof}

By \cref{lemma:inf}, and repeating the argument of Eq.~\eqref{eq:dk}, we get
\begin{align}
 \E\left[ (\langle \btheta_1^+, \btheta^+_2 \rangle^2 - \E[\langle \btheta_1^+, \btheta^+_2 \rangle^2])^2\right] \leq & 
  \E\left[ (\langle \btheta_1^+, \btheta_2^+ \rangle^2 - \E[\langle\btheta_1^0, \btheta_2^0 \rangle^2] )^2 \right] \nonumber \\
=& \frac{\E\left[ (\langle \btheta^0_1, \btheta_2^0 \rangle^2 - \E[\langle\btheta_1^0, \btheta_2^0 \rangle^2] )^2 
\bfone_{\< \bv_1, \btheta^1_0 \> \geq 0, \< \bv_1, \btheta_2^0 \> \geq 0)}\right]}
{\P(\< \bv_1, \btheta_1^0 \> \geq 0, \< \bv_1, \btheta_2^0 \> \geq 0)} \nonumber\\  	
\leq & 4
\E[(\< \btheta_1^0, \btheta_2^0 \>^2 - \E[\< \btheta^0_1, \btheta^0_2 \>^2])^2]\, .
\label{eq:ChainSquare}
	\end{align}
Therefore, by the last lemma, 
\begin{align}
	\lim_{n\to\infty}\frac{1}{n^4}
	\E\Big[\Big(\langle\btheta^+_1, \btheta^+_2 \rangle^2 - \E[\langle\btheta^+_1, \btheta^+_2 \rangle^2]\Big)^2
	\Big] = 0\, .\label{eq:ConcSquare}
	\end{align}

Let $\bnu_0 := \sqrt{n}\bv_1(\bX) / \|\bv_1(\bX)\|_2 $. 
By \cite{benaych2011eigenvalues}, we know that $\langle \btheta, \bnu_0 \rangle^2 / n^2 
\ed \langle \btheta^0_1, \bnu_0 \rangle^2 / n^2 	\toP 1 - \beta^{-2}$. 
Since $\prob_+$ is contiguous to $\prob_0$, we obtain $\langle \btheta^+_1, \bnu_0 \rangle / n 
 \toP \sqrt{1 - \beta^{-2}}$. and therefore:
\begin{align}\label{eq:spectrals}
\plim_{n\to\infty} \frac{1}{n}\|\btheta^+_1 - \bnu_0\|_2^2 = 2 - 2 \sqrt{1 - \beta^{-2}} .
\end{align}
Recall that $\gamma_{\beta, t}$ is the first positive stationary 
point of $\gamma \mapsto \Phi(\gamma, \beta, t)$. From Eqs.~\eqref{eq:ChainSquare} and
\eqref{eq:ConcSquare}, we see that 
$|\langle \btheta^+_1, \btheta^+_2 \rangle| / n = \E[\langle \btheta, \btheta^0_1 \rangle^2]^{1/2} / n 
+ o_P(1)$. 
Further by\cite[Theorem 2]{lelarge2019fundamental},  we obtain  
$|\langle \btheta, \btheta^0_1 \rangle| / n = \beta^{-2} \gamma_{\beta, 0} + o_P(1)$, 
whence $|\langle \btheta^+_1, \btheta^+_2 \rangle| / n = \beta^{-2} \gamma_{\beta, 0} + o_P(1)$. 

By Cauchy-Schwarz,
\begin{align*}
	\frac{1}{n}\|\btheta^+_1 - \btheta^+_2\|_2^2 \leq \frac{2}{n}\|\btheta_1^+ - \bnu_0\|_2^2 + \frac{2}{n}\|\btheta^+_2 - \bnu_0\|_2^2 = 4 - 4 \sqrt{1 - \beta^{-2}} + o_P(1), 
\end{align*}
hence $\langle \btheta^+_1, \btheta^+_2 \rangle / n \geq 2 \sqrt{1 - \beta^{-2}} - 1 + o_P(1)$. 
Recall that $|\langle\btheta^+_1, \btheta^+_2 \rangle| / n = \beta^{-2} \gamma_{\beta, 0} + 
o_P(1)$, then for $\beta_0$ large enough and all $\beta > \beta_0$, it holds that 
$\langle \btheta^+_1, \btheta^+_2 \rangle / n = \beta^{-2} \gamma_{\beta, 0} + o_P(1)$.
 Applying bounded convergence (since $|\langle \btheta^+_1, \btheta^+_2 \rangle / n|\le M_{\Theta}$), 
 we see that $\E[(\langle \btheta^+_1, \btheta^+_2 \rangle / n - \beta^{-2} \gamma_{\beta, 0})^2] = o_n(1)$,
  thus concluding the proof of \cref{lemma:concentration-munew} for $t = 0$. 
\end{proof}
%

Next, we will apply \cref{lemma:concentration-munew} to prove
  \cref{lemma:AMP-optimal}.  By the state evolution of the AMP algorithm,
  cf. Proposition \ref{propo:SE-Basic}, we see that 
\begin{align*}
	& \frac{1}{n} \langle \btheta^+, \hat{\bm}^k(\by(t), t) \rangle \toP \E[\E[ \Theta \mid \gamma_t^k \Theta + (\gamma_t^k)^{1/2} G]^2] = 1 - \mathsf{mmse}( \gamma_t^k), \\
	& \frac{1}{n} \| \hat{\bm}^k( \by(t), t) \|_2^2  \toP \E[\E[ \Theta \mid \gamma_t^k \Theta +  (\gamma_t^k)^{1/2} G]^2] = 1 - \mathsf{mmse}( \gamma_t^k). 
\end{align*} 
By Proposition \ref{prop:1}, we see that as $k \to \infty$, $\gamma_t^k$ converges linearly 
to $\gamma_{\beta, t}$, which further implies that $1 - \mathsf{mmse}(\gamma_t^k)$
 converges linearly to $1 - \mathsf{mmse}(\gamma_{\beta, t})
  = \beta^{-2} (\gamma_{\beta, t} - t)$. Furthermore, the convergence is uniform in 
  $t \in [0, T]$. Therefore, for all $\ep > 0$, there exists $K(\beta,T,\eps) \in \NN_{>0}$ 
  depending only on $(\beta, T, \eps, \pi_{\Theta})$, such that for all $k \geq K(\beta,T,\eps)$, 
\begin{align*}
	\big|\E[\E[ \Theta \mid \gamma_t^k \Theta + (\gamma_t^k)^{1/2} G]^2] -  \beta^{-2} (\gamma_{\beta, t} - t)\big| \leq \frac{\ep}{2}.
\end{align*}
and therefore,  for all $k \geq K(\beta,T,\eps)$, with high probability,
\begin{align}
&\left| \frac{1}{n} \< \btheta^+, \hat{\bm}^k(\by(t), t) \> - \beta^{-2} (\gamma_{\beta, t} - t)\right|
\le \eps\, ,\label{eq:ThetaPlus}\\
&\left| \frac{1}{n} \| \hat{\bm}^k( \by(t), t) \|_2^2 - \beta^{-2} (\gamma_{\beta, t} - t)\right|
\le \eps\, .\label{eq:First}
\end{align} 
By \cref{lemma:concentration-munew}, it holds that  
$\plim_{n\to\infty}
\langle \btheta_1^+, \btheta_2^+ \rangle / n =  \beta^{-2} (\gamma_{\beta, t} - t)$.
 Since
 $\btheta^+_1$ and $\btheta^+_2$ are conditionally independent given $(\bX, \by(t), \bv_1(\bX))$
 this in particular implies:
 \begin{align}
 	\frac{1}{n} \|\bm(\by(t), t)\|_2^2  =  \beta^{-2} (\gamma_{\beta, t} - t) + o_P(1)\, .
 	\label{eq:Second}
 \end{align}
Further $\btheta^+$ and $\bm(\by(t), t)$ are conditionally independent given 
$(\bX, \by(t), \bv_1(\bX))$. Hence, from Eq.~\eqref{eq:ThetaPlus}, it follows that
with high probability
 \begin{align}
 \left| \frac{1}{n} \langle \bm(\by(t), t), \hat{\bm}^{K(\beta,T,\eps)}( \by(t), t)  \rangle  
 - \beta^{-2} (\gamma_{\beta, t} - t)\right| \le 2\eps\, .\label{eq:Third}
\end{align}
Putting together \eqref{eq:First}, \eqref{eq:Second}, \eqref{eq:Third},
we get
\begin{align*}
	& \frac{1}{n} \|\bm(\by(t), t) - \hat{\bm}^{K(\beta,t,\eps)}(\by(t), t)\|_2^2  \le 10\, \eps\, ,
\end{align*}
with high probability,
thus completing the proof of \cref{lemma:AMP-optimal}.

\subsection{Proof of \cref{lemma:path-regularity}}\label{sec:proof-of-lemma:path-regularity}

The subsequent proof is analogous to the one of \cite[Lemma 4.9]{alaoui2022sampling}. 
Recall that $\Phi$ is defined in \cref{eq:Phi}, and 
(for $\beta>\beta_0(\pi_{\Theta})$), $\gamma_{\beta, t}$ is the 
unique global maximizer of $\gamma \mapsto \Phi(\gamma, \beta, t)$ over $\gamma \in (t, \infty)$. 
Taking the partial derivative of $\Phi(\gamma, \beta, t)$ with respect to $\gamma$, we obtain that 
\begin{align*}
	\frac{\partial}{\partial \gamma} \Phi(\gamma, \beta, t) = 
	\frac{\gamma - t}{2\beta^2} - \frac{1}{2} + \frac{1}{2}\mathsf{mmse}(\gamma),
\end{align*}
where we recall that, for $(\Theta,G) \sim \pi_{\Theta}\otimes \normal(0,1)$,
\begin{align*}
	\mmse(\gamma) = \E \big[(\Theta - \E[\Theta \mid \gamma \Theta + \sqrt{\gamma} G])^2\big].
\end{align*}
Therefore, $\gamma_{\beta, t}$ is a solution to the following fixed point equation. 
\begin{align}\label{eq:72}
	\gamma_{\beta, t}  = \beta^2\E[\E[\Theta \mid \gamma_{\beta, t}\Theta + \sqrt{\gamma_{\beta, t}}G]^2] + t. 
\end{align} 
For any $t_1 < t_2$,  we have
\begin{align}
	& \lim_{n \to \infty} \frac{1}{n} \E\left[ \|\bm(\by(t_2), t_2)    - \bm(\by(t_1), t_1)\|_2^2  \right] =\nonumber  \\
	 = & \lim_{n \to \infty} \frac{1}{n} \left\{ \E\left[ \|\btheta - \bm(\by(t_1), t_1)\|_2^2  \right] - \E\left[ \|\btheta - \bm( \by(t_2), t_2)\|_2^2 \right]\right\} \nonumber \\
	= & \frac{\gamma_{\beta, t_2} - t_2 - \gamma_{\beta, t_1} + t_1}{\beta^2}. \label{eq:73}
\end{align}
By \cref{lemma:AMP-optimal} we know that for all $t \geq 0$, with high probability $\|\bm(\by(t), t) - \hat{\bm}^k(\by(t), t)\|_2^2 / n \leq \ep_k$, for some deterministic constants $\ep_k$ satisfying $\ep_k \to 0^+$ as $k \to \infty$.
 Therefore, using the concentration of 
 $\|\hbm^k(\by(t_2), t_2) - \bm^k(\by(t_1), t_1)\|_2^2/n$, we get
\begin{align}
	\plim\limits_{n \to \infty} \frac{1}{n} \|\bm(\by(t_2), t_2) - \bm(\by(t_1), t_1)\|_2^2 = \frac{\gamma_{\beta, t_2} - t_2 - \gamma_{\beta, t_1} + t_1}{\beta^2}.\label{eq:PconvMt} 
\end{align}
Note that $t \mapsto \bm(\by(t), t)$ is a martingale. Hence, for any fixed constant $c$,
 the process $Y_{n, t} := (M_{n, t} - c)_+$ is a positive submartingale, where
  $M_{n,t} = \|\bm(\by(t), t) - \bm(\by(t_1), t_1)\|_2 / \sqrt{n}$. By
   Doob's maximal inequality, we then see that
\begin{align*}
	\P\left( \sup_{t \in [t_1, t_2]} Y_{n,t} \geq a \right) \leq \frac{1}{a}\E[Y_{n, t_2}] \leq \frac{1}{a} \E[Y_{n, t_2}^2]^{1/2}
\end{align*}
for any $a > 0$. Setting $c = \sqrt{\gamma_{\beta, t_2} - t_2 - \gamma_{\beta, t_1} + t_1} / \beta$,
 we have $\plim_{n \to \infty} M_{n, t_2}^2 = c^2$
 by Eq.~\eqref{eq:PconvMt}. 
 Then for any fixed $a > 0$, we obtain:
\begin{align*}
	\limsup_{n \to \infty} \P\left( \sup_{t \in [t_1, t_2]} M_{n, t} \geq c + a \right) \leq & \limsup\limits_{n \to \infty} \P\left( \sup_{t \in [t_1, t_2]} Y_{n,t} \geq a \right) \\
	\leq & \frac{1}{a} \lim_{n \to \infty} \E\left[ (M_{n, t_2} - c)^2 \right]^{1/2} = 0. 
\end{align*}
A lower bound can be derived analogously. Thus, 
\begin{align*}
	\plim_{n \to \infty}\sup_{t \in [t_1, t_2]} M_{n,t}^2 = \frac{\gamma_{\beta, t_2} - t_2 - \gamma_{\beta, t_1} + t_1}{\beta^2},
\end{align*}
which yields
\begin{align*}
	\plim_{n \to \infty}\sup_{t \in [t_1, t_2]} \frac{1}{n}\|\bm(\by(t), t) - \bm(\by(t_1), t_1))\|_2^2 = & \plim_{n \to \infty} \frac{1}{n} \|\bm( \by(t_2), t_2) - \bm( \by(t_1), t_1)\|_2^2 \\
	=& \frac{\gamma_{\beta, t_2} - t_2 - \gamma_{\beta, t_1} + t_1}{\beta^2}.
\end{align*}
Therefore, in order to prove the lemma, it suffices to show the existence of $C_{\sreg} > 1$ depending uniquely on $\beta$, such that 
\begin{align*}
	\frac{|\gamma_{\beta, t_2} - t_2 - \gamma_{\beta, t_1} + t_1|}{\beta^2} \leq C_{\sreg}|t_1 - t_2|, 
\end{align*}
which follows from Proposition \ref{prop:1}. This concludes the proof of the lemma.

\subsection{Proof of \cref{lemma:general-AMP-hits-strong-signal-region}}\label{sec:proof-of-lemma:general-AMP-hits-strong-signal-region}

Before proving  \cref{lemma:general-AMP-hits-strong-signal-region}, we establish a 
simple estimate on the conditional variance.
\begin{lemma}\label{lemma:discrete-MSE}
	There exists a constant $C_{\sconv} > 0$ depending only on $\pi_{\Theta}$, such that 
	\begin{align}\label{eq:discrete-MSE}
		\E\big[\Var(\Theta \mid \beta^2  \Theta   + \beta G )\big] \leq 
		C_{\sconv}^{-1}\exp(-4C_{\sconv} \beta^2) / 2. 
	\end{align}
	Without loss, we can and will assume that $C_{\sconv} < 1$. 
\end{lemma}
\begin{proof}
We denote by $\{x_1, x_2, \cdots, x_s\}$ the support of $\pi_{\Theta}$ 
and assume without loss of generality $x_1 < x_2 < \cdots < x_s$. Define 
$\htheta: \RR \to \{x_1, x_2, \cdots, x_s\}$ by
\begin{align*}
	\htheta(y):= \argmin\big(|x-y|:\; x\in\supp(\pi_{\Theta})\big)\, .
\end{align*}
In case of ties, we choose  the smallest value. We immediately see that 
$\htheta(y) = x_i$ if and only if $(x_{i - 1} + x_i) / 2 < y \leq (x_i + x_{i + 1}) / 2$
 (with the convention that $x_0 = -\infty$ and $x_{s + 1} = +\infty$). Let
  $Y = \Theta + \beta^{-1} G$, then
\begin{align*}
	\E[\Var[\Theta \mid \beta^2  \Theta  + \beta G ]]
	\leq & \E[(\Theta - \htheta(Y))^2].
\end{align*}
Let $\delta_\Theta = \min\{|x_i - x_{i + 1}| / 2: i \in [s - 1]\}$. We then have 
\begin{align*}
	& \E[(\Theta - \htheta(Y))^2 \mid \Theta = x_i] \\
	 \leq & 0 \times \P\left( Y \in \Big( (x_{i - 1} + x_i) / 2, (x_i + x_{i + 1})  / 2 \Big] \mid \Theta = x_i\right)  \\
	 &+ 4M_{\Theta}^2 \P \left( Y \in \Big( (x_{i - 1} + x_i) / 2, (x_i + x_{i + 1})  / 2 \Big]^c \mid \Theta = x_i\right) \\
	 \leq &  \frac{16M_{\Theta}^2}{\delta_{\Theta} \beta \sqrt{2\pi}} e^{-\delta_{\Theta}^2 \beta^2 / 8},
\end{align*}
where to arrive at the last inequality we make use of \cref{lemma:gaussian-tail}. 
Combining the above bounds, we obtain that $\E[(\Theta - \htheta(Y))^2] \leq  
\frac{16M_{\Theta}^2}{\delta_{\Theta} \beta \sqrt{2\pi}} e^{-\delta_{\Theta}^2 \beta^2 / 8}$.
 Using this fact, we conclude that there exists a constant $C_{\sconv} > 0$ that is 
 a function of $\pi_{\Theta}$ only, such that \cref{eq:discrete-MSE} holds. 
\end{proof}

\begin{proof}[Proof of \cref{lemma:general-AMP-hits-strong-signal-region}]
By the state evolution of Bayes AMP, Proposition \ref{propo:SE-Basic}, we have
\begin{align}\label{eq:68}
	& \frac{1}{n}\|\bD_{\gamma_{t}^{k}}(\hbm^k(\by(t),t))\|_F^2 \toP \E[\Var[\Theta \mid \gamma_t^k \Theta +  (\gamma_t^k)^{1/2} G]]. 
\end{align}
Since $\gamma_t^0 = \beta^2 - 1$ and  $\gamma_t^k \geq \gamma_t^0$
(this follows from instance by the fact that
$\gamma\mapsto \mmse(\gamma)$ is non-increasing, see the discussion at the beginning
of Section \ref{sec:ProofNonSymmetric}), we conclude that for $\beta>2$,
\begin{align}\label{eq:69}
	\E[\Var[\Theta \mid  \gamma_t^k \Theta + (\gamma_t^k)^{1/2} G]] \leq \E[\Var[\Theta \mid \beta^2  \Theta / 4  + \beta G / 2]].
\end{align}
By \cref{lemma:discrete-MSE} below, we obtain that there exists a constant $C_{\sconv} > 0$
 depending uniquely on $\pi_{\Theta}$, such that 
\begin{align}\label{eq:70}
	\E[\Var[\Theta \mid \beta^2  \Theta / 4  + \beta G / 2]] \leq C_{\sconv}^{-1}\exp(-C_{\sconv} \beta^2) / 2.
\end{align}
Equation \eqref{eq:general-m} follow from \cref{eq:68,eq:69,eq:70}. 

By definition $b_t^k = \beta^2 \E[\Var[\Theta \mid \gamma_t^k \Theta + (\gamma_t^k)^{1/2} G]]$, 
which by \cref{eq:69,eq:70} is no larger than $C_{\sconv}^{-1}\exp(-C_{\sconv} \beta^2) / 2$, thus completing the proof of \cref{eq:general-b}. 

As for \cref{eq:general-Lipschitz1,eq:general-Lipschitz2}, we will in fact show a stronger
 result and prove that these two inequalities hold for all $k \leq k_0(\beta) + 1$,
  via induction over $k$. We already observed that, with probability $1 - o_n(1)$ we have 
$\|\bX\|_{\op} \leq \beta+\|\bW\|_{\op}\le \beta+2$ \cite{Guionnet}, and will work
on this high-probability event.

  For the base case $k = 0$, the claim directly follows as 
  $\hat{\bm}^0(\by_1, t) = \hat{\bm}^0(\by_2, t) = \E[\Theta \mid \gamma^0_t \Theta + 
  (\gamma^0_t)^{1/2} G = \bnu]$ and hence $\hat{\bp}^0(\by_1, t) = \hat{\bp}^0(\by_2, t)$. 
  Now suppose for all $k \leq k_1$ and all $\by_1$, $\by_2$, we have
\begin{align*}
		& \frac{1}{\sqrt{n}} \|\hat{\bm}^{k}(\by_1, t) - \hat{\bm}^{k}(\by_2, t)  \|_2 \leq \frac{\Lipc({\beta, k})}{\sqrt{n}} \|\by_1 - \by_2\|_2,  \\
		& \frac{1}{\sqrt{n}} \| \hat{\bp}^{k}(\by_1, t) - \hat{\bp}^{k}(\by_2, t)  \|_2 \leq \frac{\Lipc({\beta, k})}{\sqrt{n}} \|\by_1 - \by_2\|_2, 
\end{align*}
where $\Lipc({\beta, k}) > 0$ is a function of $(\beta, k)$ only. We then prove that the above statement also holds for $k = k_1 + 1$. 

Recalling the definition $\sF(x;\gamma) :=  \E[\Theta \mid \gamma \Theta + \gamma^{1/2}G = z]$,
a computation of the derivatives shows that the mappings 
\begin{align*}
 z \mapsto \sF(z;\gamma_t^{k + 1}), \;\;\;\;\;\
	& z \mapsto \Phi_{\alpha_t^{k + 1}}^{-1}\big(\sF(z;\gamma_t^{k + 1})\big)\, ,
\end{align*}
are $M_{\Theta}^2$-Lipschitz and $M_{\Theta}$-Lipschitz, respectively, where
$\supp(\pi_{\Theta})\subseteq [-M_{\Theta},M_{\Theta}]$.
 As a result, on the event $\|\bX\|_{\op} \le \beta+2$, we have
\begin{align*}
	& \frac{1}{\sqrt{n}}\|\hat{\bm}^{k_1 + 1}(\by_1, t_1) - \hat{\bm}^{k_1 + 1}(\by_2, t_2)\|_2 \\
	 \leq & \frac{M_{\Theta}^2}{\sqrt{n}} \|\beta \bX(\hat{\bm}^{k_1}(\by_1, t) - \hat{\bm}^{k_1}(\by_2, t)) + \by_1 - \by_2 - b_t^{k_1} (\hat{\bm}^{k_1 - 1}(\by_1, t) - \hat{\bm}^{k_1 - 1}( \by_2, t))\|_2 \\
	 \leq & \frac{M_{\Theta}^2(\beta^2 + \beta + 1)\Lipc({\beta, k_1})}{\sqrt{n}}\|\by_1 - \by_2\|_2 + \frac{M_{\Theta}^2}{\sqrt{n}}\|\by_1 - \by_2\|_2 + \frac{C_{\sconv}^{-1} M_{\Theta}^2 \Lipc({\beta, k_1 - 1})}{\sqrt{n}}\|\by_1 - \by_2\|_2. 
\end{align*}
Similarly, under the $\bp$-parameterization we have 
\begin{align*}
	& \frac{1}{\sqrt{n}}\|\hat{\bp}^{k_1 + 1}(\by_1, t) - \hat{\bp}^{k_1 + 1}( \by_2, t)\|_2 \\
	 \leq & \frac{M_{\Theta}}{\sqrt{n}} \|\beta \bX(\hat{\bm}^{k_1}(\by_1, t) - \hat{\bm}^{k_1}(\by_2, t)) + \by_1 - \by_2 - b_t^{k_1} (\hat{\bm}^{k_1 - 1}(\by_1, t) - \hat{\bm}^{k_1 - 1}(\by_2, t))\|_2 \\
	 \leq & \frac{M_{\Theta}(\beta^2 + \beta + 1)\Lipc({\beta, k_1})}{\sqrt{n}}\|\by_1 - \by_2\|_2 + \frac{M_{\Theta}}{\sqrt{n}}\|\by_1 - \by_2\|_2 + \frac{C_{\sconv}^{-1} M_{\Theta} \Lipc({\beta, k_1 - 1})}{\sqrt{n}}\|\by_1 - \by_2\|_2. 
\end{align*}
As a result, we see that setting $\Lipc({\beta, k_1 + 1}) = M_{\Theta}^2((\beta^2 + \beta + 1) \Lipc({\beta, k_1}) + 1 + C_{\sconv}^{-1} \Lipc({\beta, k_1 - 1}))$ concludes the proof of the induction step. This further completes the proof of \cref{eq:general-Lipschitz1} and \cref{eq:general-Lipschitz2}, thus finishing the proof of the lemma.

\end{proof}

\subsection{Proof of \cref{lemma:MainLip}}\label{sec:proof-of-lemma:main}

We first state a simplified version of \cref{lemma:C2}. More precisely, for $q \in (0,1)$, note that $\sqrt{q \log (e / q)} \leq 3q^{1/4}$. If we substitute this result 
into \cref{eq:C2}, normalize $\bt_1, \bt_2$, and set $\xi = \Delta^{1/6}$, $q = \Delta^{2/3}$, then we obtain the next corollary.   
\begin{corollary}
Under the conditions of \cref{lemma:C2}, for all $\Delta > 0$ and $M > 0$, we have
	\begin{align}\label{eq:C2-additional}
		\P\left( \sup_{\substack{\bt_1, \bt_2 \in [0,M]^n,\\
		 \|\bt_1\|_2^2 / n \leq \Delta, \|\bt_2\|_2^2 / n \leq \Delta}}  \|\diag(\bt_1) \bW \diag(\bt_2)\|_{\op} \geq 4C' M^{5/3} \Delta^{1/6} \right) \leq Ce^{-cn \Delta^{2/3}M^{-4/3}}.
	\end{align}
\end{corollary}
We denote by $\cuE_{\beta, n}^{(4)}$ the event depicted by \cref{eq:C2-additional}, with 
$\Delta = \Delta(\beta)$. More precisely, let
\begin{align}\label{eq:E4}
	\cuE_{\beta, n}^{(4)} := \left\{ \sup_{\bt_1, \bt_2 \in [0,M_{\Theta}]^n, \|\bt_1\|_2^2 / n \leq \Delta(\beta), \|\bt_2\|_2^2 / n \leq \Delta(\beta)}  \|\diag(\bt_1) \bW \diag(\bt_2)\|_{\op} \geq 4C' M_{\Theta}^{5/3} \Delta(\beta)^{1/6} \right\}.
\end{align}
Throughout the proof, we will make use of the following functions:
\begin{align}\label{eq:functions}
	& \Delta(\beta) := C_{\sconv}^{-1}e^{-C_{\sconv} \beta^2} + M_{\Theta}^2 \cdot (2\beta^2 + 2\beta + 4 + 2C_{\sconv}^{-1}e^{-C_{\sconv} \beta^2})e^{ -C_{\sconv}\beta^2 / 4 }, \nonumber \\
	& \rho(\beta) := \beta^2 M_{\Theta}^2 \Delta(\beta) + 4\beta C' M_{\Theta}^{5/3} \Delta(\beta)^{1/6}, \\
	& F(\beta) := \rho(\beta) +  M_{\Theta}^2 \cdot C_{\sconv}^{-1}e^{-C_{\sconv} \beta^2}.\nonumber
\end{align}
We can and will choose $\beta_0$ large enough such that $F(\beta) \leq 1 / 2$ holds for all 
$\beta \geq \beta_0$. To simplify notations, we define
\begin{align*}
	& \bm^k = \hat{\bm}^k(\by_1, t), \\
	& \tilde{\bm}^k = \hat{\bm}^k( \by_2, t), \\
	& \bp^k = \Psi^{-1}_{\gamma_t^k}(\bm^k), \\
	& \tilde{\bp}^k = \Psi^{-1}_{\gamma_t^k}(\tilde{\bm}^k).
\end{align*}
We will choose $r(\beta)$ (depending uniquely on $\pi_{\Theta},\beta$)
small enough so that $2r(\beta) \cdot (\Lipc(\beta) + 1)\cdot (M_{\Theta}^2 + 1) \leq 2 
e^{- C_{\sconv} \beta^2 / 4}$. Notice that indeed the choice of $r(\beta)$ can only depend on $(\beta, \pi_{\Theta})$. 
By \cref{lemma:general-AMP-hits-strong-signal-region}, 
for all $\by_1,\by_2\in \Ball^n(\by(t),r(\beta))$, 
we know that
\begin{align}\label{eq:diff-p}
	\frac{1}{\sqrt{n}}\|\bp^{k_{\ast}} - \tilde{\bp}^{k_{\ast}}\|_2 \leq \frac{\Lipc(\beta)}{\sqrt{n}} \|\by_1 - \by_2\|_2 \leq 2\exp\left(- \frac{1}{4}C_{\sconv} \beta^2\right)
\end{align} 
for $k_{\ast} \in \{ k_0(\beta), k_0(\beta) \pm 1\}$. Without loss, we can and will assume that $\Lipc(\beta) \geq 2M_{\Theta}$. 

We define $\cuE_{\beta, L, \delta, \ep, n} = \cuE^{(1)}_{L, \delta, \eps, n} \cap \cuE_{\beta, L, \delta, n}^{(2)} \cap \cuE^{(3)}_{\beta, L, \delta, \eps, n}$.
The subsequent proof will be based on the following lemma:
\begin{lemma}\label{lemma:support-main}
On the set $\cuE_{\beta, L, \delta, \ep, n}$, if in addition we have
\begin{align}\label{eq:diff-pB}
	\frac{1}{\sqrt{n}}\|\bp^{k_{\ast}} - \tilde{\bp}^{k_{\ast}}\|_2 \leq \frac{\Lipc(\beta)}{\sqrt{n}} \|\by_1 - \by_2\|_2  \leq 2\exp\left(- \frac{1}{4}C_{\sconv} \beta^2\right)
\end{align}
holds for all $k_{\ast} \in \{k, k + 1, k + 2\}$ with $k_0(\beta) - 1 \leq k \leq K(\beta, T, \eps) - 3$, then it also holds for $k_{\ast} = k + 3$. Furthermore, the following inequality holds for all $k_0(\beta) - 1 \leq k \leq K(\beta, T, \eps) - 3$: 
\begin{align}\label{eq:49}
	\frac{1}{\sqrt{n}} \|\bp^{k + 3} - \tilde{\bp}^{k + 3}\|_2 \leq \frac{\rho(\beta)}{\sqrt{n}} \|\bp^{k + 2} - \tilde{\bp}^{k + 2}\|_2 + \frac{\rho(\beta)}{\sqrt{n}} \|\bp^{k + 1} - \tilde{\bp}^{k + 1}\|_2 + \frac{1}{\sqrt{n}}\|\by_1 - {\by_2}\|_2, 
\end{align}
where we recall that $\rho$ is defined in \cref{eq:functions}.
\end{lemma}
\cref{lemma:support-main} and \cref{eq:diff-p} 
imply the following upper bound via induction argument: 
\begin{align}\label{eq:LipAlmostFinal}
	\frac{1}{\sqrt{n}} \|\bp^{K(\beta, T, \eps)} - \tilde\bp^{K(\beta, T, \eps)}\|_2 \leq \frac{1 + 2\Lipc(\beta)}{1 - 2\rho(\beta)} \times \frac{1}{\sqrt{n}}\|\by_1 - \by_2\|_2. 
\end{align}
Define
\begin{align}\label{eq:Sbeta}
	\Lips(\beta) := \frac{M_{\Theta} \cdot (1 + 2\Lipc(\beta) )}{1 - 2\rho(\beta)}.
\end{align}
The claim of the lemma follows from \cref{eq:LipAlmostFinal} using
the fact that $\Psi_{\alpha_t^{K(\beta, T, \eps)}}$ has Lipschitz 
constant $M_{\Theta}$.

The remainder of this section is dedicated to proving \cref{lemma:support-main}.

\begin{proof}[Proof of \cref{lemma:support-main}]
The condition of \cref{lemma:support-main} assumes that for all $k_{\ast} \in \{k, k + 1, k + 2\}$, 
\begin{align}\label{eq:diff-m}
	\frac{1}{\sqrt{n}} \|\bp^{k_{\ast}} - \tilde{\bp}^{k_{\ast}}\|_2 \leq 2\exp\left(- \frac{1}{4}C_{\sconv} \beta^2\right). 
\end{align}
Next, we will make use of the Jacobian matrices given in \cref{eq:general-jacob-p} to provide a preliminary upper bound for $\|\bp^{k + 3} - \tilde{\bp}^{k + 3}\|_2$. 

On the set $\cuE_{\beta, L, \delta, \ep, n}$, for all $\bm, \bm' \in [a_{\Theta},b_{\Theta}]^n$ and $k_0(\beta) \leq k \leq K(\beta, T, \eps) - 3$, it holds that 
\begin{align*}
	& \|\beta \bD(\bm) \bX \bD(\bm')\|_{\op} \leq M_{\Theta}^2 \cdot ( \beta^2 + 2\beta), 
	\qquad \|b_t^k \bD(\bm) \bD(\bm')\|_{\op} \leq M_{\Theta}^2 \cdot C_{\sconv}^{-1}\exp(-C_{\sconv} \beta^2), \\
	& \|\bD(\bm)\|_{\op} \leq M_{\Theta},
\end{align*}
where we used \cref{eq:general-b}.  
Combining these upper bounds and \cref{eq:general-jacob-m}, we obtain a crude upper bound for $\|\bp^{k + 3} - \tilde{\bp}^{k + 3}\|_2$:
\begin{align*}
	& \frac{1}{\sqrt{n}} \|\tilde{\bp}^{k + 3} - \bp^{k + 3}\|_2 \\
	 \leq & \frac{M_{\Theta}^2 \cdot (\beta^2 + \beta + 1)}{\sqrt{n}} \|\tilde\bp^{k + 2} - {\bp}^{k + 2}\|_2 + \frac{M_{\Theta}^2 \cdot C_{\sconv}^{-1}\exp(-C_{\sconv} \beta^2)}{\sqrt{n}} \|\tilde{\bp}^{k + 1} - \bp^{k + 1}\|_2 \\
	 & + \frac{M_{\Theta}^2}{\sqrt{n}}\|\by_1 - {\by}_2\|_2 \\
	\overset{(i)}{\leq} &  M_{\Theta}^2 \cdot (2\beta^2 + 2\beta + 4 + 2C_{\sconv}^{-1}\exp(-C_{\sconv} \beta^2))\exp\left( -\frac{1}{4}C_{\sconv}\beta^2 \right),
\end{align*}
where to obtain \emph{(i)}, we use the following facts: (1) $\|\tilde{\bp}^{k + i} - \bp^{k + i}\|_2 
/ \sqrt{n} \leq 2 e^{-C_{\sconv} \beta^2 / 4}$ for all $i \in \{1,2\}$; (2) 
$\|\by_1 - \by_2\|_2 / \sqrt{n} \leq 2r(\beta) \leq 2 e^{-C_{\sconv} \beta^2 / 4}$. 
Since $\by(t) \in \Ball^n(\by(t),r(\beta))$, we can also control the difference between $\bp^{k + 3}$, $\tilde\bp^{k + 3}$ and  $\hat{\bp}^{k + 3}( \by(t), t )$ following exactly the same manner, and produce exactly the same upper bound. 

Before completing the proof, it
is useful to establish the following lemma. 
\begin{lemma}\label{lemma:D-p}
For any $\pi_{\Theta}$ such that $\supp(\pi_{\Theta})\subseteq[-M_{\Theta},M_{\Theta}]$ and
any $\gamma > 0$, the mapping 
	\begin{align*}
		Q(p) := \Var[\Theta \mid \gamma \Theta + \sqrt{\gamma} G = \Gamma_{\gamma}^{-1}(p)]
	\end{align*}
	is $3M_{\Theta}^2$-Lipschitz continuous. 
\end{lemma}
\begin{proof}[Proof of \cref{lemma:D-p}]
Let $h = \Gamma_{\gamma}^{-1}(p)$. Taking the derivative of $Q(\cdot)$, we obtain via chain rule
\begin{align*}
	& \frac{\dd Q}{\dd p} =  \frac{\dd Q}{ \dd h} \cdot \frac{\dd h}{\dd p} \\
	= & \left( \E[\Theta^2 (\Theta - \E[\Theta \mid \gamma \Theta + \sqrt{\gamma} G = h]) \mid \gamma \Theta + \sqrt{\gamma} G = h] - \right.\\
	&\left. 2\E[\Theta \mid \gamma \Theta + \sqrt{\gamma} G = h] \Var[\Theta \mid \gamma \Theta + \sqrt{\gamma} G = h] \right)  \cdot \Var[\Theta \mid \gamma \Theta + \sqrt{\gamma} G = h]^{-1/2}.
\end{align*}
Applying Cauchy–Schwarz inequality and the bounded support assumption, 
\begin{align*}
	& \left| \E[\Theta^2 (\Theta - \E[\Theta \mid \gamma \Theta + \sqrt{\gamma} G = h]) \mid \gamma \Theta + \sqrt{\gamma} G = h] \right| \leq M_{\Theta}^2 \cdot \Var[\Theta \mid \gamma \Theta + \sqrt{\gamma} G = h]^{1/2}, \\
	& \left| \E[\Theta \mid \gamma \Theta + \sqrt{\gamma} G = h] \Var[\Theta \mid \gamma \Theta + \sqrt{\gamma} G = h] \right|  \leq  M_{\Theta}^2 \cdot \Var[\Theta \mid \gamma \Theta + \sqrt{\gamma} G = h]^{1/2}. 
\end{align*}
Putting together the above analysis, we conclude that $\|\frac{\dd Q}{\dd p}\|_{\infty} \leq 3M_{\Theta}^2$, thus completing the proof of the lemma. 

\end{proof}
By \cref{eq:E2}, on the set $\cuE_{\beta, L, \delta, \ep, n}$, it holds that $\|\bD_{\gamma_{t}^{k}}(\hbm^k(\by(t),t))\|_F^2 / n \leq C_{\sconv}^{-1}\exp(-C_{\sconv} \beta^2)$. 
Invoking this result, triangle inequality, \cref{lemma:D-p}, and Cauchy-Schwartz inequality,
 we can conclude that for all $\bzeta^{k + 3}$ that lies on the line segment connecting 
 $\bp^{k + 3}$ and $\tilde{\bp}^{k + 3}$, it holds that
\begin{align*}
	&\frac{1}{n}\| \bD_{\gamma_t^{k + 3}}( \Psi_{\gamma_t^{k + 3}}( \bzeta^{k + 3}))\|_F^2 \\
	\leq & \frac{1}{n}\| \bD_{\gamma_t^{k + 3}}( \Psi_{\gamma_t^{k + 3}}( \hat{\bp}^{k + 3}(\by(t), t)))\|_F^2 + \frac{3M_{\Theta}^2}{n}\|\hat{\bp}^{k + 3}(\by(t), t) - \bzeta^{k + 3}\|_1 \\
	\leq & C_{\sconv}^{-1}\exp(-C_{\sconv} \beta^2) + M_{\Theta}^2 \cdot (2\beta^2 + 2\beta + 4 + 2C_{\sconv}^{-1}\exp(-C_{\sconv} \beta^2))\exp\left( -\frac{1}{4}C_{\sconv}\beta^2 \right) \\
	= & \Delta(\beta)\, ,
\end{align*}
where we recall that $\Delta(\beta)$ is defined in \cref{eq:functions}. Similarly, we can derive that for all $\bzeta^{k + 2}$ that is on the line segment connecting $\bp^{k + 2}$ and $\tilde{\bp}^{k + 2}$, 
$$\frac{1}{n}\| \bD_{\gamma_t^{k + 2}}( \Psi_{\gamma_t^{k + 2}}( \bzeta^{k + 2}))\|_F^2 \leq \Delta(\beta).$$ 
In addition, note that for all $\bzeta^{k + 2}, \bzeta^{k + 3}$ as above, it holds that 
\begin{align*}
	\max\left\{ \|\bD_{\gamma_t^{k + 2}}( \Psi_{\gamma_t^{k + 2}}( \bzeta^{k + 2}))\|_{\infty}, \, \|\bD_{\gamma_t^{k + 3}}( \Psi_{\gamma_t^{k + 2}}( \bzeta^{k + 2}))\|_{\infty}\right\} \leq M_{\Theta}.
\end{align*} 
Therefore, if we view the diagonal elements of matrices 
$\bD_{\gamma_t^{k + 2}}( \Psi_{\gamma_t^{k + 2}}( \bzeta^{k + 2}))$ and
 $\bD_{\gamma_t^{k + 3}}( \Psi_{\gamma_t^{k + 3}}( \bzeta^{k + 2}))$ as vectors, then they 
 belong to the set $\{\bx \in [0, M_{\Theta}]^n: \|\bx\|_2^2 / n \leq \Delta(\beta)\}$. 
Hence, recalling the definition of event  $\cuE_{\beta, n}^{(4)}$ in Eq.~\cref{eq:E4}, we see that for all $\bzeta^{k + 2} $ and $\bzeta^{k + 3}$,
 the following inequalities hold on $\cuE_{\beta, L, \delta, \eps, n} \cap \cuE_{\beta, n}^{(4)}$:
\begin{align*}
	\|\beta \bD(\bzeta^{k + 3}) \bX \bD(\bzeta^{k + 2})\|_{\op} &
	\leq  \frac{\beta^2}{n} \|\bD(\bzeta^{k + 3}) \btheta \btheta^{\sT} \bD(\bzeta^{k + 2})\|_{\op} + \beta \|\bD(\bzeta^{k + 3}) \bW \bD(\bzeta^{k + 2})\|_{\op} \\
	&\leq  \beta^2 M_{\Theta}^2 \Delta(\beta) + 4\beta C' M_{\Theta}^{5/3} \Delta(\beta)^{1/6}, \\
	\|b_t^k \bD(\bzeta^{k + 3}) \bD(\bzeta^{k + 1}) \|_{\op} &\leq  M_{\Theta}^2 \cdot C_{\sconv}^{-1}\exp(-C_{\sconv} \beta^2) .
\end{align*}
Putting together the above inequalities and \cref{eq:general-jacob-p}, we obtain that 
\begin{align}\label{eq:C3-1}
	& \frac{1}{\sqrt{n}} \|\bp^{k + 3} - \tilde{\bp}^{k + 3}\|_2 \nonumber \\
	\leq & \frac{\beta^2 M_{\Theta}^2 \Delta(\beta) + 4\beta C' M_{\Theta}^{5/3} \Delta(\beta)^{1/6}}{\sqrt{n}} \|\bp^{k + 2} - \tilde{\bp}^{k + 2}\|_2 + \nonumber\\
	& \frac{M_{\Theta}^2 \cdot C_{\sconv}^{-1}\exp(-C_{\sconv} \beta^2)}{\sqrt{n}} \|\bp^{k + 1} - \tilde{\bp}^{k + 1}\|_2 + \frac{M_{\Theta}}{\sqrt{n}}\|\by_1 - \by_2\|_2 \nonumber \\
	\leq &  \frac{\rho(\beta)}{\sqrt{n}} \|\bp^{k + 2} - \tilde{\bp}^{k + 2}\|_2 + \frac{M_{\Theta}^2 \cdot C_{\sconv}^{-1}\exp(-C_{\sconv} \beta^2)}{\sqrt{n}} \|\bp^{k + 1} - \tilde{\bp}^{k + 1}\|_2 + \frac{M_{\Theta}}{\sqrt{n}}\|\by_1 - \by_2\|_2  \\
	\overset{(d)}{\leq} & \frac{\Lipc(\beta) \cdot \rho(\beta) + \Lipc(\beta) \cdot M_{\Theta}^2 \cdot C_{\sconv}^{-1}e^{-C_{\sconv} \beta^2} + M_{\Theta}}{\sqrt{n}} \cdot \|\by_1 - \by_2\|_2   \nonumber \\
	{\leq} \nonumber &   \frac{\Lipc(\beta) \cdot F(\beta) + M_{\Theta}}{\sqrt{n}} \cdot \|\by_1 - \by_2\|_2 \nonumber \\
	\overset{(e)}{\leq} & \frac{\Lipc(\beta)}{\sqrt{n}} \cdot \|\by_1 - \by_2\|_2, \nonumber
\end{align}
where in step \emph{(d)} we used \cref{eq:diff-pB}  with $k_{\ast} \in\{k + 1, k + 2\}$, and in step \emph{(e)} we make use of the following facts: (1) $\Lipc(\beta) \geq 2M_{\Theta}$; (2) $F(\beta) \leq 1/ 2$ for all $\beta \geq \beta_0$. 

Recall that $2r(\beta) \cdot (\Lipc(\beta) + 1)\cdot (M_{\Theta}^2 + 1) \leq 2 e^{- C_{\sconv} \beta^2 / 4}$. 
Therefore, we can conclude that \cref{eq:diff-pB} holds for $k_{\ast} = k + 3$. 
In addition, for $\beta_0$ large enough clearly we have $M_{\Theta}^2 \cdot C_{\sconv}^{-1}\exp(-C_{\sconv} \beta^2) < \rho(\beta)$ holds for all $\beta \geq \beta_0$. As a result, we can deduce from \cref{eq:C3-1} that \cref{eq:49} holds for all desired $k$, thus completing the proof of \cref{lemma:support-main}.  

\end{proof}

\subsection{Proof of \cref{lemma:proof-of-eq}}\label{sec:proof-of-lemma:proof-of-eq}

Throughout  this proof, we work with $\bX, \btheta^0_1,\btheta^0_2$ 
with distribution $\prob_0$ defined in Section \ref{sec:ProofSymmetric}.
We will lighten notations by writing  $\btheta_i:=\btheta^0_i$.

We write the posterior distribution as
\begin{align}
\mu_{t}(\dd\btheta)& = \frac{1}{Z(t)} e^{H_t(\btheta)}\, \pi^{\otimes n}_{\Theta}(\de\btheta)\, ,\\
	H_t(\btheta) & :=  \frac{\beta}{2} \< \btheta, \bX \btheta \> 
	-\frac{\beta^2}{4n}\|\btheta\|_2^4 + \<\by(t),\btheta\>-\frac{t}{2}\|\btheta\|^2\, .
\end{align}
In this proof, we will never consider the joint distribution of these objects at two distinct
values of $t$. Hence, we can carry out derivations with
\begin{align}
\by(t) = t\, \btheta_*+\sqrt{t}\, \bz\, ,
\end{align}
for a fixed $\bz\sim\normal(0,\id_n)$.  
Further, we will write $\mu_t(F(\btheta_1,\btheta_2)):=\int F(\btheta_1,\btheta_2)\, 
\mu^{\otimes 2}_t(\dd\btheta)$.
Finally, we define
\begin{align}
U_t(\btheta)& = \frac{2}{n}\frac{\partial\phantom{t}}{\partial t} H_t(\btheta)\\
&= \frac{1}{n}\Big\{ 2\<\btheta,\btheta_*\>+\frac{1}{\sqrt{t}} \<\bz,\btheta\>-\|\btheta\|_2^2\Big\}\, .
\end{align}

Using Gaussian integration by parts and \cref{lemma:nishimori}, we have
\begin{align}
 \E\left[ \mu_t\big(U_t(\btheta)\big) \right]
	&=\frac{1}{n}\E\left[ \mu_t\big( \langle \btheta_1, \btheta_2 \rangle  \big)\right]\, ,
	 \label{eq:31} \\
	 \E\left[ \mu_t\Big( U_t(\btheta_1) \langle\btheta_1, \btheta_2 \rangle \Big) \right] & =
	 \frac{1}{n}\E\left[ \mu_t \big(\langle\btheta_1, \btheta_2 \rangle^2\big) \right]. \label{eq:32}
\end{align}

	Recall that ${\rm supp}(\pi_{\Theta})\subseteq [-M_{\Theta},M_{\Theta}]$, and therefore
	$|\langle \btheta_1, \btheta_2 \rangle / n| \leq M_{\Theta}^2$, which further yields
\begin{align}\label{eq:30}
 \left| \E\Big[  U_t(\btheta_1)
	\frac{1}{n} \langle \btheta_1, \btheta_2 \rangle  \Big] - 
	\E\Big[ U_t(\btheta) \Big]
	\E\Big[ \frac{1}{n} \langle \btheta_1, \btheta_2 \rangle\Big]\right| 
		\leq  M_{\Theta}^2 \cdot \E\left[  \Big|U_t(\btheta)  -
		 \E\big[ U_t(\btheta) \big]\Big|   \right].
\end{align}
Combining \cref{eq:30,eq:31,eq:32} gives
\begin{align}\label{eq:33}
	\frac{1}{n^2} \E\left[  (\langle\btheta_1, \btheta_2\rangle - \E[\langle \btheta_1, \btheta_2 \rangle ] )^2  \right] \leq M_{\Theta}^2 \cdot
	\E\left[  \Big|U_t(\btheta)  -
		 \E\big[ U_t(\btheta) \big]\Big|   \right]\, .
\end{align}
Next, we will prove that the right hand side of \cref{eq:33} is $o_n(1)$ for all $t > 0$, thus 
completing the proof of the lemma. 

We define the free energy density:
\begin{align*}
	\phi(t) := \frac{1}{n} \log  \Big\{\int e^{H_t(\btheta)}\, \pi^{\otimes n}_{\Theta}(\de\btheta)
	\Big\}\, .
\end{align*}
We can compute the first and second derivatives of $\phi$: 
\begin{align}
	& \frac{\partial\phi}{\partial\sqrt{t}}(t) = \sqrt{t}\mu_t \big(U_t(\btheta) \big), \label{eq:phip} \\
	& \frac{\partial^2\phi}{\partial(\sqrt{t})^2}(t) = nt \, {\rm Var}_{\mu_t} \big(U_t(\btheta)\big) + \frac{1}{n} \mu_t\Big(
	 2\langle \btheta_*, \btheta \rangle  - \|\btheta\|_2^2\Big). \nonumber
\end{align}
Therefore, defining $\psi(r) := \phi(r^2)$, we obtain
that $r\mapsto \bar\psi(r):=\psi(r)+3M_{\Theta}^2r^2/2$ is convex for $r\in(0,\infty)$.
Applying \cref{lemma:diff-of-derivative} to the functions $\bar\psi(r)$ and 
$\E\bar\psi(r)$, for $0 < \ep < r / 2$ we have
\begin{align}\label{eq:37}
	\E\left[ |\psi'(r) - \E[\psi'(r)]| \right] \leq &  \E[\psi'(r + \ep) - \psi'(r - \ep)] + \frac{3}{\ep} \sup\limits_{|r' - r| \leq \ep} \E\left[ |\psi(r') - \E[\psi(r')]| \right] + 
	6M_{\Theta}^2\ep. 
\end{align}
The next lemma proves that $ \sup_{|r' - r| \leq \ep} \E\left[ |\psi(r') - \E[\psi(r')]| \right]$
 is small. 
\begin{lemma}\label{lemma:uniform-conc}
	There exists a constant $C(t,\beta,\pi_{\Theta})> 0$ which is a function of $(t,\beta,  \pi_{\Theta})$ only,
	and is bounded compact intervals $[t_1,t_2]\subseteq (0,\infty)$
	 such that 
	\begin{align*}
		\E\left[ |\phi(t) - \E[\phi(t)]| \right] \leq C(t,\beta,\pi_{\Theta}) n^{-1/2}. 
	\end{align*} 
\end{lemma}
\begin{proof}
Letting $\bX = \beta\btheta_*\btheta^{\sT}/n+\bW$, consider the mapping 
\begin{align*}
	f: (\bW, \bz) \mapsto \phi(t).
\end{align*}
We denote by $W_{ij}$ the $(i,j)$-th entry of $\bW$. 
 The following upper bounds on the partial derivatives are straightforward:
\begin{align}\label{eq:34}
\begin{split}
	& \left| \frac{\partial}{ \partial \sqrt{n} W_{ij}} f(\bW, \bz) \right| \leq \beta M_{\Theta}^2n^{-3/2}, \\
	& \left| \frac{\partial}{\partial z_i} f(\bW, \bz)\right| \leq 2rM_{\Theta}n^{-1}.
\end{split} 
\end{align}
Hence by Gaussian concentration, we obtain 
\begin{align}\label{eq:35}
	\E_{\bW, \bz}\left[ (\phi(t) - \E_{\bW, \bz}[\phi(t)]) \right] \leq C_1n^{-1}.
\end{align}
Here and below, we denote by  $C_i$ constants that depend only on $\beta,M_{\Theta},r$ and 
and are bounded over compacts.

Finally, we show that $\E_{\bW, \bz}[\phi(t)]$, as a function of $\btheta$,
 concentrates around its expectation. This follows from  the estimate:
\begin{align*}
	\left| \frac{\partial}{\partial \theta_i} \E_{\bW, \bz}[\phi(t)] \right| \leq C_2 n^{-1}.
\end{align*}
Using Efron Stein's inequality, we get
\begin{align}\label{eq:36}
	\E_{\btheta}[(\E_{\bW, \bz}[\phi(t)] - \E_{\bW, \bz, \btheta}[\phi(t)])^2] \leq C_3 n^{-1}. 
\end{align}
The proof of \cref{lemma:uniform-conc} follows from  \cref{eq:35} and \cref{eq:36}. 
\end{proof}

We now conclude the proof of \cref{lemma:proof-of-eq}.  We have
\begin{align}
	\psi'(r) =\left.\frac{\partial\phi}{\partial\sqrt{t}}(t)\right|_{t=r^2} =
	 \left.\sqrt{t}\mu_t \big(U_t(\btheta) \big) \right|_{t=r^2}\, . 
\end{align}
Therefore letting $t_{\pm}: =\sqrt{r\pm \eps}$,  \cref{eq:31} and  \cref{eq:37} imply
\begin{align}
	\E\left[ |\mu_t \big(U_t(\btheta) \big) - \E[\mu_t \big(U_t(\btheta) \big)]| \right] \leq &  
	\frac{1}{n}\E[\mu_{t_+}(\<\btheta_1,\btheta_2\>) - 
	\mu_{t_-}(\<\btheta_1,\btheta_2\>)] + \frac{C_4}{\ep} n^{-1/2}+ 
	6M_{\Theta}^2\ep. 
\end{align}
Proposition \ref{prop:1} implies that the mapping $t \mapsto \gamma_{\ast}(\beta, t)$
 is locally  Lipschitz continuous on $(0, \infty)$. 
 Since 
$\lim_{n \to \infty}\E\left[ \mu_t( \langle\btheta_1, \btheta_2 \rangle / n ) \right] = 
\gamma_{\ast}(\beta, t)$, this yields
\begin{align}
\lim_{n\to\infty}\frac{1}{n}\E[\mu_{t_+}(\<\btheta_1,\btheta_2\>) - 
	\mu_{t_-}(\<\btheta_1,\btheta_2\>)] \le \delta(\ep)\, ,
	\end{align}
for some $\delta(\ep)\downarrow 0$ as $\ep\to 0$.
 
 Since $\ep$ is arbitrary, we obtain
\begin{align}
\lim_{n\to\infty}
 \E\left[ |\mu_t \big(U_t(\btheta) \big) - \E[\mu_t \big(U_t(\btheta) \big)]| \right] = 0\,.
 \end{align}

 The proof is completed by showing that, for all $0<t_1<t_2$
\begin{align}
\lim_{n\to\infty}\int_{t_1}^{t_2}
 \E\Big[ \Var_{\mu_t}(U_t(\btheta))\Big] \, \de t = 0\, .\label{eq:VarClaim}
\end{align}
We  notice that
\begin{align}
\frac{\de\phantom{t} }{\de t} \mu_t(U_t(\btheta)) = \frac{n}{2}\Var_{\mu_t}(U_t(\btheta))
-\frac{1}{2n t^{3/2}} \mu_t(\<\bz,\btheta\>)\, ,
\end{align}
and therefore
\begin{align}
\int_{t_1}^{t_2}
 \E\Big[ \Var_{\mu_t}(U_t(\btheta))\Big] \, \de t = \frac{2}{n}
\big\{\E\mu(U_{t_2}(\btheta))- \E\mu(U_{t_1}(\btheta))\big\}
+ \frac{1}{n^2}\int_{t_1}^{t_2} \frac{1}{t^{3/2}} \E \mu_t(\<\bz,\btheta\>)\, \de t\, .
\end{align}
On the other hand,  using Eq.~\eqref{eq:31} and $\|\btheta\|_{\infty}\le M_{\Theta}$, we get
\begin{align}
 \Big|\E \mu_t\big(U_t(\btheta)\big) \Big|& \le M_{\Theta}^2\, ,\\
 \Big|\E \mu_t(\<\bz,\btheta\>)\Big| & \le M_{\Theta}\sqrt{n}\E\|\bz\|_2\le  2M_{\Theta}n\, .
 \end{align}
Substituting above, we get
\begin{align}
\int_{t_1}^{t_2}
 \E\Big[ \Var_{\mu_t}(U_t(\btheta))\Big] \, \de t \leq \frac{4}{n}M_{\Theta}^2 
 + \frac{2}{n^2}M_{\Theta} \frac{t_2-t_1}{t_1^{3/2}}\, ,
  \end{align}
 which implies the claim \eqref{eq:VarClaim}.
 

\subsection{Proof of \cref{lemma:lower-bound-var}}
\label{sec:proof-lemma:lower-bound-var}

\begin{proof}[Proof of \cref{lemma:lower-bound-var}]
	We denote by $\pi_{x, t}$ the posterior distribution of $\Theta$ upon observing $t \Theta + \sqrt{t} G = x$, where $(\Theta, G) \sim \pi_{\Theta} \otimes \normal(0, 1)$. 
	It is not hard to see that 
	\begin{align*}
		\pi_{x, t} (\theta) = e^{-V_{x, t}(\theta)}, \qquad V_{x, t}(\theta) = -U(\theta) + \langle x, \theta \rangle - \frac{t\|\theta\|_2^2}{2} - \log \int e^{-U(\theta) + \langle x, \theta \rangle - t \|\theta\|_2^2 / 2} \dd \pi_{x, t}(\theta). 
	\end{align*}
	By assumption, we know that $\|V_{x, t}''\|_{\infty} \leq C + t$. We denote by $\mathcal{M}_{x, t}$ the median of $\pi_{x, t}$. Therefore, for any $a \geq 0$
	\begin{align*}
		\Var[\Theta \mid t \Theta + \sqrt{t} G = x] = & \frac{1}{2} \int (\theta - \theta')^2 \dd \pi_{x, t}(\theta) \dd \pi_{x, t}(\theta') \\
		\geq & \frac{1}{2} \int (\theta - \theta')^2 \mathbbm{1}_{\theta' \leq \mathcal{M}_{x, t}} \mathbbm{1}_{\theta \geq \mathcal{M}_{x, t}} \dd \pi_{x, t}(\theta) \dd \pi_{x, t}(\theta') \\
		\geq & \frac{1}{2} \int (\theta - \mathcal{M}_{x, t})^2 \mathbbm{1}_{\theta \geq \mathcal{M}_{x, t}}  \dd \pi_{x, t}(\theta) \\
		\geq & \frac{a^2}{2} \PP_{\theta \sim \pi_{x, t}}(\theta > \mathcal{M}_{x, t} + a) = \frac{a^2}{2} \left( \frac{1}{2} - \PP_{\theta \sim \pi_{x, t}}(\theta \in [\mathcal{M}_{x, t}, \mathcal{M}_{x, t} + a]) \right). 
	\end{align*}
	We next upper bound $\PP_{\theta \sim \pi_{x, t}}(\theta \in [\mathcal{M}_{x, t}, \mathcal{M}_{x, t} + a])$. Note that for any $\bar a \geq 0$, 
	\begin{align*}
		\PP_{\theta \sim \pi_{x, t}}(\theta \in [\mathcal{M}_{x, t}, \mathcal{M}_{x, t} + a]) = \int_{\mathcal{M}_{x, t}}^{\mathcal{M}_{x, t} + a} e^{-V(\theta)} \dd \theta \leq  \frac{\int_{\mathcal{M}_{x, t}}^{\mathcal{M}_{x, t} + a} e^{-V(\theta)} \dd \theta}{\int_{\mathcal{M}_{x, t}}^{\mathcal{M}_{x, t} + \bar a} e^{-V(\theta)} \dd \theta} = \frac{\int_{\mathcal{M}_{x, t}}^{\mathcal{M}_{x, t} + a} e^{-V(\theta) + V(\mathcal{M}_{x, t})} \dd \theta}{\int_{\mathcal{M}_{x, t}}^{\mathcal{M}_{x, t} + \bar a} e^{-V(\theta) + V(\mathcal{M}_{x, t})} \dd \theta}. 
	\end{align*}
	Since $\|V''\|_{\infty} \leq C + t$, we conclude that the numerator in the last term above is upper bounded by
	\begin{align*}
		e^{(C + t) a^2 / 2}\int_{\mathcal{M}_{x, t}}^{\mathcal{M}_{x, t} + a} e^{-V'(\mathcal{M}_{x, t})(\theta - \mathcal{M}_{x, t})} \dd \theta = \frac{e^{(C + t) a^2 / 2}(1 - e^{-V'(\mathcal{M}_{x, t}) a})}{V'(\mathcal{M}_{x, t})}, 
	\end{align*} 
	and the denominator is lower bounded by 
	\begin{align*}
		e^{-(C + t) \bar a^2 / 2} \int_{\mathcal{M}_{x, t}}^{\mathcal{M}_{x, t} + \bar a} e^{-V'(\mathcal{M}_{x, t})(\theta - \mathcal{M}_{x, t})} \dd \theta = \frac{e^{-(C + t) \bar a^2 / 2} (1 - e^{-V'(\mathcal{M}_{x, t}) \bar a})}{V'(\mathcal{M}_{x, t})}. 
	\end{align*}
	If $V'(\mathcal{M}_{x, t}) = 0$, then 
	\begin{align*}
		\PP_{\theta \sim \pi_{x, t}}(\theta \in [\mathcal{M}_{x, t}, \mathcal{M}_{x, t} + a]) \leq \frac{a e^{(C + t) a^2 / 2}}{\bar a e^{-(C + t) \bar a^2 / 2}}. 
	\end{align*}
	Setting $a = 1 / (100\sqrt{C + t})$ and $\bar a = 1 / (10 \sqrt{C + t})$, we conclude that $\PP_{\theta \sim \pi_{x, t}}(\theta \in [\mathcal{M}_{x, t}, \mathcal{M}_{x, t} + a]) \leq e^{1/100} / 10$, hence $\Var[\Theta \mid t \Theta + \sqrt{t} G = x] \geq (10^5(C+t))^{-1}$. 
	
	If $V'(\mathcal{M}_{x, t}) < 0$, then by \cref{lemma:a-ratio}, it holds that 
	\begin{align*}
		\PP_{\theta \sim \pi_{x, t}}(\theta \in [\mathcal{M}_{x, t}, \mathcal{M}_{x, t} + a]) \leq e^{(C + t)(a^2 + \bar a^2) / 2} \cdot \frac{1 - e^{-V'(\mathcal{M}_{x, t}) a}}{1 - e^{-V'(\mathcal{M}_{x, t}) \bar a}} \leq \frac{a e^{(C + t)(a^2 + \bar a^2) / 2}}{\bar a}. 
	\end{align*}
	Once again taking $a = 1 / (100\sqrt{C + t})$ and $\bar a = 1 / (10 \sqrt{C + t})$, we conclude that $\Var[\Theta \mid t \Theta + \sqrt{t} G = x] \geq (10^5(C+t))^{-1}$. 
	
	If $V'(\mathcal{M}_{x, t}) > 0$, then we switch to another lower bound: 
	\begin{align*}
		\Var[\Theta \mid t \Theta + \sqrt{t} G = x] = & \frac{1}{2} \int (\theta - \theta')^2 \dd \pi_{x, t}(\theta) \dd \pi_{x, t}(\theta') \\
		\geq & \frac{1}{2} \int (\theta - \theta')^2 \mathbbm{1}_{\theta \leq \mathcal{M}_{x, t}} \mathbbm{1}_{\theta' \geq \mathcal{M}_{x, t}} \dd \pi_{x, t}(\theta) \dd \pi_{x, t}(\theta') \\
		\geq & \frac{1}{2} \int (\theta - \mathcal{M}_{x, t})^2 \mathbbm{1}_{\theta \leq \mathcal{M}_{x, t}}  \dd \pi_{x, t}(\theta) \\
		\geq & \frac{a^2}{2} \PP_{\theta \sim \pi_{x, t}}(\theta < \mathcal{M}_{x, t} - a) = \frac{a^2}{2} \left( \frac{1}{2} - \PP_{\theta \sim \pi_{x, t}}(\theta \in [\mathcal{M}_{x, t} - a, \mathcal{M}_{x, t}]) \right). 
	\end{align*}
	Note that for any $\bar a \geq 0$, 
	\begin{align*}
		\PP_{\theta \sim \pi_{x, t}}(\theta \in [\mathcal{M}_{x, t} - a, \mathcal{M}_{x, t}]) = \int_{\mathcal{M}_{x, t} - a}^{\mathcal{M}_{x, t}} e^{-V(\theta)} \dd \theta \leq  \frac{\int_{\mathcal{M}_{x, t} - a}^{\mathcal{M}_{x, t}} e^{-V(\theta)} \dd \theta}{\int_{\mathcal{M}_{x, t} - \bar a}^{\mathcal{M}_{x, t}} e^{-V(\theta)} \dd \theta} = \frac{\int_{\mathcal{M}_{x, t} - a}^{\mathcal{M}_{x, t}} e^{-V(\theta) + V(\mathcal{M}_{x, t})} \dd \theta}{\int_{\mathcal{M}_{x, t} - \bar a}^{\mathcal{M}_{x, t}} e^{-V(\theta) + V(\mathcal{M}_{x, t})} \dd \theta}. 
	\end{align*}
	The numerator above is upper bounded by 
	\begin{align*}
		e^{(C + t) a^2 / 2}\int_{\mathcal{M}_{x, t} - a}^{\mathcal{M}_{x, t}} e^{-V'(\mathcal{M}_{x, t})(\theta - \mathcal{M}_{x, t})} \dd \theta = \frac{e^{(C + t) a^2 / 2}(-1 + e^{V'(\mathcal{M}_{x, t}) a})}{V'(\mathcal{M}_{x, t})}.  
	\end{align*}
	The denominator on the other hand is lower bounded by 
	\begin{align*}
		e^{-(C + t) \bar a^2 / 2} \int_{\mathcal{M}_{x, t} - \bar a}^{\mathcal{M}_{x, t}} e^{-V'(\mathcal{M}_{x, t})(\theta - \mathcal{M}_{x, t})} \dd \theta = \frac{e^{-(C + t) \bar a^2 / 2} (-1 + e^{V'(\mathcal{M}_{x, t}) \bar a})}{V'(\mathcal{M}_{x, t})}. 
	\end{align*}
	Putting together the above bounds and \cref{lemma:a-ratio}, we conclude that for $V'(\mathcal{M}_{x, t}) > 0$, it holds that
	\begin{align*}
		\PP_{\theta \sim \pi_{x, t}}(\theta \in [\mathcal{M}_{x, t} - a, \mathcal{M}_{x, t}]) \leq e^{(C + t)(a^2 + \bar a^2) / 2} \cdot \frac{1 - e^{V'(\mathcal{M}_{x, t}) a}}{1 - e^{V'(\mathcal{M}_{x, t}) \bar a}} \leq \frac{a e^{(C + t)(a^2 + \bar a^2) / 2}}{\bar a}.
	\end{align*}
	The proof is complete by setting $a = 1 / (100\sqrt{C + t})$ and $\bar a = 1 / (10 \sqrt{C + t})$. In this case, $\Var[\Theta \mid t \Theta + \sqrt{t} G = x] \geq (10^5(C+t))^{-1}$.

\end{proof}

\end{appendices}
	
\end{document}